%% file: main.tex
\documentclass[final]{siamltex}
\usepackage{latexsym, graphicx, epsfig, amsmath, amsfonts,amssymb,bm,mathrsfs}
\usepackage{subfigure}
\usepackage{color}
\newtheorem{remark}{Remark}[section]
\usepackage{url}

\def\c{\mathbf{c}}

\def\be{\begin{equation}}
\def\ee{\end{equation}}

\def\x{\mathbf{x}}

\def\c{\mathbf{c}}

\def\R{{\mathbb R}}

\def\u{{\bf u}}
\def\tu{\widetilde{\bf u}}
\def\hu{\widehat{\bf u}}
\def\tH{\widetilde{H}}

\def\c{\mathbf{c}}

\newcommand{\figref}[1]{Fig.~\ref{#1}}

\title{Structure-preserving Method for Reconstructing Unknown Hamiltonian Systems from Trajectory Data}

\author{Kailiang Wu\and Tong Qin \and Dongbin
       Xiu\thanks{Department of Mathematics,
		The Ohio State University, Columbus, OH 43210, USA.
		{\tt wu.3423@osu.edu, qin.428@osu.edu, xiu.16@osu.edu.}
		Funding: This work was partially supported by AFOSR FA9550-18-1-0102.}
}

\begin{document}
\maketitle
\begin{abstract}
We present a numerical approach for approximating unknown Hamiltonian systems using observational data. A distinct feature of the proposed method is that it is structure-preserving, in the sense that it enforces the conservation of the reconstructed Hamiltonian. This is achieved by directly approximating the underlying unknown Hamiltonian, rather than the 
right-hand-side of the governing equations. We present the technical details of the proposed algorithm and its error estimate in a special case, along with a practical de-noising 
procedure to cope with noisy data.  A set of numerical examples are presented to demonstrate the structure-preserving property 
and effectiveness of the algorithm.
\end{abstract}
\begin{keywords}
data-driven discovery, Hamiltonian systems, equation approximation, structure-preserving method
\end{keywords}

\input Introduction

\input SetupNew
\input Method

\input Examples

\input Conclusion
\input Appendix



\bibliographystyle{siamplain}
\bibliography{neural,LearningEqs,Hamiltonian,identification,least-squares}

\end{document}

%% file: Introduction.tex
\section{Introduction} \label{sec:intro}

Data-driven discovery of physical laws has received
an increasing amount of attention recently.
While
earlier attempts such as
\cite{bongard2007automated,schmidt2009distilling} used
symbolic regression to select the proper physical laws and
determine the underlying dynamical systems, more recent efforts tend
to treat the problem as an approximation problem. In this approach,
%
the sought-after governing equation is treated as
an unknown target function relating the data of the state variables to their temporal
derivatives. Methods along this line of approach usually seek to
exactly recover the equations by using certain sparse approximation 
techniques (e.g., \cite{tibshirani1996regression}) from a large set of dictionaries; see, for example, \cite{brunton2016discovering}.
Many studies have been conducted to effectively deal with noise in data
\cite{brunton2016discovering, schaeffer2017sparse}, corruptions in data
\cite{tran2017exact}, limited data \cite{schaeffer2017extracting}, 
partial differential equations \cite{rudy2017data,
	schaeffer2017learning}, etc, and in conjunction with other
      methods such as
model selection approach
\cite{Mangan20170009}, Koopman theory \cite{brunton2017chaos}, and
Gaussian process regression \cite{raissi2017machine}, to name a few.
Standard approximations using polynomials without seeking exact
recovery can also be effective (cf.~\cite{WuXiu_JCPEQ18}).
%
More recently, there is a surge of work that tackles the problem using machine
learning methods, particularly via neural networks
\cite{raissi2017physics1,raissi2017physics2}, to systems involving
ordinary differential equations (ODEs) \cite{E2017,raissi2018multistep,rudy2018deep,qin2018data, qin2019neural} and partial differential equations (PDEs)
\cite{MardtPWN_Nature18, long2017pde, KhooLuYing_2018, han2018solving,raissi2018deep,long2018pde, wu2019data}. 
Neural network structures such as 
residual network
(ResNet) were shown to be highly suitable for this type of problems
\cite{qin2018data}.

Hamiltonian systems are an important class of governing equations in
science and engineering. One of the most important properties of
Hamiltonian systems is conservation of Hamiltonian, 
usually a nonlinear function of state variables, along
trajectories. Research efforts have been devoted to estimating
Hamiltonian function of a given system from measurements;  
cf.~\cite{Shabani2011,wiebe2014hamiltonian, fujita2018construction,
  Bairey2019, li2018density}.
However, few studies exist to reconstruct an unknown Hamiltonian dynamical systems from trajectory data of state variables.

In this paper, we present a numerical approach to reconstruct an unknown
Hamiltonian system from its trajectory data. The focus of this paper
is on the conservation of the reconstructed
Hamiltonian along the solution trajectories.
The current method is an extension of the method proposed in
\cite{WuXiu_JCPEQ18}, which seeks accurate approximation of unknown
governing equations using orthogonal polynomials.
However, instead of approximating the governing equations directly, as in
\cite{WuXiu_JCPEQ18} and most of other existing studies, our current
method seeks to 
approximate the unknown Hamiltonian first and then reconstruct the
approximate governing
equations using the approximate Hamiltonian.
The approximation of the unknown Hamiltonian is conducted
using orthogonal polynomials and with controllable numerical errors.
Since in most practical situations, the
Hamiltonian takes the form of smooth functions, polynomial
approximation can achieve high order accuracy with modest degree polynomials.
The resulting approximate governing equations, which are derived from
the reconstructed Hamiltonian, can then automatically satisfy
the conservation of the reconstructed Hamiltonian, which is an accurate
approximation of the true Hamiltonian.
This structure preserving (SP) property---the conservation of
Hamiltonians along trajectories---is a
distinctly new feature of our present method, not found in most of the
existing studies. Along with a detailed exposition of the algorithm, we
also provide an error estimate of the method in a special case and use a set of
numerical examples to demonstrate the properties of the method. 



%% file: SetupNew.tex
\section{Preliminaries} \label{sec:setup}

In this section, we introduce some basics about Hamiltonian systems
and the setup of our data-driven discovery of Hamiltonian systems.

\subsection{Hamiltonian Systems}

Let us consider a Hamiltonian system
\begin{equation}\label{eq:HamPQ}
\begin{aligned}
& \frac{d {\bm p}}{ d t } = - \nabla_{\bm q} H({\bm p},{\bm q}), \\
& \frac{d {\bm q}}{ d t } = \nabla_{\bm p} H({\bm p},{\bm q}),
\end{aligned}
\end{equation} 
where $\bm p$ and $\bm q$ are column vectors in $\mathbb R^d$, and
$H(\bm p,\bm q)$ is a  continuously differentiable scalar  function
called Hamiltonian. The Hamiltonian often represents the total energy of the system. 
It is not unique and is defined up to an arbitrary constant.

Let ${\bf u} := ( \bm p^\top, \bm q^\top )^\top$ be the state variable
vector. The Hamiltonian system \eqref{eq:HamPQ} can be equivalently written as 
\begin{equation}\label{eq:Ham}
\frac{d {\bf u}}{ d t } =  {\bf J}^{-1} \nabla_{\bf u} H({\bf u}), \\
\end{equation} 
where $\nabla_{\bf u}$ stands for full gradient and 
the matrix $\bf J$ takes the form
\begin{equation}\label{eq:J}
{\bf J} = \begin{pmatrix}
{\bf 0}_d & {\bf I}_d\\
-{\bf I}_d & {\bf 0}_d
\end{pmatrix},
\end{equation}
with ${\bf I}_d$ and ${\bf 0}_d$ being identity matrix and zero matrix
of size $d\times d$, respectively. Hereafter, we will  use 
$\nabla$ in place of $\nabla_{\bf u}$, unless confusion arises otherwise.

The Hamiltonian system \eqref{eq:Ham} is an autonomous system
and conserves the Hamiltonian along the integral curves (c.f., \cite{marsden2013introduction}). That is, 
the solution of the Hamiltonian system \eqref{eq:Ham} satisfies 
$$
H ({\bf u} ( t; {\bf u}_0 )) = H ( {\bf u}_0 ),\quad~\forall t\ge 0,~
\forall {\bf u}_0,
$$
where ${\bf u}_0$ is the initial state of the system at $t=0$, and
${\bf u} ( t; {\bf u}_0 )$ stands for the solution ${\bf u}$ at time
$t$ with an initial state ${\bf u}_0$.

\subsection{Data and Problem Setup}

We assume that the equations of the Hamiltonian system
\eqref{eq:HamPQ}, or \eqref{eq:Ham}, are not known. Our data-driven
approach for approximating the
unknown equations requires the availability of a set of data pairings
between the solution states and their corresponding time derivatives,
in the form of
  \begin{equation}
\label{pairing}
	\left\{ {\bf x}_k, \dot{\bf x}_k \right\},\qquad k=1,\ldots, K,
\end{equation}
where $\dot{\x}$ denotes the time derivative of $\x$, and $K$ is the total number of data pairs.


\subsubsection{Direct Data Collection}

Let $D\in\R^{2d}$ be a bounded domain.
It is the domain-of-interest,
inside which we seek to construct an accurate approximation to the
unknown governing equations \eqref{eq:Ham}.
Our data set \eqref{pairing} shall be collected in the domain $D$.

When time derivatives of the state variables are readily available, either
directly measured from experiments or computed via certain numerical
simulation techniques, the data collection procedure is
straightforward.
Let $M\geq 1$ be the number of solution trajectories, originated from $\u_0^{(1)},\dots,
\u_0^{(M)}$ initial states.
Let $0=t_0^{(m)}<t_1^{(m)}<\cdots <t_{J_m}^{(m)}$ be
a sequence of time instances on the $m$-th trajectory, for $m=1,\dots, M$.  We assume that the
state variables data and their derivative data
are available on these time instances, i.e., for $0\leq  j \le J_m $ and $ 1\leq m \leq M$,
\begin{align} 
\label{x0}
{\bf x}^{(m)}(t_{j}^{(m)})& = {\bf u}(t_{j}^{(m)};\u_0^{(m)}) + {\bm \epsilon}_{j}^{(m)},\\
\label{dx0}
{{\bf \dot{x}}^{(m)}(t_j^{(m)})} &= \frac{d}{dt}{\bf u}(t_{j}^{(m)};\u_0^{(m)}) + {\bm \tau}_{j}^{(m)},
\end{align}
where
${\bm \epsilon}_j^{(m)}$ and ${\bm \tau}_j^{(m)}$ are
errors/noises in the data for the state variables and their time
derivatives, respectively.

Once all the data pairs are collected, they are grouped in the set
\eqref{pairing}, where we omit the subscripts and superscripts for
notational convenience. This is because our method for approximating the governing
equation using the data does not utilize the trajectory or time instance
information associated with each data pair.



\subsubsection{Time Derivatives Approximation and
  De-noising}\label{sec:der}

In many practical situations, time derivative data are not
available. Consequently, one only possesses trajectory data for the
state variables.
In this case, it is necessary to estimate the time derivatives of the
state variables via a numerical procedure.

Again, let $M$ be the number of trajectories where {\em only the state
  variable data} are
available. For notational convenience we let $0=t_0<t_1<\cdots <t_{J}$ be
a sequence of the same time instances on all trajectories, where the state
variable data
\begin{equation} 
{\bf x}^{(m)}(t_{j}) = {\bf u}(t_{j};\u_0^{(m)}) + {\bm
  \epsilon}_{j}^{(m)}, \qquad j=0,\dots, J,
\end{equation}
are available. Again, ${\bm\epsilon}_{j}^{(m)}$ stands for errors/noises
in the data. To numerically estimate the time derivatives, it is
necessary to require $J\geq 1$, i.e., there need to be at least two data entries
of the state variables along each trajectory in order to estimate the time derivatives.


For noiseless data, i.e. ${\bm
	\epsilon}_{j}^{(m)}=\mathbf{0}$, time derivatives can be
      computed by straightforward 
numerical differentiation.  For example,
for equally spaced time instances with uniform step-size $\Delta t$, a second-order finite difference 
\begin{equation}\label{eq:finitediff}
{\bf \dot x}^{(m)}(t_{j})
= \frac{ {\bf x}^{(m)}(t_{j+1}) -  {\bf x}^{(m)}(t_{j-1})   }{2 \Delta
  t}, \qquad 1\le j\le J-1,
\end{equation}
with proper one-sided second-order finite difference at the end points
$j=0$ and $j=J$. This requires at least three data entries on each
trajectory, i.e., $J\geq 2$, and induces errors of $O(\Delta t^2)$. Higher order approximations requires more
data points on each trajectory.

For noisy state variable data with ${\bm
	\epsilon}_{j}^{(m)}\neq\mathbf{0}$, direct numerical differentiation is less robust,
as the errors in estimating  ${\bf \dot x}^{(m)}(t_{j}) $ would scale as 
$\sim {\bm \epsilon}_{j}^{(m)}/\Delta t$.  
Several techniques have been developed for numerical differentiation of
noisy data. See, for example,
\cite{knowles2014methods,wagner2015regularised,doi:10.1137/0708026,chartrand2011numerical,knowles1995variational}.  
In this paper, we employ a straightforward de-noising approach, which
has been shown to be effective for equation recovery (\cite{WuXiu_JCPEQ18}).
We first construct a least squares polynomial approximation of the
trajectory using the available data on $\bf{x}$, and then analytically differentiate the
least squares fitted polynomial to obtain an estimate of the time derivatives.
More specifically, for each trajectory,
the least squares polynomial approximation is to find a polynomial vector
${\bm {\mathcal P}}_m \in [\mathbb P^1_Q]^{2d}$ such that
\begin{equation}\label{eq:LS}
 {\bm {\mathcal P}}_m = \mathop{\rm argmin}\limits_{ {\bm {\mathcal P}} \in [\mathbb P_Q^1]^{2d} }  \sum_{j=0}^{Q} \left\| 
{\bm {\mathcal P}}(t_j) - {\bf x}^{(m)}(t_{j})
\right\|_2^2, \qquad m=1,\dots,M,
\end{equation} 
where $\| \cdot \|_2$ denotes vector 2-norm and $\mathbb P_Q^1$
denotes the space of one-dimensional polynomials of degree at most $Q$, with $1\le Q \le J$. 
Once the least squares fitting problem is solved, the time derivatives
can be approximated by differentiating the polynomials, i.e.,
\begin{equation}\label{ederivative}
{\bf \dot x}^{(m)}(t_{j}) := \frac{d }{dt} {\bm {\mathcal P}}_m(t_j) \approx \frac{d}{dt} {\bf { u}}(t_j;{\bf u}_0^{(m)}). 
\end{equation}

This approach also provides a filter to de-noise the noisy trajectory
data. For noisy data, we advocate the use of the filtered trajectory data to
replace the original noisy data, i.e.,
\begin{equation}\label{filter}
{\bf x}^{(m)}(t_{j}) \leftarrow {\bm {\mathcal P}}_m(t_j).
\end{equation}
Our numerical experiments indicate that this filtering procedure can
improve the learning accuracy for
noisy data. This is similar to the results from \cite{WuXiu_JCPEQ18}.




\begin{remark}
If the true trajectories ${\bf u}(t;\u_0^{(m)})$ are non-smooth, estimating time derivatives using global 
approximation \eqref{eq:LS} might not be sufficiently accurate. In this case,
piecewise approximation should be considered.
\end{remark}


%

%% file: Method.tex
\section{The Main Method} \label{sec:method}


With the data pairs \eqref{pairing}, our goal is now to accurately
approximate the unknown Hamiltonian system \eqref{eq:Ham}.
Let ${\bf f} := \nabla H $, which is the unknown right-hand-side of \eqref{eq:Ham}.  We seek an
accurate approximation $\widetilde {\bf f} \approx {\bf f}$ such that
\begin{equation}\label{appsystem}
\frac{ d {\bf u} }{ dt } = {\bf J}^{-1}  \widetilde {\bf f} ( {\bf u} )
\end{equation}
is an accurate approximation of the true system \eqref{eq:Ham}. Our
key goal is to ensure the approximate system is also Hamiltonian, in
the sense that $\widetilde {\bf f} = \nabla \widetilde H$,  where 
$\widetilde H$ becomes an approximation to the true (and unknown)
Hamiltonian. The existing methods for equation recovery seek to
approximate the right-hand-side of the true system directly and
therefore do not enforce the conservation of Hamiltonian.
%
%

\subsection{Algorithm}

To preserve the Hamiltonian, we propose to directly approximate the
unknown Hamiltonian first and then derive the approximate governing
equations from the approximate Hamiltonian.
%

Let us assume the unknown Hamiltonian $H \in {\mathbb H^1_\omega}(D)$,
which is a weighted Sobolev space on domain $D \subset \mathbb R^{2d}$ equipped with inner product 
$$
(g,h)_{ {\mathbb H^1_\omega} } = \int_D \big( g h + \nabla g \cdot \nabla h \big)  {\rm d} \omega (\x),
$$
where $\omega(\x)$ is a (probability) measure defined on $D$. 

Let $\mathbb W \subset {\mathbb H^1_\omega}(D)$ be a  finite
dimensional subspace. We then define its associated gradient function space as 
\begin{equation}\label{eq:gradspaceV}
\mathbb V := \left\{ \nabla h: h \in {\mathbb W}  \right\}, \qquad N=\dim \mathbb V \ge 1.
\end{equation}
Let $\{ {\bm \psi}_j({\bf x}) \}_{j=1}^N$ be a  basis for $\mathbb
V$. Then, for each $j=1,\dots, N$, there exists a function 
$\phi_j \in \mathbb W$ such that 
\begin{equation}\label{psiphi}
{\bm \psi}_j = \nabla \phi_j, \qquad 1\le j \le N.
\end{equation}



We then seek $\widetilde H \in \mathbb W$ as an approximation to the
true Hamiltonian $H$ and $\widetilde {\bf f} = \nabla \widetilde H \in
\mathbb V$ as an approximation to $\nabla H$. Assume that $N < K$, i.e., the dimension of the linear subspace is smaller than
the total number of available data pairs \eqref{pairing},   we then
define the following  least squares problem 
\begin{equation}\label{LSQ}
 \nabla \widetilde H = \mathop{\rm argmin}\limits_{ \nabla h \in
   \mathbb V } \sum_{k=1}^K 
\left\| {\bf J} \dot{\x}_k - \nabla h \left( \x_k  \right)   \right\|^2_2, 
\end{equation}
where $\| \cdot \|_2$ denotes vector 2-norm. 

With the basis \eqref{psiphi}, $\nabla \widetilde H $ can be expressed as
\begin{equation}\label{expan1}
\nabla \widetilde H ({\bf x}) = \sum_{j =1}^N c_j {\bm \psi}_j ( {\bf x} ) = 
\sum_{j =1}^N c_j \nabla \phi_j ( {\bf x} ). 
\end{equation}
This provides a class of approximate Hamiltonians which differ only in an additive constant $C$:  
\begin{equation}\label{expan2}
\widetilde H ( {\bf x} ) = C +  \sum_{j =1}^N c_j  \phi_j ( {\bf x} ) =: C + \widetilde H_0 ({\bf x}),
\end{equation}
where the constant $C$ can be arbitrarily chosen and does not affect the resulting approximate Hamiltonian system \eqref{appsystem}. In particular, when taking $C=0$, we use 
	the notation $\widetilde H_0$ for $\widetilde H $. 
The problem \eqref{LSQ} is then
equivalent to the following problem for the unknown coefficients
$
{\bf c} = (c_1,\dots, c_N)^\top,
$
\begin{equation}\label{LSQc}
\mathop{\rm min}\limits_{ \c \in \mathbb R^{N} } \| {\bf A} {\bf c} - {\bf b} \|_2,
\end{equation}
where 
\begin{equation}\label{eq:DefA}
{\bf A} = ( a_{ij} )_{1\le i, j \le N}, \qquad 
{\bf b} = (b_1,\dots, b_{N})^\top, 
\end{equation}
with 
\begin{align*}
&a_{ij} = \frac{1}{K } 
\sum_{k=1}^K \Big( \nabla \phi_{i} \left( \x_k  \right)
\cdot \nabla \phi_{j}  \big( \x_k \big) \Big), \quad 1\le i,j \le N,  \\
&b_i = \frac{1}{K} \sum_{k=1}^K \Big( \left( {\bf J} \,\dot{\x}_k \right) 
\cdot \nabla \phi_{i} \big( {\x}_k  \big) \Big), \quad 1\le i \le N.
\end{align*}

This is an over-determined system of equations and can be readily solved.
Upon solving this least squares type problem, we obtain $\widetilde H$ and
subsequently $\widetilde {\bf f} = \nabla \widetilde H$, which gives
us the approximate system of equations \eqref{appsystem}. It is
trivial to see that the system preserves the approximate Hamiltonian
$\widetilde H$ in the following sense.
\begin{theorem}\label{thm:Hconst}
Let $\widetilde {\bf u}(t; {\bf u}_0  ) $ be the solution  of the
system \eqref{appsystem} with initial state ${\bf u}_0$, then,
\begin{equation}\label{tildeHconse}
\widetilde H ( \widetilde {\bf u} ( t; {\bf u}_0 )) = \widetilde H ( {\bf u}_0 ),\qquad\forall t\ge 0,\quad \forall {\bf u}_0.
\end{equation}
\end{theorem}


\subsection{Analysis}

We now present error analysis for the proposed algorithm in a special case. Our analysis is based on a few basic results from \cite{Cohen2013} for
least squares polynomial approximations, which requires the following assumptions on the basis functions and the data.

\subsubsection{Assumptions}
The basis functions $\{ \phi_j
\}_{j=1}^N$ are assumed to be orthonormal in the following sense 
$$
\int_D \nabla \phi_i (\x) \cdot \nabla \phi_j (\x ) d \omega (\x) = \delta_{ij} . 
$$
Note that this assumption is only needed for the theoretical
analysis. The practical
computation of $\nabla \widetilde H$ can be conducted by using any basis
of $\mathbb V$, for the solution $\nabla \widetilde H$ does not depend
on the basis. (Also, any non-orthogonal basis can be orthogonalized
via Gram-Schmidt procedure.)
 We remark that 
the choice of basis affects the stability of the least squares problem \eqref{LSQc}. 
The actual
computation of $\nabla \widetilde H$ can be made using any known basis of $\mathbb W$, since the solution 
to the problem \eqref{LSQ} is independent of the chosen basis. 
Thus, the error estimates in Section \ref{sec:bound} also hold for any other bases of $\mathbb W$.

We assume that data for the state variable $\x_k$, $k=1,2,\dots,K$, 
are i.i.d.\ drawn from a probability measure $\omega (\x)$ on
$D$. 
This is a standard assumption, made mostly to
  facilitate theoretical analysis. See, for example, \cite{Cohen2013}.

\subsubsection{Stability}

The following stability result holds for the least squares problem \eqref{LSQc}.

\begin{lemma}\label{thm:Prob}
Consider the problem \eqref{LSQc}, it holds that, for  $0<\delta < 1$,
\begin{equation}\label{Pr1}
{\rm Prob} \big\{ \| {\bf A} - {\bf I} \| > \delta \big\} \le 2 N\exp\left( - \frac{  \beta_\delta K}{{\mathscr K}_N}\right),
\end{equation}
where $\beta_\delta :=  (1+\delta) \log (1+\delta)- \delta>0$, and 
\begin{equation}\label{eq:KN}
{\mathscr K}_N := 
\sup_{\x \in D} \sum_{j =1}^N \left\|  \nabla \phi_j (\x) \right\|_2^2.
\end{equation}
\end{lemma}

\begin{proof}
The proof is a direct extension of the proof of Theorem 1 in
\cite{Cohen2013} (see also \cite{Cohen2019} for a correction).
\end{proof}

\begin{remark}
The function $\sum_{j =1}^N \left\|  \nabla \phi_j (\x) \right\|_2^2$ is 
the ``diagonal'' of the reproducing kernel of $\mathbb V$. It is independent of the choice of the orthonormal basis and only depends on the space $\mathbb V$ and the measure $\omega$. 
\end{remark}

The following result is a direct consequence of Lemma \ref{thm:Prob} with $\delta = \frac12$.
\begin{corollary}\label{coro:Prob}
The least squares problem \eqref{LSQc} is stable in the following sense: for any $r>0$, 
\begin{equation}\label{Pr2}
{\rm Prob} \bigg\{ \| {\bf A} - {\bf I} \| > \frac12\bigg\} \le 2 K^{-r}, 
\end{equation}
provided that 
\begin{equation}\label{condition1}
{\mathscr K}_N \le \lambda \frac{K}{\log K}, \quad \mbox{with}~~ \lambda := \frac{ 3 \log(3/2) -1 }{2+2 r}.
\end{equation}
\end{corollary}

\subsubsection{Error bound}\label{sec:bound}

To analyze errors in the proposed algorithm, we consider only
noiseless data case.
For noisy data, the analysis is considerably more involved and will be
pursued in a separate study.

For noiseless data, we consider the more practical case when only 
state variable data are available and time derivative data are computed
numerically. 
Since the state variable data are noiseless, the only errors in the
data set \eqref{pairing} are the numerical approximation errors for
the time derivatives, as discussed Section \ref{sec:der}.
We assume the approximation errors ${\bm \tau}_k:=\dot{\bf x}_k - {\bf
  J}^{-1} \nabla H \big( {\bf x}_k \big)$ are uniformly bounded, i.e., 
\begin{equation}\label{assump1}
 \left\| {\bm \tau}_k  \right\|_2 
=  \left\| \dot{\bf x}_k - {\bf J}^{-1} \nabla H \big( {\bf x}_k \big) \right\|_2 \le \tau_\infty, \quad \forall {\x}_k \in D,
\end{equation}
where  
$\tau_\infty<+\infty$ is an assumed bound that  
 depends on the regularity of ${\bf u}(t)$ in the
time interval $[0, J\Delta t]$ and the accuracy of the numerical
differentiation method.

\begin{theorem}\label{thm:error1}
	Assume  
	\begin{equation} \label{H_bound}
	\| \nabla H ({\bf x}) \|_{2} \le L<+\infty, \qquad {\bf x} \in D\quad a.e.  
	\end{equation} 
For any $r>0$, under the condition \eqref{condition1}, it holds that 
$$
\mathbb E \Big( \big\| \nabla H - {\bf T}_L (  \nabla \widetilde H ) \big\|_{2,\mathbb L^2_\omega}^2  \Big) 
\le \left( 1 + \frac{8 \lambda}{\log K} \right) \left\| \nabla  H -     \Pi_{\mathbb V} ( \nabla H ) \right\|^2_{2,\mathbb L_\omega^2} + \frac{8 L^2 }{K^r} + 8 N \tau_\infty^2, 
$$
where the expectation $\mathbb E$ is taken over the random sequences of $\{ {\bf x}_k \}_{k=1}^K$, 
$\lambda $ is defined in \eqref{condition1}, $L$ is the bound defined in \eqref{H_bound},  ${\bf T}_L(\x)$ is defined by 
$$
{\bf T}_L(\x) = \frac{ L }{ \max\{ \| \x \|_2, L \} } \x,
$$
and  
$\Pi_{\mathbb V} (\nabla H)$ denotes the orthogonal projector of $\nabla H$ onto $\mathbb V$, i.e., 
the best approximation to $\nabla H$ in $\mathbb V$,
$$
\Pi_{\mathbb V} (\nabla H) := \mathop{\rm argmin}\limits_{ \nabla h \in \mathbb V } \int_D \left\| \nabla H - \nabla h \right\|_2^2 d \omega. 
$$
\end{theorem}
\begin{proof}
  See Appendix \ref{app:proofthm2}.
\end{proof}

As a direct consequence, we have the following corollary.

\begin{corollary}\label{thm:error2}
	Assume $
	\widetilde{L} := \max\{ \| \nabla H \|_{2,L^\infty}, \| \nabla \widetilde H \|_{2,L^\infty}   \} < +\infty.
	$
	Then, for any $r>0$, under the condition \eqref{condition1},
        the following result holds,
	$$
	\mathbb E \Big( \big\| \nabla H -  \nabla \widetilde H \big\|_{2,\mathbb L^2_\omega}^2  \Big) 
	\le \left( 1 + \frac{8 \lambda}{\log K} \right) \left\| \nabla
          H -     \Pi_{\mathbb V} ( \nabla H ) \right\|^2_{2,\mathbb
          L_\omega^2} + \frac{8 \widetilde{L}^2 }{K^r} + 8 N \tau_\infty^2. 
	$$ 
\end{corollary}

We now discuss error bound for the reconstructed Hamiltonian
\begin{equation}\label{expan2d}
\widetilde H (\x ) = C +  \sum_{j =1}^N c_j  \phi_j ( {\bf x} )=: C+ \widetilde H_0 (\x ).
\end{equation}
Note that Hamiltonian is not unique and is defined up to an additive constant $C$. 
Therefore, the error between $\widetilde H (\x )$ and $ H (\x )$
should be understood in the 
quotient space ${\mathbb H^1_\omega}(D)/\mathbb R$.

\begin{theorem}\label{thm:error3}
	Assume $D$ is a bounded connected open subset of $\mathbb
        R^{2d}$ with Lipschitz boundary and 
	let $d \omega = \frac{1}{\int_D d \x } d \x$. 
	Then, there exists a real constant $C$ such that 
	\begin{equation} \label{eq12}
	 \big \| C+ \widetilde H_0 (\x ) - H(\x)  \big\|_{L^2_\omega }^2 
	\le  C_{D,d}   \big\| \nabla H -  \nabla \widetilde H \big\|_{2,\mathbb L^2_\omega}^2,
	\end{equation}
	where $C_{D,d}$ is a constant depending only on the domain $D$ and the dimensionality $d$.
	 Furthermore, under the assumptions of Corollary
	 \ref{thm:error2}, we have
	 $$\mathbb E \Big( \big \| C+ \widetilde H_0 (\x ) - H(\x)  \big\|_{L^2_\omega }^2  \Big) 
	 \le C_{D,d} \Bigg( \left( 1 + \frac{8 \lambda}{\log K} \right) \left\| \nabla  H -     \Pi_{\mathbb V} ( \nabla H ) \right\|^2_{2,\mathbb L_\omega^2} + \frac{8 L^2 }{K^r} + 8 N \tau_\infty^2 \Bigg).
	$$
\end{theorem}

\begin{proof}
	Let us take 
	$$
	C=\int_D \Big( H(\x) - \widetilde H_0 (\x) \Big) d \omega.  
	$$
	Using the Poincar\'e inequality (cf.~\cite{leoni2017first}) we
        obtain \eqref{eq12}. 
	The proof is then completed upon combining \eqref{eq12} with Corollary
	\ref{thm:error2}.
\end{proof}

%% file: Examples.tex
\section{Numerical Examples} \label{sec:examples}

In this section we present numerical examples to demonstrate the
properties and effectiveness of the proposed method.

In all the test cases, we generate synthetic trajectory data by
solving the underlying
Hamiltonian systems using a high resolution numerical
solver. More specifically, we use the classical fourth-order explicit 
Runge-Kutta method (cf. \cite[p.~131]{griffiths2010numerical}) with a
very small time step of size $0.0001 \Delta t$.  
The proposed numerical method is then applied to the data to produce
the corresponding approximate Hamiltonian systems, whose solutions are
then compared against the solutions to the true
Hamiltonian systems to examine numerical errors. Note that in all of
the tests the only available data are on the solution state
variables. The time derivatives of the states are estimated
numerically using the procedure discussion in Section \ref{sec:der}.

For convenience, we assume the computational domain $D$ to be a hypercube. 
Without loss of generality, we employ polynomial basis functions in all the numerical examples. 
Specifically, 
we set the finite dimensional subspace $\mathbb W$ as ${\mathbb P}_n^{2d}$, the linear space of $2d$-dimensional polynomials of total degrees up to $n\ge 1$. That is, 
$$
{\mathbb P}_n^{2d} = {\rm span}\{ {\bf x}^{\bf i} = x_1^{i_1} \cdots x_{2d}^{i_{2d}},~|{\bf x}|\le n \},
$$
where ${\bf i}=(i_1,\dots,i_{2d})$ is multi-index with $|{\bf
  i}|=i_1+\dots+i_{2d}$. In all examples, we use 
 the tensor products of univariate
  Legendre polynomials as a basis on the hypercube domain $D$, which are commonly used in many practical
  applications. See, for example, \cite{Szego39}.
Although the Legendre polynomials 
 do not satisfy the orthogonality defined at the beginning of Section 3.2, 
our numerical results indicate that they are a good choice.   
Note that the solution 
to the least squares problem \eqref{LSQ} is independent of  basis choice.

The gradient function space $\mathbb V$ is defined  via
\eqref{eq:gradspaceV}, and we have
$$
N=\dim \mathbb V = \dim \mathbb W - 1 = {n+2d \choose 2d} - 1.
$$
The basis functions of $\mathbb V$ are set as $\bm \psi_j (\x)= \nabla
\phi_j (\x)  $, $j=1,\dots,N$, where $\phi_j$ are the Legendre polynomials in $\mathbb W$. 

For noiseless data, we employ second-order finite difference method to
estimate the time derivatives. For noisy data, we use the polynomial
least squares de-noising  \eqref{eq:LS} with a polynomial degree of
$Q=5$. The detail of the time derivative estimation is discussed in See Section \ref{sec:der}.

Once the approximate system \eqref{appsystem} is constructed, we
simulate its trajectories $\tu$ for some arbitrarily chosen initial state
$\u_0^*$, which is not in the training data \eqref{pairing}, 
 and then compare the errors against the trajectories $\u$ produced
by the exact Hamiltonian system from the same initial state $\u_0^*$. All errors are reported as relatively
errors in the following form for any $t\ge 0$: 
$$
\frac{\|\tu(t; \u_0^*) - \u(t; \u_0^*)\|_2}{\|\u(t;\u_0^*)\|_2}.
$$

\subsubsection*{Example 1: Single pendulum}

The Hamiltonian of an ideal single pendulum with unit mass is its total energy
$$
H(p,q)=\frac{1}{2l^2} p^2 + g\,l\, ( 1-\cos q ),
$$
where $l$ is the length of the pendulum, $q$ is the angular displacement of the pendulum from its downward
equilibrium position, $p$ the angular momentum, and $g=9.8$ the gravitational constant. 
The true Hamiltonian formulation of the dynamics is
\begin{equation}
\label{eq:example1}
\begin{cases}
\vspace{3pt}
\dot{p}= - g\,l\, \sin q,\\
\dot{q}= \dfrac{p}{l^2}.
\end{cases}
\end{equation}

We set $l=1$ and the computational domain $D= (-2\pi, 2\pi)\times (-\pi,\pi)$.
The data pairs \eqref{pairing} consist of 
$M=500$ short trajectories, each of which is generated by random
initial state in  $D$ and contains
$J=40$ steps. 
Hereafter, the random
initial states are independently drawn from the uniform distribution over $D$.   
All data are then perturbed by a multiplicative factor
$(1+\eta)$, where $\eta$ is i.i.d. uniform
distributed in $[-0.08, 0.08] $. This corresponds to $\pm 8\%$
relative noise   
in all data.

The Hamiltonian $\widetilde{H}(\cdot)$ is approximated with
polynomials of degree up to $n=6$. The numerical solution of the
approximate Hamiltonian system is denoted as $\tu(t;\u_0)$.
To assess the accuracy of the algorithm, we set an arbitrarily chosen
initial state $\mathbf{u}_0^*=(-3.876, -1.193)^\top$ and solve both
the approximate solution $\tu(t;\u_0^*)$ and the exact solution
$\u(t;\u_0^*)$. For cross-comparison, we also implemented the equation
approximation algorithm from \cite{WuXiu_JCPEQ18}, which directly
approximates the right-hand-side of the unknown governing equations
and	thus, in general, does not preserve any Hamiltonian.
We denote this solution  
$\hu(t; \u_0^*)$.

\begin{figure}[htbp]
	\centering
	\subfigure[Evolution of relative errors]{\includegraphics[width=0.48\textwidth]{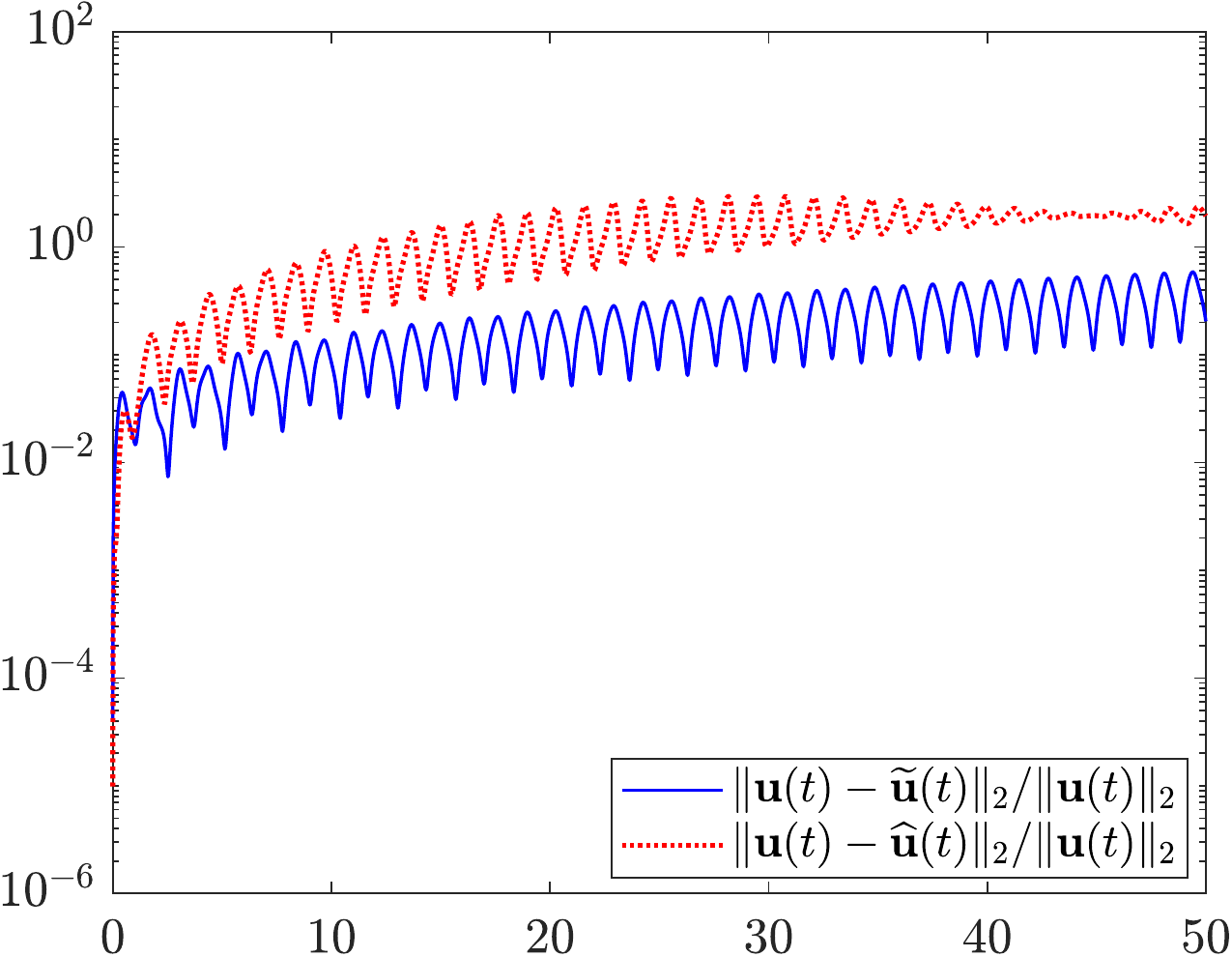}}
	\subfigure[Evolution of Hamiltonian]{\includegraphics[width=0.48\textwidth]{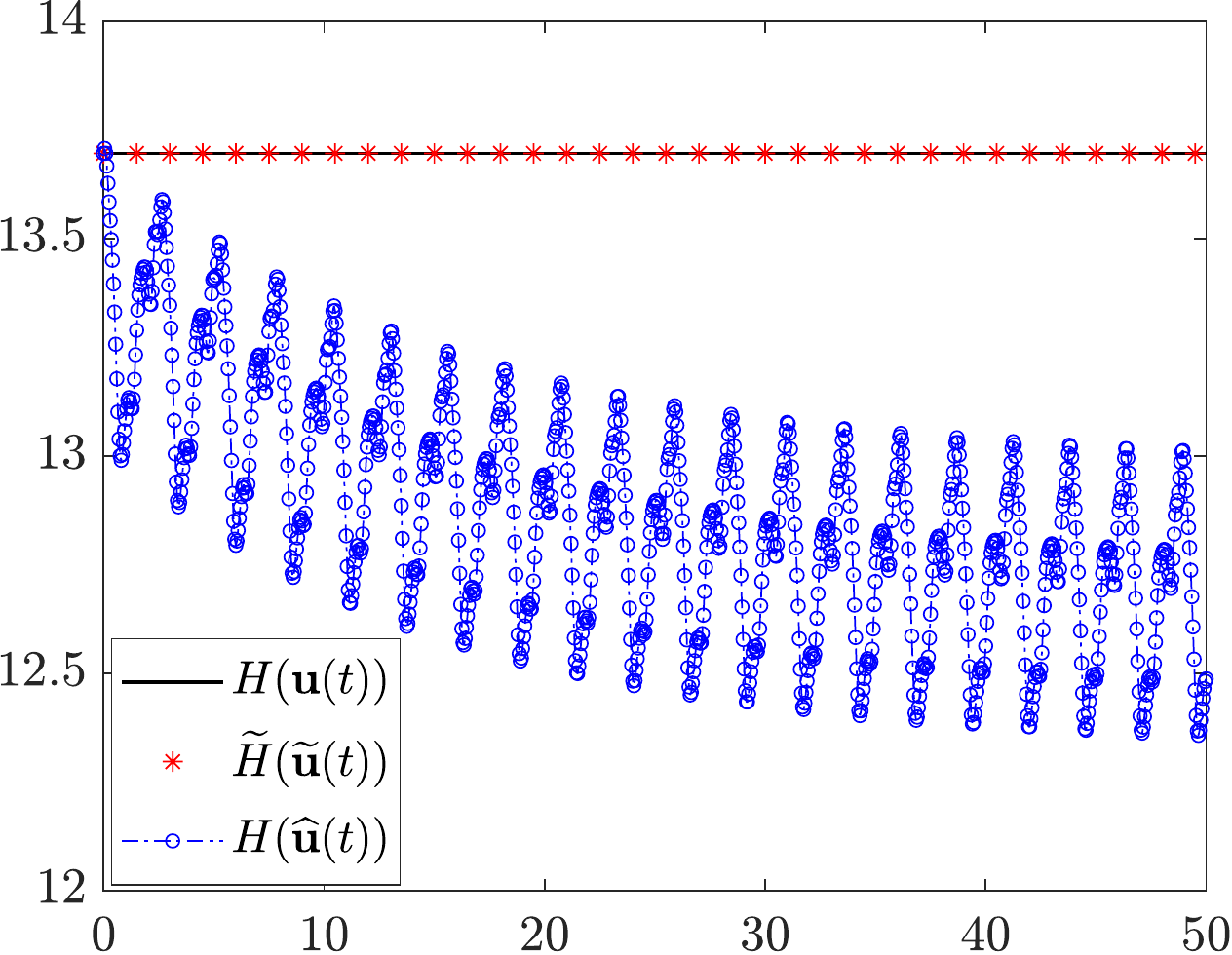}}
	\caption{\small
		Example 1: Solutions of the reconstructed system with an initial state $\mathbf{u}_0^*=(-3.876, -1.193)^\top$. Left:
		relative errors against the true solution by the SP
		method ($\tu$) and non-SP method ($\hu$); Right: time
		evolution of the Hamiltonian.
	}\label{fig:ex1_HErrU}
\end{figure}

\begin{figure}[htbp]
	\centering
	\subfigure[$\tau = 1 \times 10^{-3}$]{\includegraphics[width=0.32\textwidth]{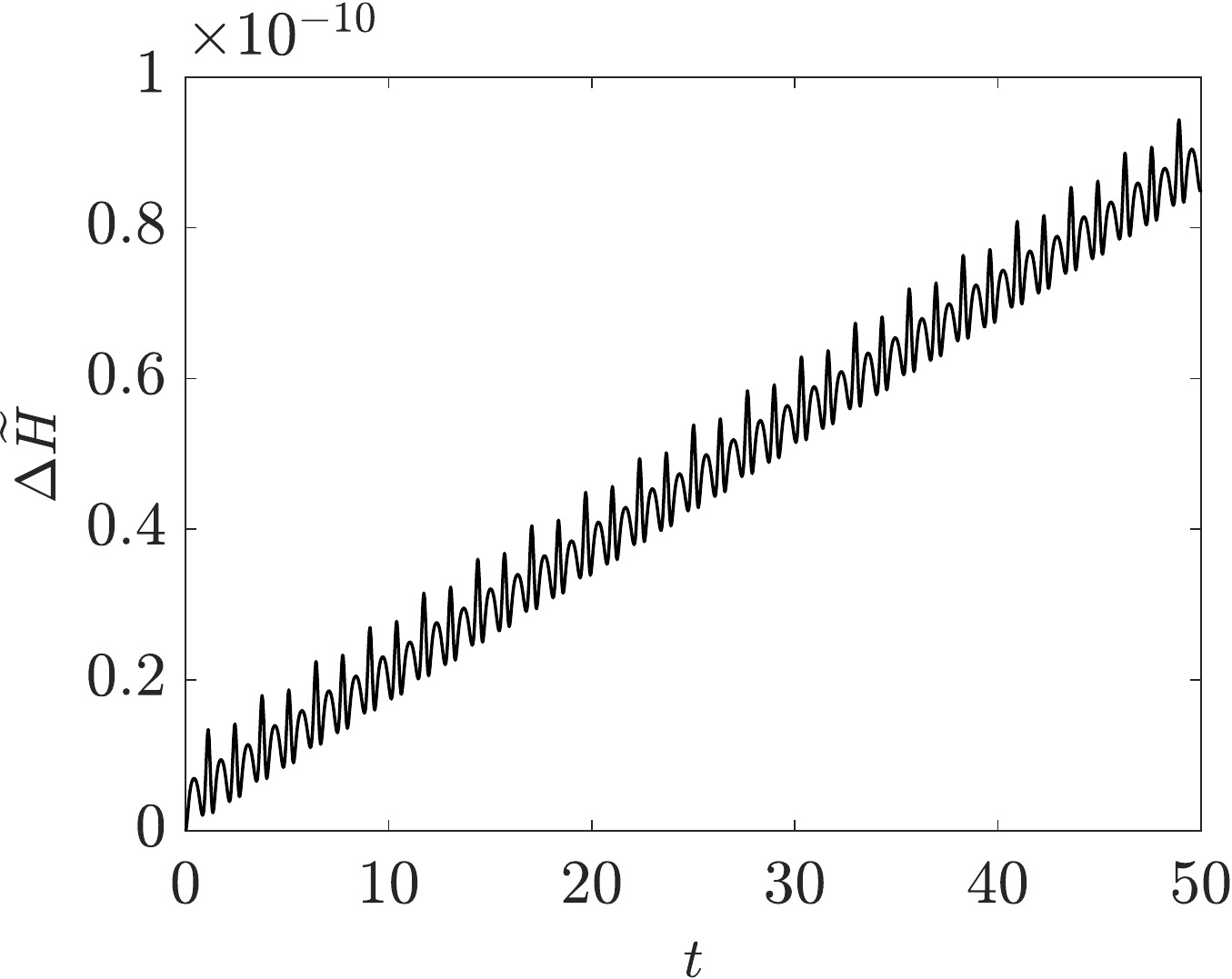}}
	\subfigure[$\tau = 5 \times 10^{-4}$]{\includegraphics[width=0.32\textwidth]{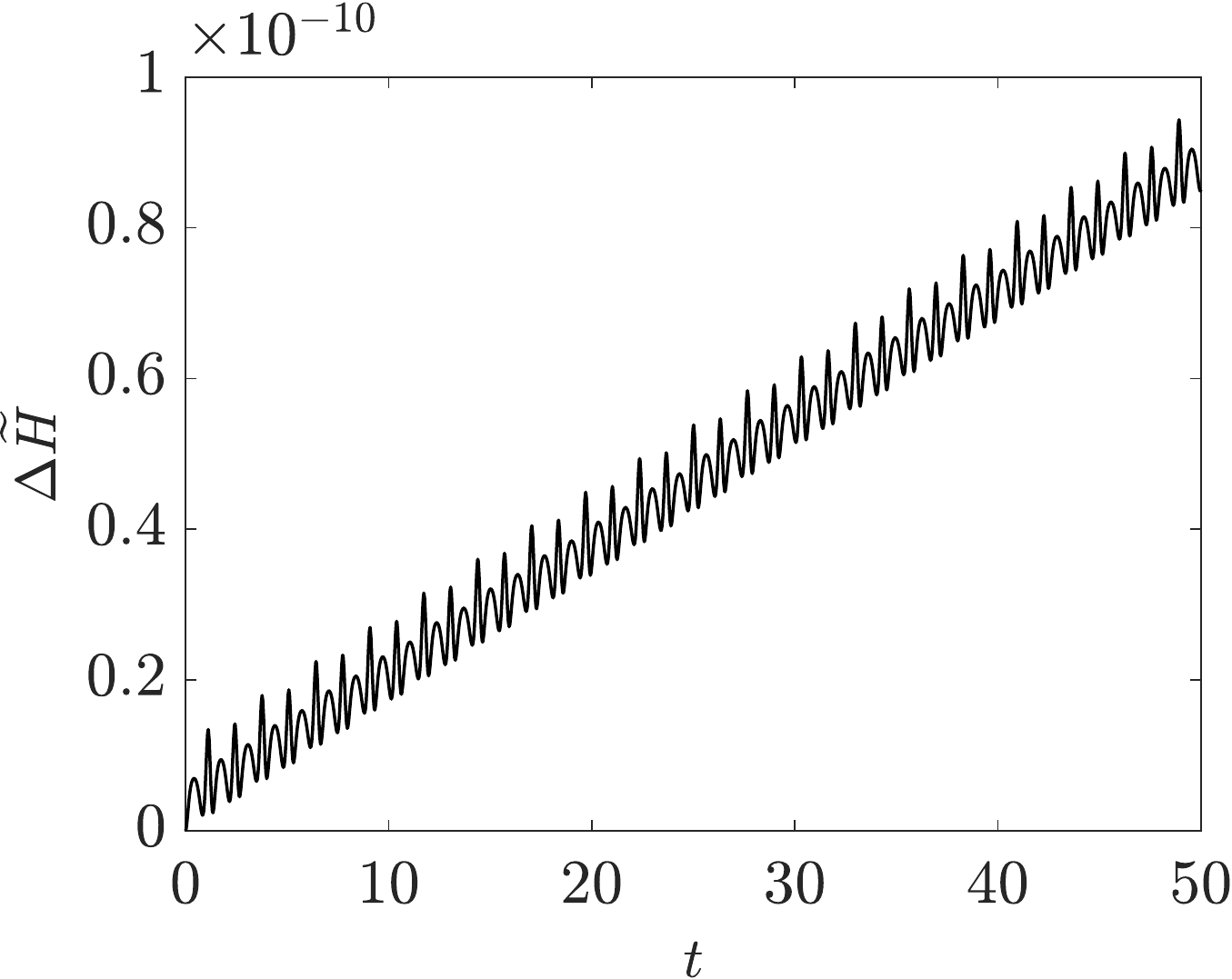}}
	\subfigure[$\tau = 2.5 \times 10^{-4}$]{\includegraphics[width=0.32\textwidth]{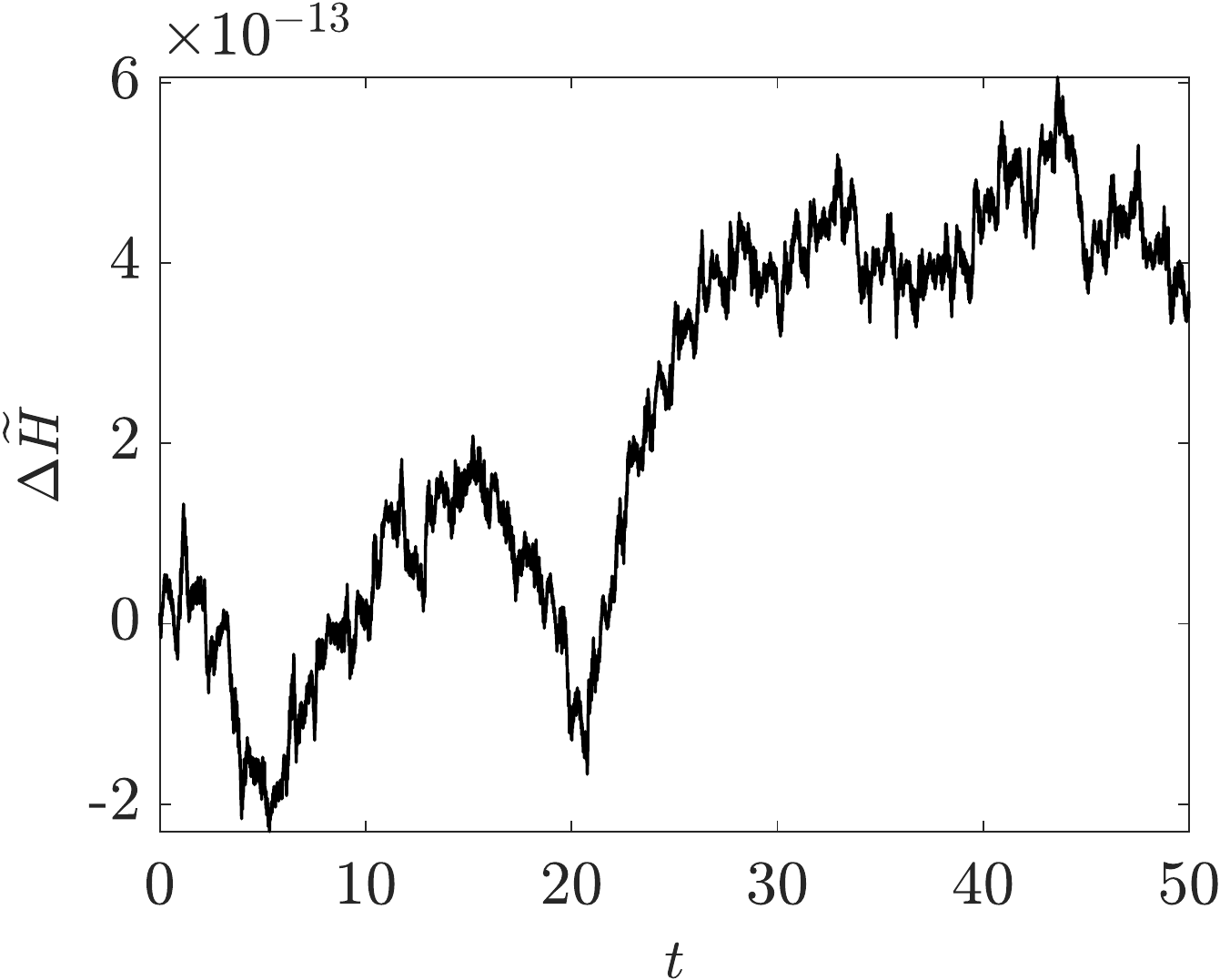}}
	\caption{\small 
		Example 1: evolution of the Hamiltonian deviation $\Delta \widetilde H(t) $ of the reconstructed system computed by the fourth-order Runge-Kutta solver with different time step-sizes $\tau$.
	}\label{fig:ex1_DeltaH}
\end{figure}

In \figref{fig:ex1_HErrU}(a), we plot the evolution of the relative errors in
the numerical solutions. We clearly observe that the errors in our
structure-preserving (SP) algorithm is notably smaller than the non-SP
algorithm from \cite{WuXiu_JCPEQ18}. In \figref{fig:ex1_HErrU}(b), we
examine the time evolution of the Hamiltonians. The exact Hamiltonian
is obviously conserved along the exact solution trajectory,
i.e., $H(\u(t; \u_0^*)) = H(\u_0^*)$.
As expected from Theorem \ref{thm:Hconst}, the approximate Hamiltonian
	$\tH$ of the recovered system \eqref{appsystem} is also exactly preserved along its trajectory. The only (small) errors in the computed Hamiltonian 
	may (merely) arise from the ODE solver, which 
	we employ to numerically solve the reconstructed system; see Figure \ref{fig:ex1_DeltaH} for the Hamiltonian deviation $\Delta \widetilde H(t) := \widetilde H( \widetilde {\bf u}( t;{\bf u_0}) ) 
	-  \widetilde H ( \widetilde {\bf u}_0 ) $ computed by the classical fourth-order explicit Runge-Kutta solver with different time step-sizes $\tau$. 
	We clearly observe that the errors in the Hamiltonian deviation decrease quickly   
	 as we reduce $\tau$, and the errors are close to the level of round-off error when $\tau=2.5 \times 10^{-4}$. 
	 We also list the $L^\infty$, $L^2$, and total variation norms of 
	 the computed $\Delta \widetilde H(t)$ in Table \ref{tab:ex1}, which 
	 shows that the errors in the Hamiltonian deviation converge to zero 
	 at a order related to the employed ODE solver. 
	 These results further confirm that the recovered 
	 system \eqref{appsystem} does exactly preserve 
	 the approximate Hamiltonian
	 $\tH$.  
Note that the non-SP method from \cite{WuXiu_JCPEQ18}, albeit quite accurate, generally does not
preserve or relate to any Hamiltonian. 
For the present test case, we have examined that the system recovered by the non-SP method, denoted by 
\begin{equation*}
\begin{cases}
\vspace{3pt}
\dot{p}= g(p,q),\\
\dot{q}= h(p,q),
\end{cases}
\end{equation*}	
	 is not a Hamiltonian system, because it does not satisfy 
	 $\nabla_p g + \nabla_q h = 0 $, so that there is no Hamiltonian $\widehat H(p,q)$ satisfying 
	 $-\nabla_q \widehat H = g$ and $\nabla_p \widehat H = h$. Note that the data are noisy, and the approximate functions $g(p,q)$ and $h(p,q)$ are {\it not} univariate, {\it not} as the  functions in the true system \eqref{eq:example1}.

\begin{figure}[htbp]
	\centering
	\subfigure[$p(t)$ of the learned SP system]{\includegraphics[width=0.48\textwidth]{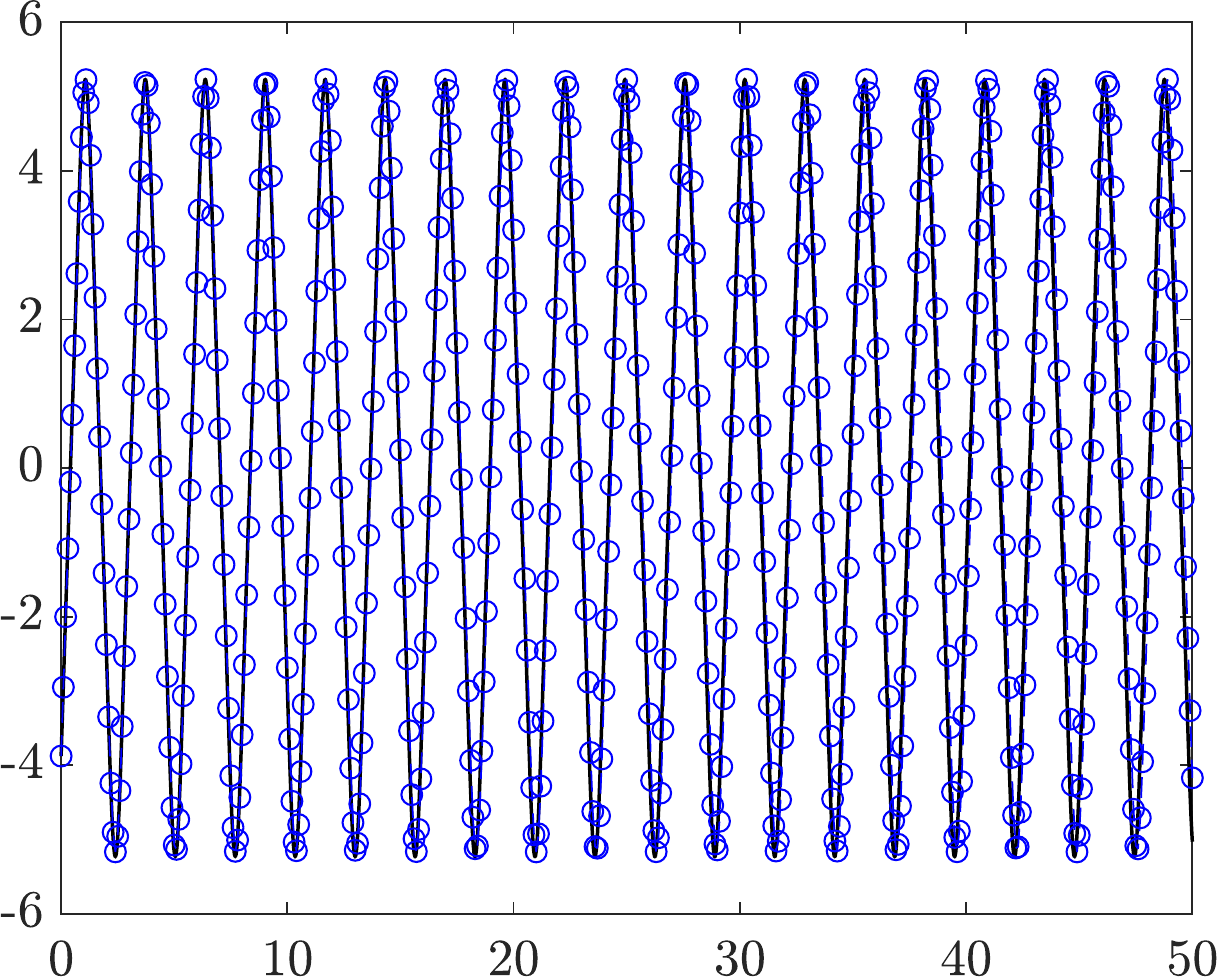}}
	\subfigure[$q(t)$ of the learned SP system]{\includegraphics[width=0.48\textwidth]{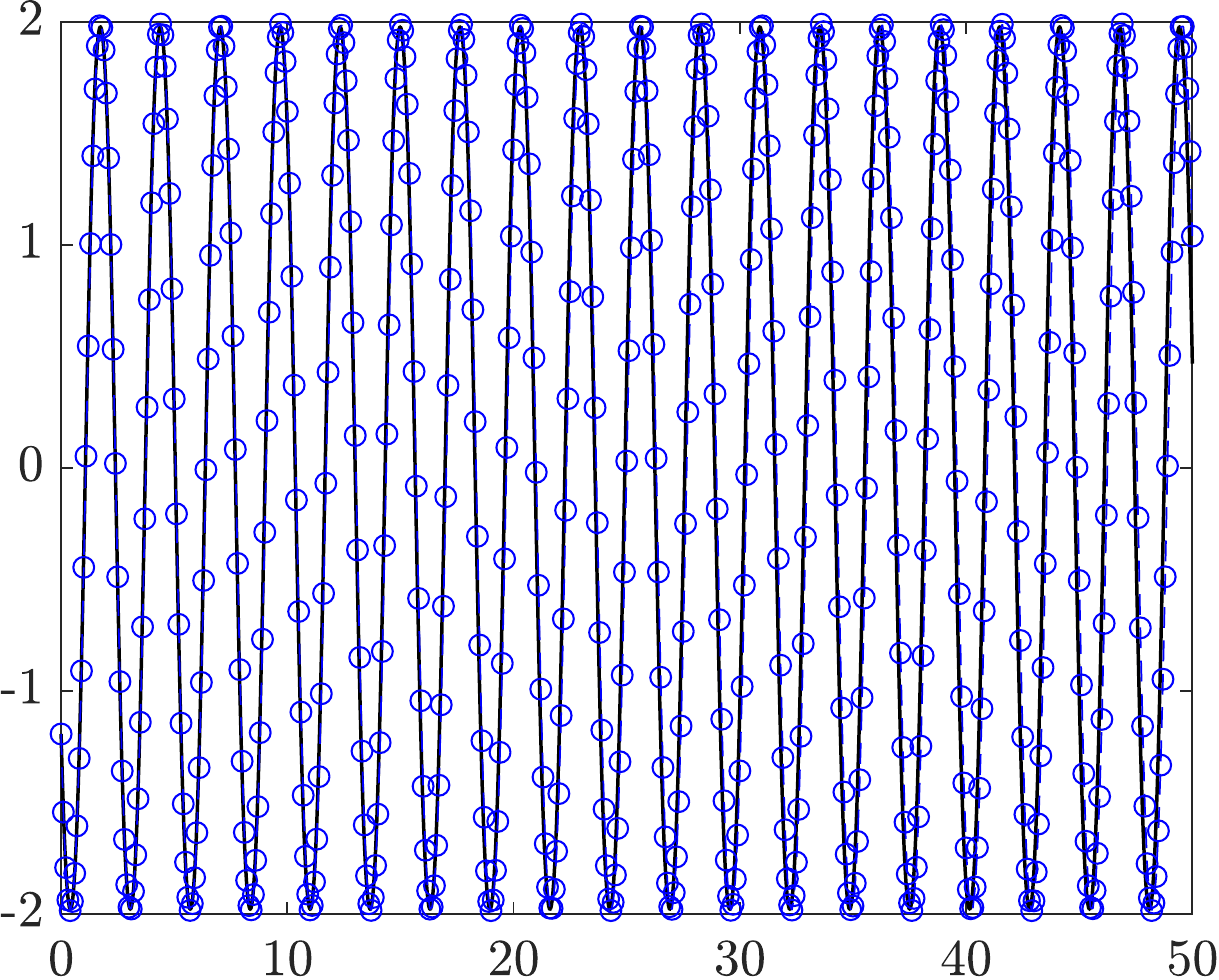}}
	\subfigure[$p(t)$ of the learned non-SP system]{\includegraphics[width=0.48\textwidth]{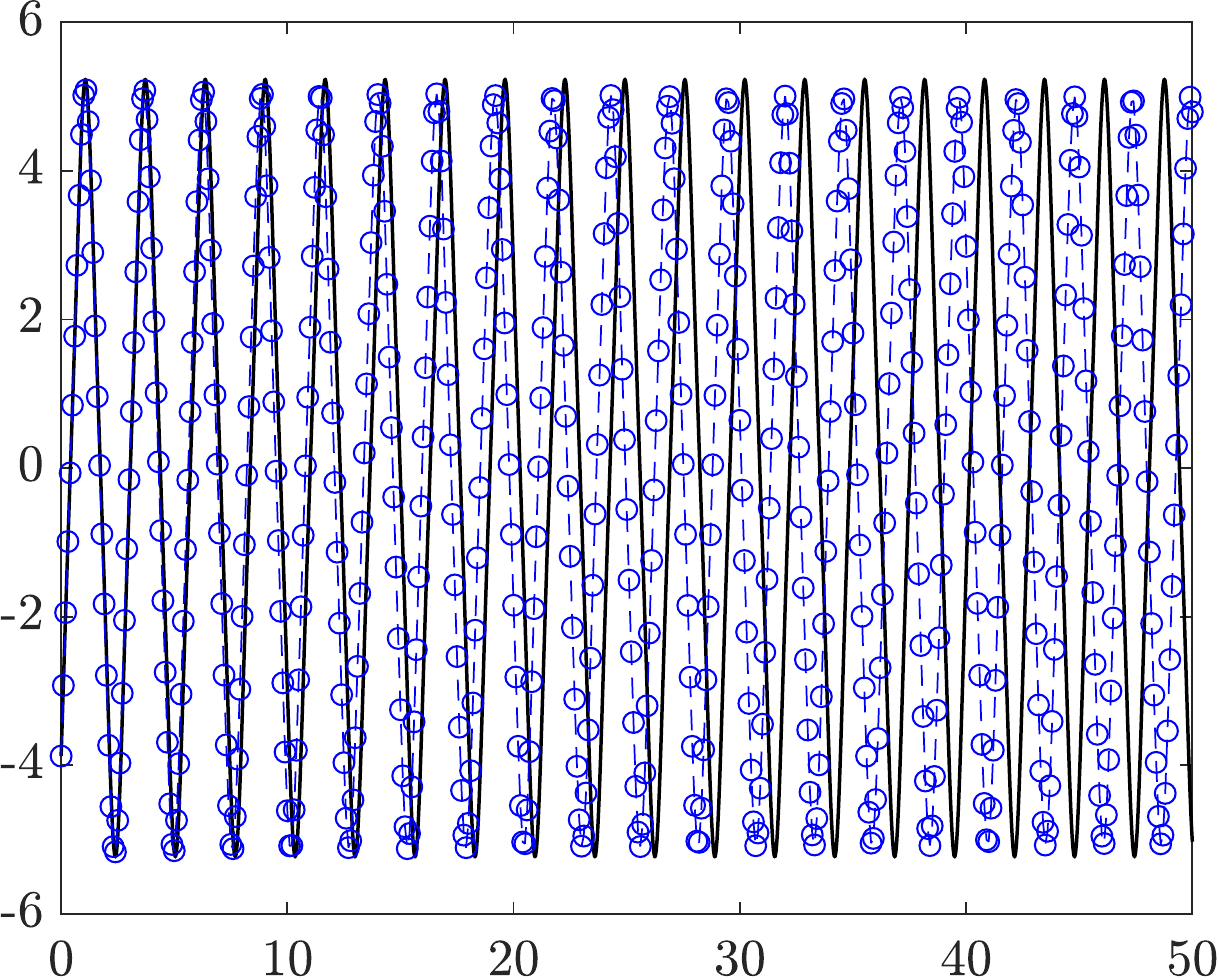}}
	\subfigure[$q(t)$ of the learned non-SP system]{\includegraphics[width=0.48\textwidth]{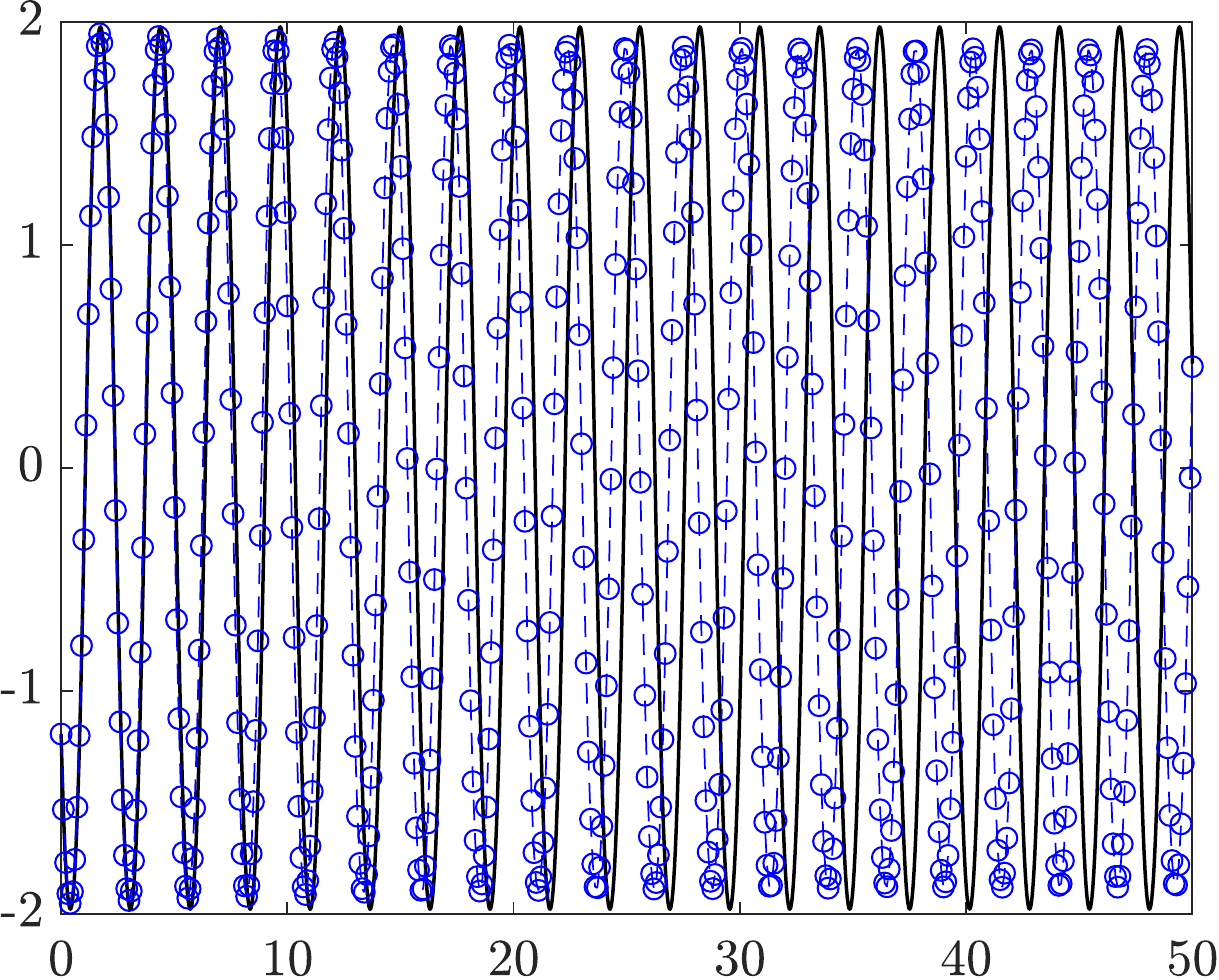}}	
	\caption{\small
		Example 1: Long-term solution of the reconstructed
                system by the SP algorithm (top plots) and non-SP
                algorithm (bottom plots), with 
		initial state $\mathbf{u}_0^*=(-3.876, -1.193)^\top$. (Solid lines represent solution of the true system.)
	}\label{fig:ex1_U}
\end{figure}

The advantage of the proposed SP algorithm is more notable in 
\figref{fig:ex1_U}, we present system predictions over longer
time. The SP method is able to accurately capture the phase of the
solution much better than the non-SP method in \cite{WuXiu_JCPEQ18}.

\begin{figure}[htbp]
	\centering
	\subfigure[$p(t)$ of the SP method without de-noising]{\includegraphics[width=0.48\textwidth]{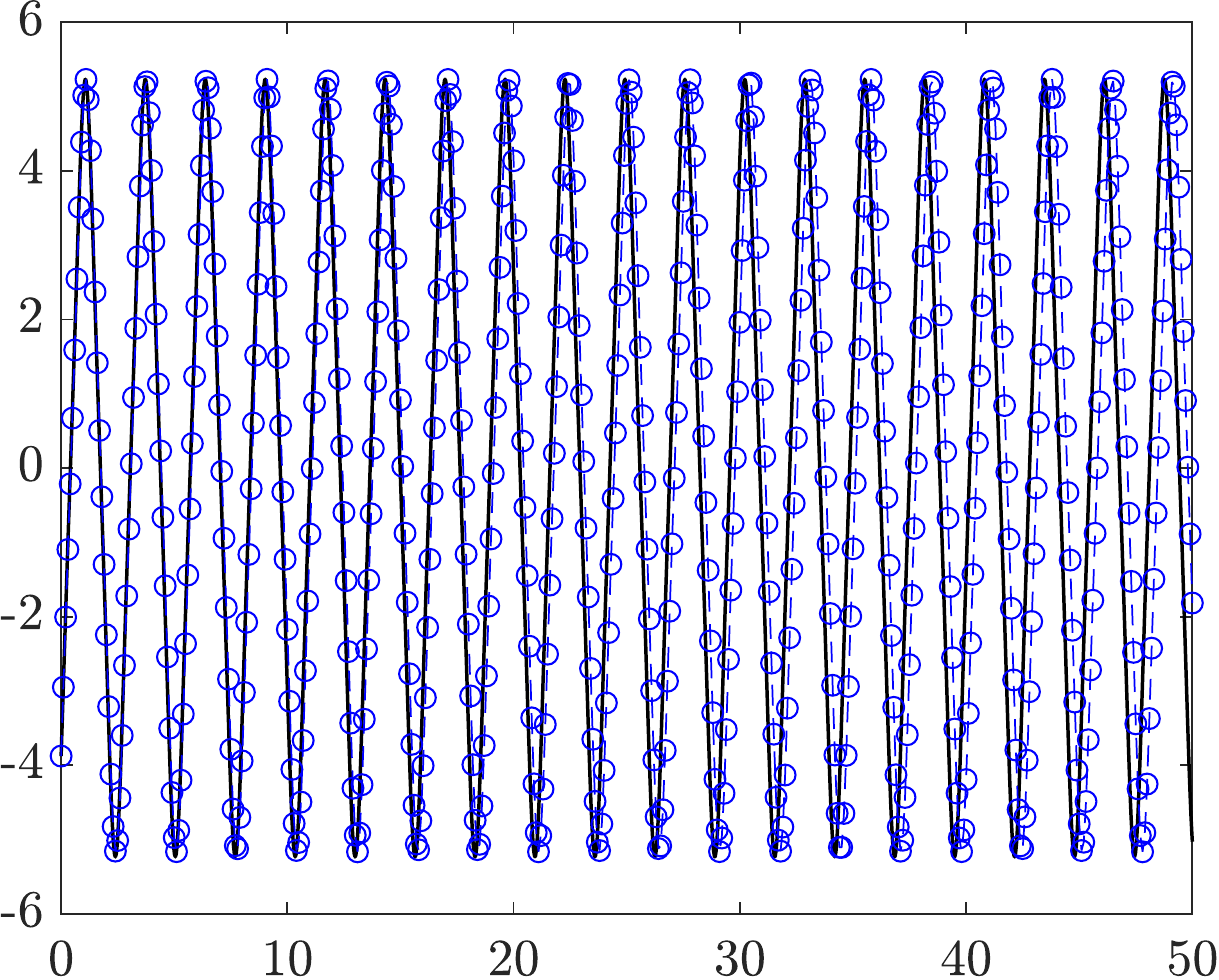}}
	\subfigure[$q(t)$ of the SP method without de-noising]{\includegraphics[width=0.48\textwidth]{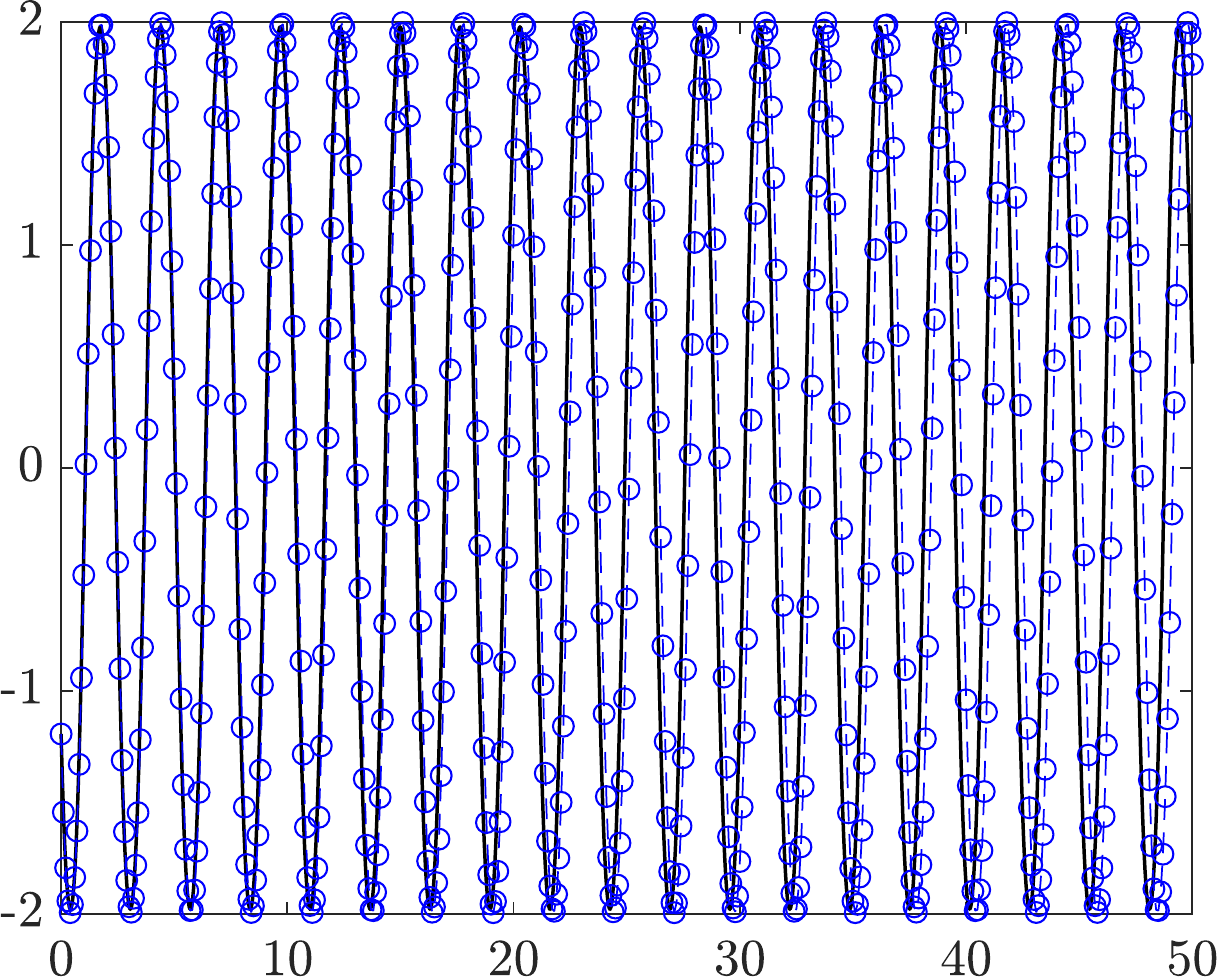}}	
	\caption{\small
		Example 1:  Long-term solution of the reconstructed
                systems by the SP algorithm  without the de-noising procedure \eqref{filter}. (Solid lines represent solution of the true system.)
	}\label{fig:ex1_U_NoF}
\end{figure}

\begin{figure}[htbp]
	\centering
	\includegraphics[width=0.6\textwidth]{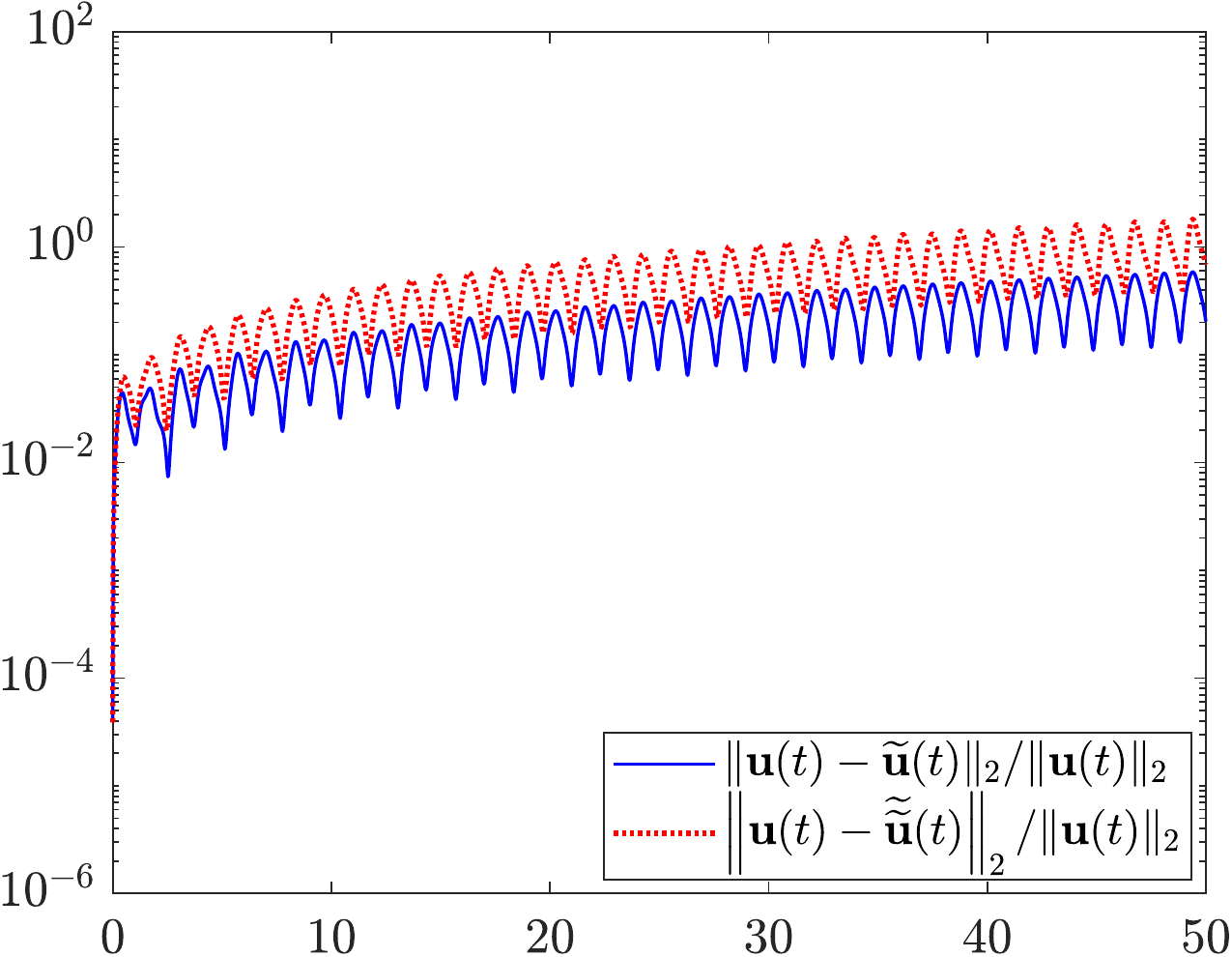}
	\caption{\small
		Example 1: 
		Evolution of relative errors against the true solution by the SP
		method with the de-noising procedure ($\tu$) and without the procedure ($\widetilde{\tu}$), respectively.
	}\label{fig:ex1_HErrU_filter}
\end{figure}

Note that the de-noising procedure
\eqref{filter} has been applied in the computation. For comparison, we
also apply the proposed SP learning method without using the
  de-noising procedure \eqref{filter}.
  The results are plotted in \figref{fig:ex1_U_NoF}. Direct comparison
  of the numerical errors obtained by the two approaches is shown \figref{fig:ex1_HErrU_filter}. It is evident
  that the results obtained without the de-noising procedure are less
  accurate than those by using de-noising.

\begin{table}[htbp]
	\centering
	\caption{\small  
		Example 1: The errors and convergence rates of the Hamiltonian deviation $\Delta \widetilde H(t) $
		computed by the fourth-order Runge-Kutta solver with different time step-sizes $\tau$. 
			}\label{tab:ex1}
	\begin{tabular}{c|c|c|c|c|c|c}
		\hline
		$\tau$ & $L^\infty$-errors& order & $L^2$-errors & order & total variation & order \\
		\hline
		$8 \times 10^{-3}$        &2.7493e-7& -- & 1.0668e-6& --      &  4.0899e-8  &-- \\
		$4 \times 10^{-3}$    & 2.1146e-8&    3.70 & 8.3070e-8& 3.68 & 1.2806e-9 &  5.00 \\
		$2 \times 10^{-3}$   & 1.4455e-9&  3.87 & 5.7047e-9& 3.86 & 4.0064e-11  &   5.00 \\
		$1 \times 10^{-3}$    & 9.4345e-11&  3.94 & 3.7331e-10 &  3.93  & 1.2632e-12 &  4.99  \\
		$5 \times 10^{-4}$  &5.9859e-12&    3.98 & 2.4324e-11 &  3.94 & 1.7065e-13  & 2.89 \\		
		\hline
	\end{tabular}
\end{table}

\subsubsection*{Example 2}

We now consider the following Hamiltonian system  
\begin{equation}
	\label{eq:example2}
	\begin{cases}
		\dot{p}=  4 \alpha_2 q^3 \exp \left( - \alpha_1 p^2 - \alpha_2 q^4 \right),\\
		\dot{q}= - 2 \alpha_1 p \exp \left( - \alpha_1 p^2 - \alpha_2 q^4 \right),
	\end{cases}
\end{equation}
whose Hamiltonian is
$$
H(p,q) = \exp \left( - \alpha_1 p^2 - \alpha_2 q^4 \right).
$$

\begin{figure}[htbp]
	\centering
	\subfigure[Evolution of relative errors]{\includegraphics[width=0.48\textwidth]{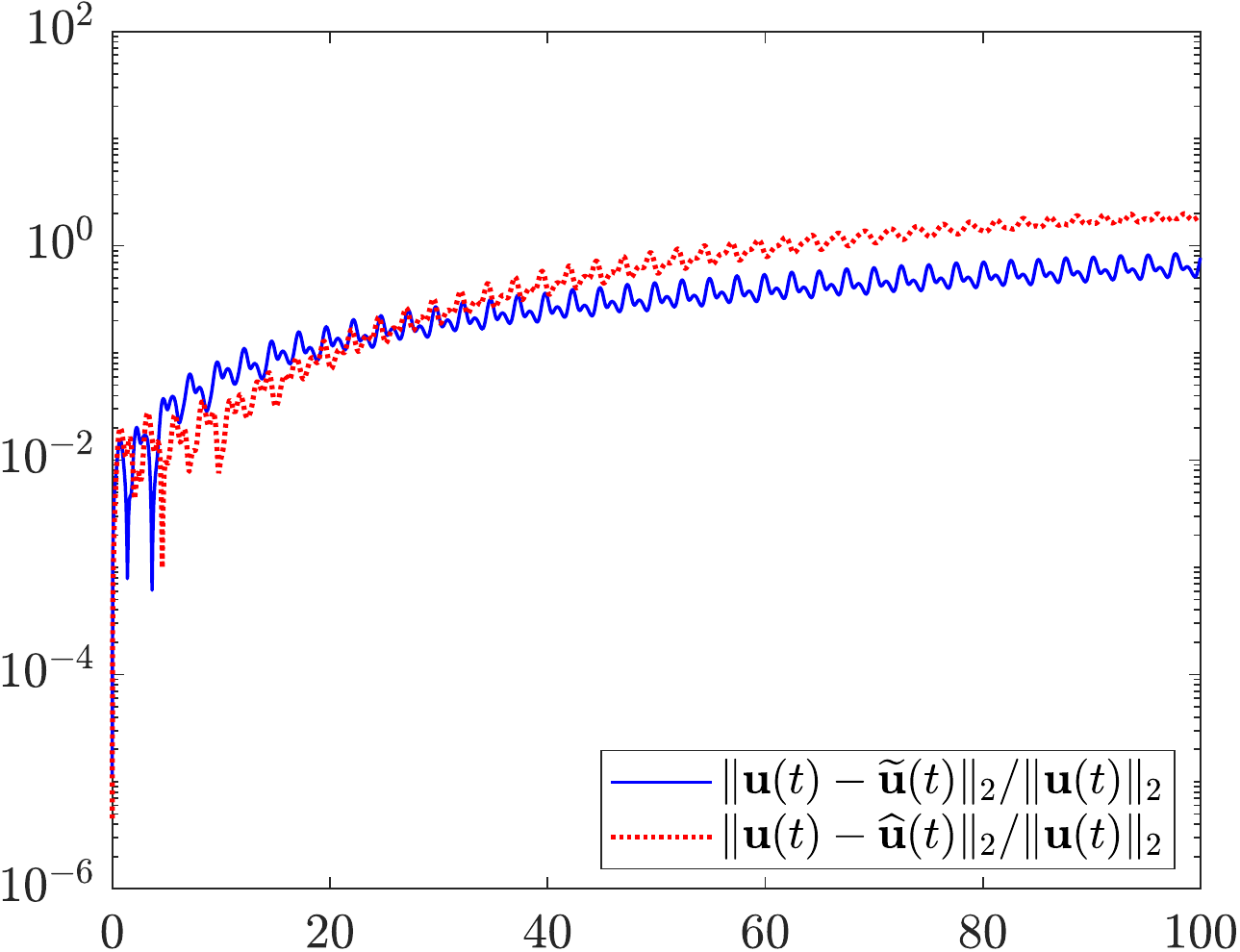}}
	\subfigure[Evolution of Hamiltonian]{\includegraphics[width=0.48\textwidth]{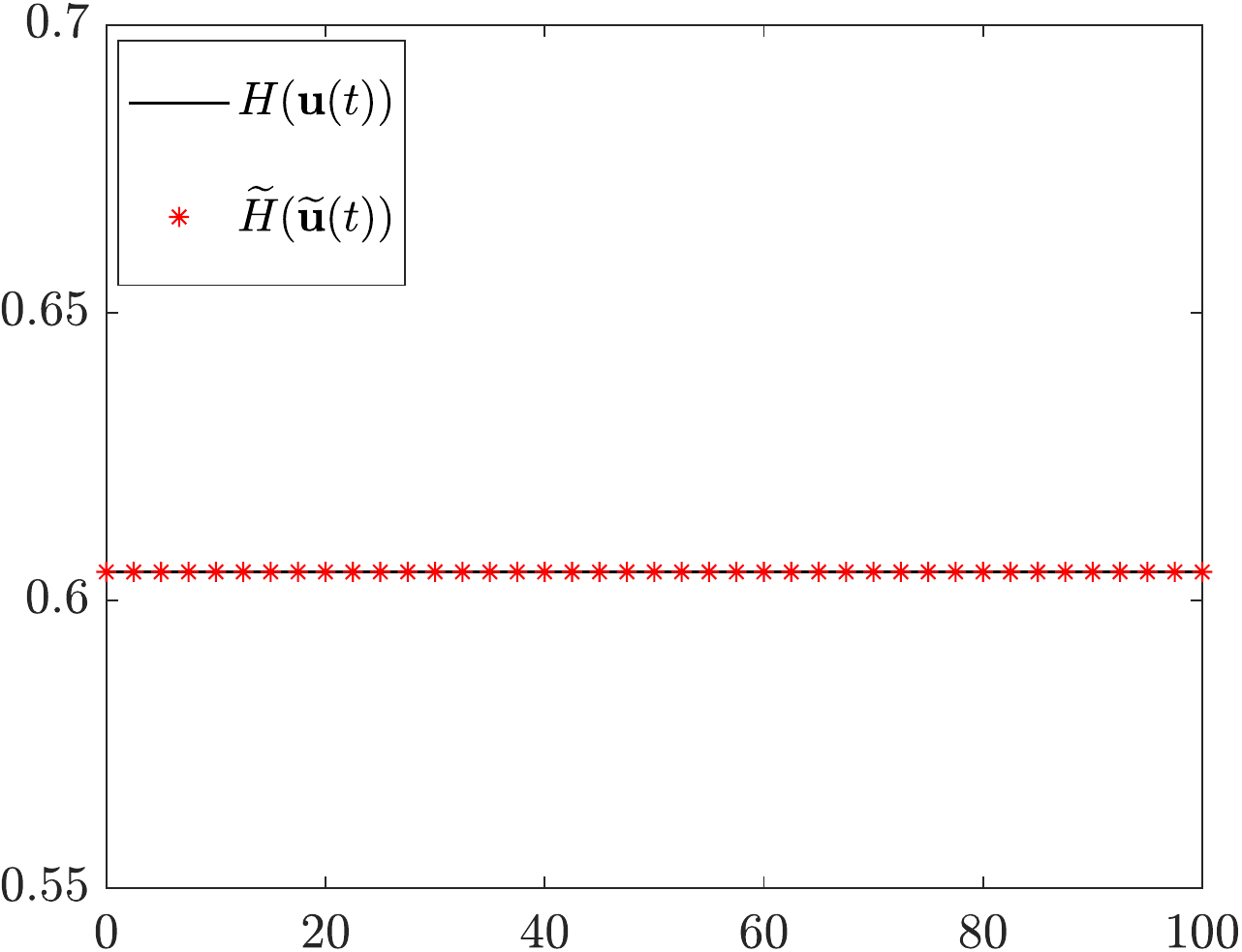}}
	\caption{\small
		Example 2: Solutions of the reconstructed system with an
                initial state ${\bf u}_0^*=( 0.6,0.6 )^\top$. Left:
		relative errors against the true solution by the SP
                method ($\tu$) and non-SP method ($\hu$); Right: time
                evolution of the Hamiltonian.
	}\label{fig:ex2_HErrU}
\end{figure}

\begin{figure}[htbp]
	\centering
	\subfigure[$p(t)$ of the learned SP system]{\includegraphics[width=0.48\textwidth]{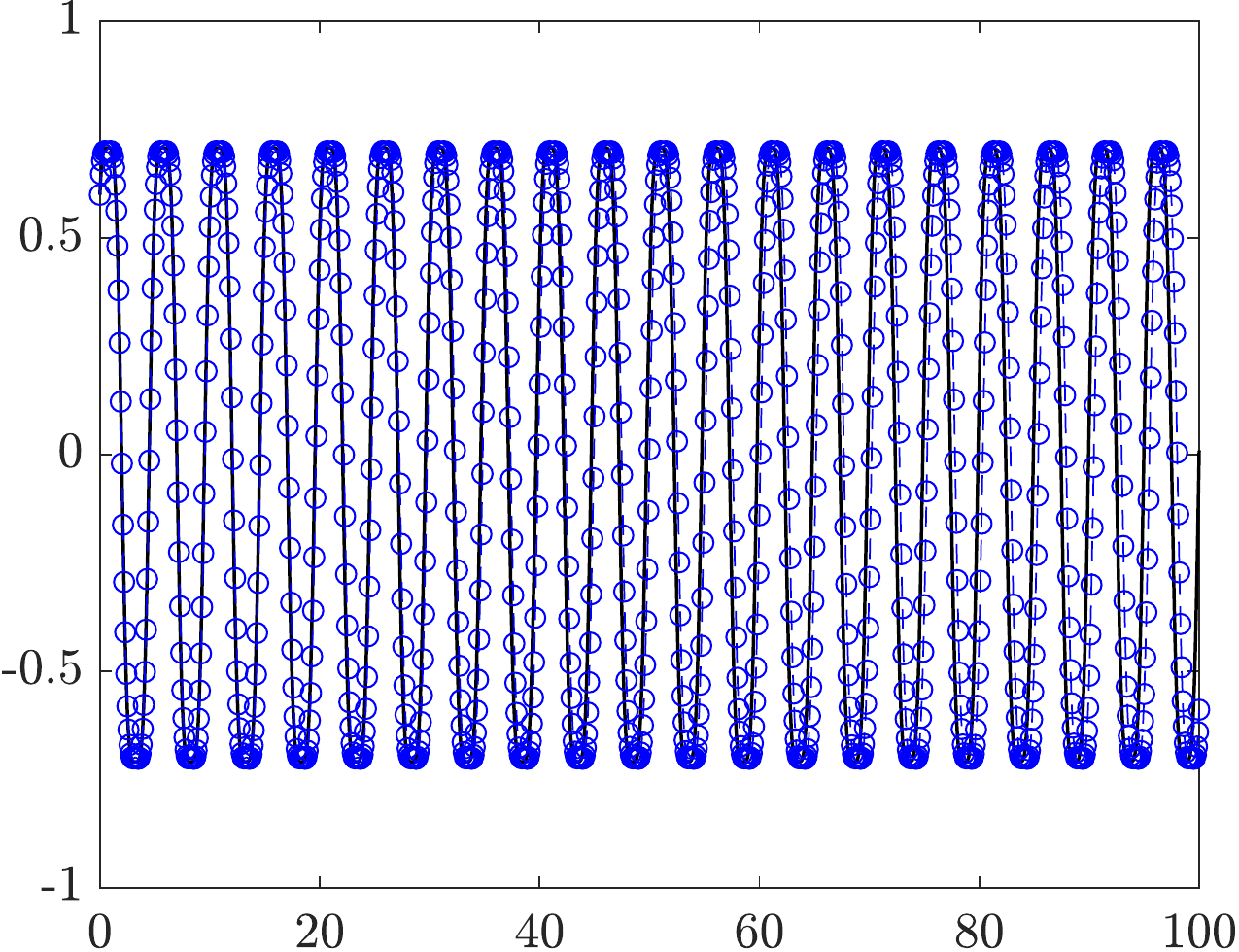}}
	\subfigure[$q(t)$ of the learned SP system]{\includegraphics[width=0.48\textwidth]{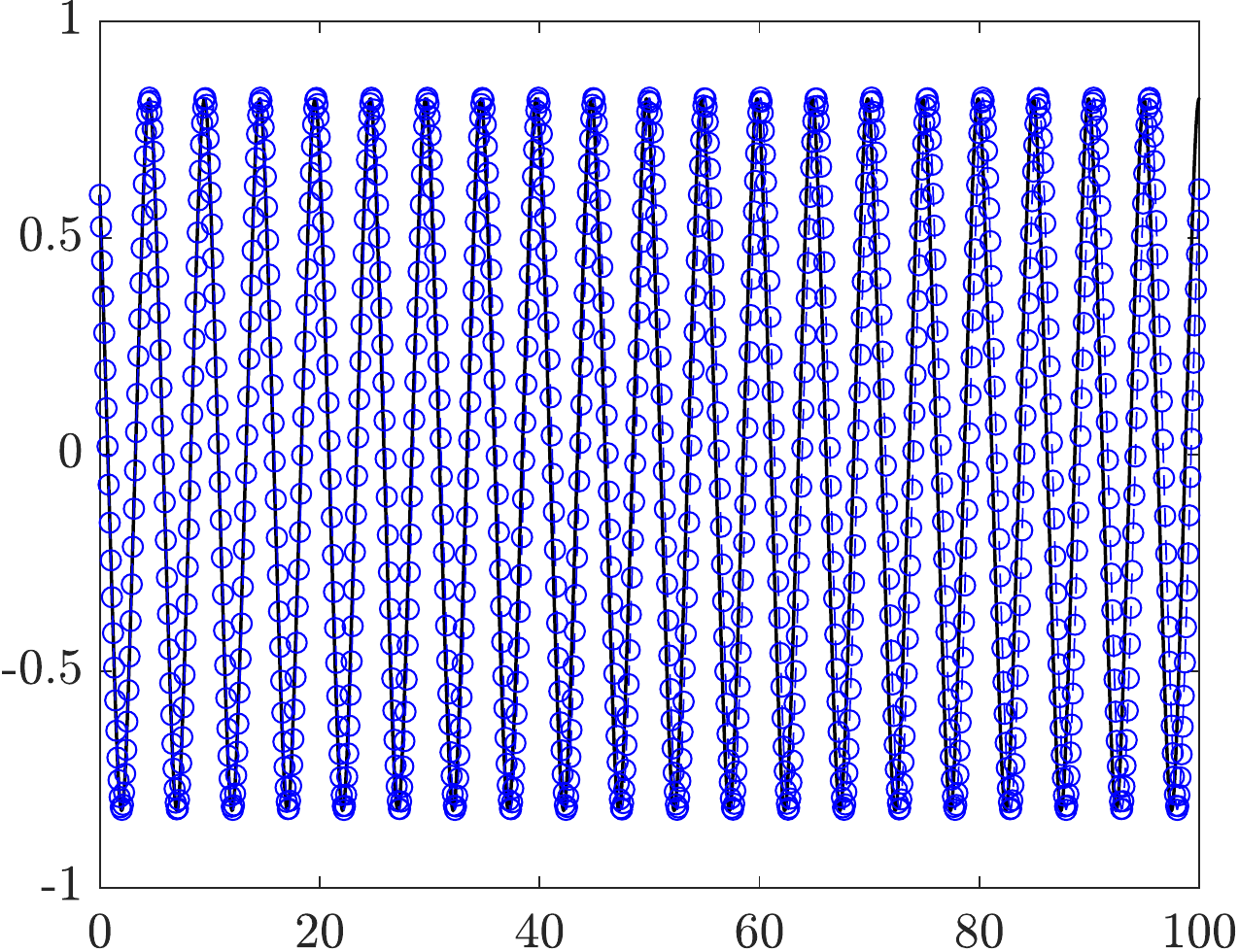}}
	\subfigure[$p(t)$ of the learned non-SP system]{\includegraphics[width=0.48\textwidth]{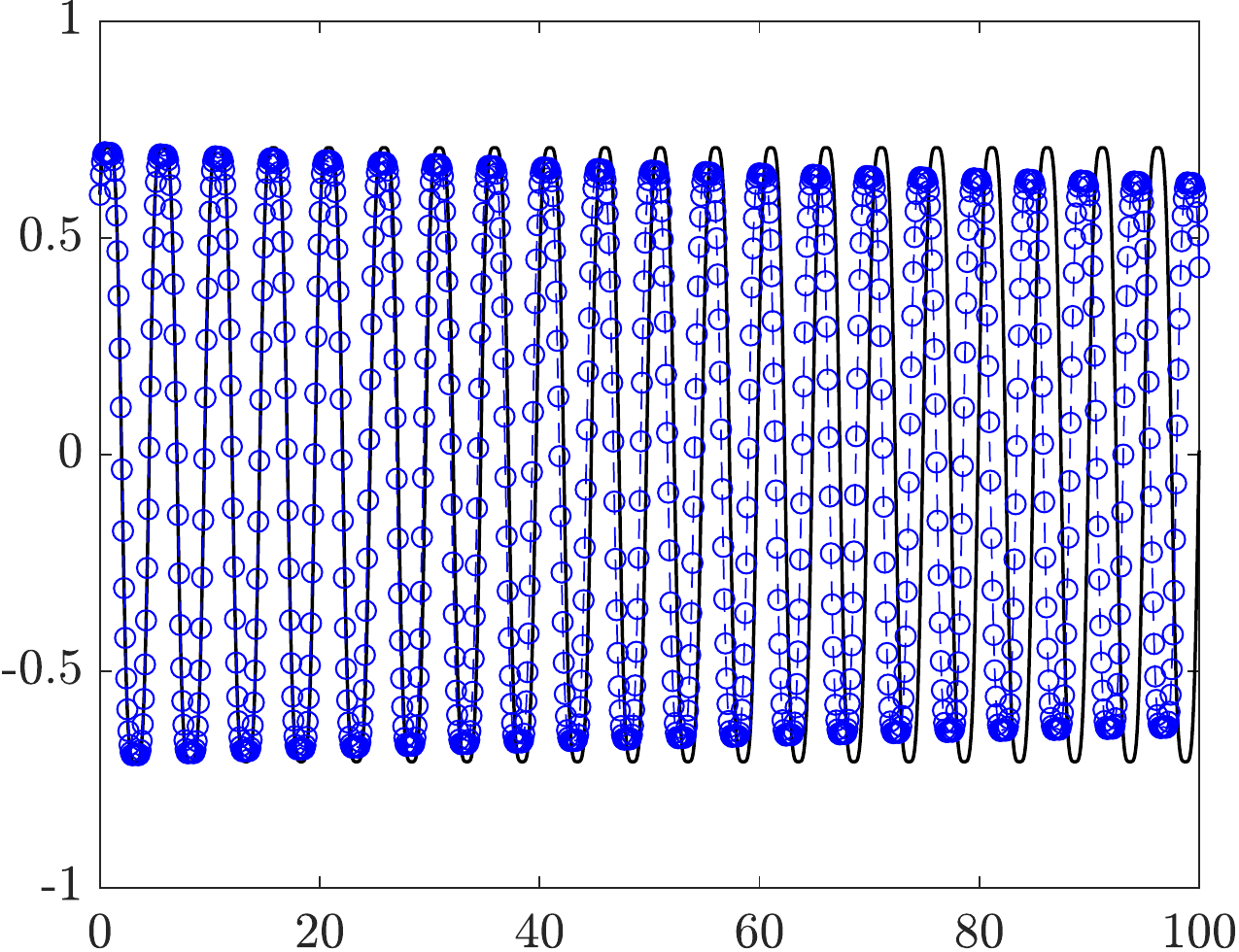}}
	\subfigure[$q(t)$ of the learned non-SP system]{\includegraphics[width=0.48\textwidth]{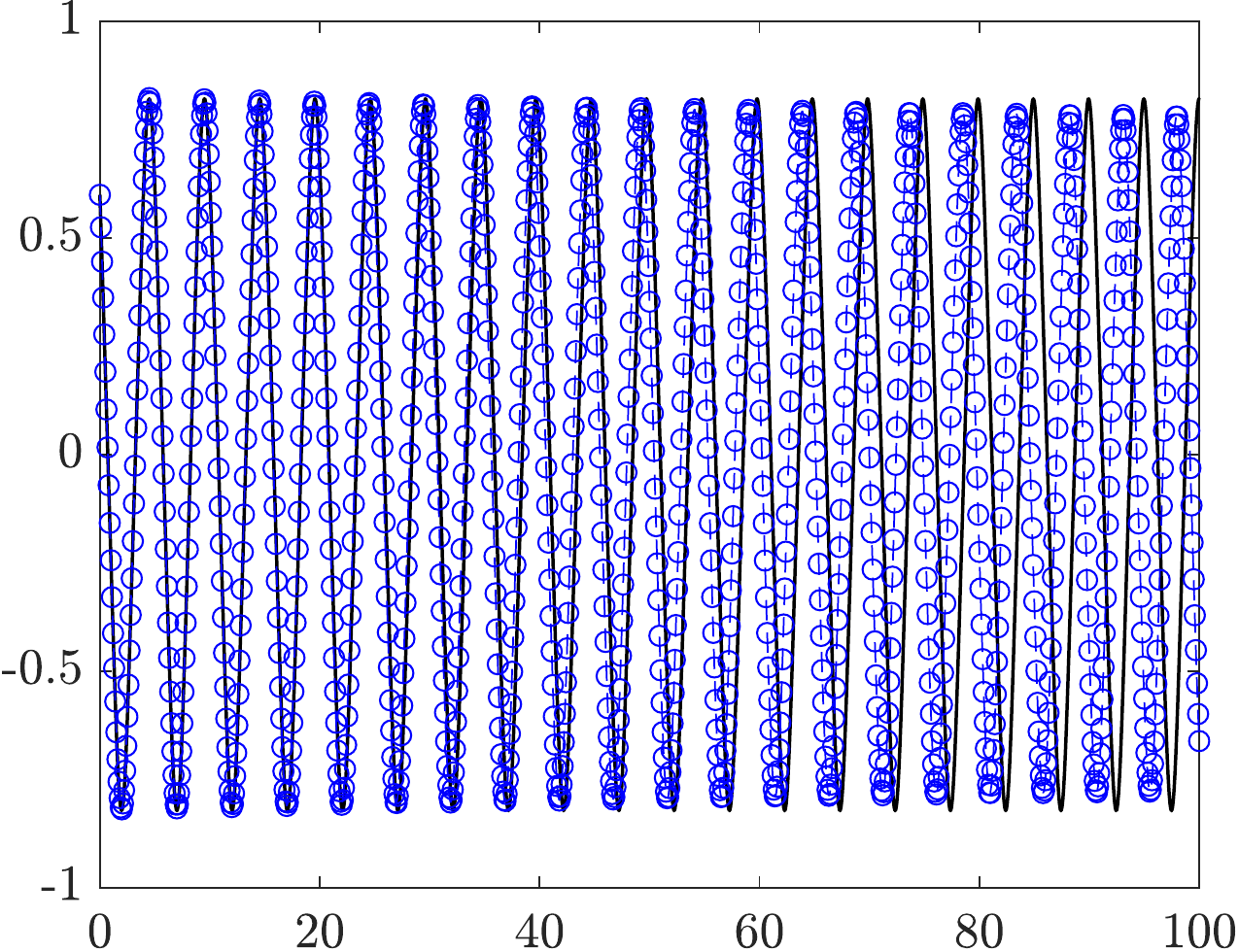}}	
	\caption{\small
		Example 2: Comparison of solutions of the learned SP and non-SP systems with the solution of 
		the true system for the initial state ${\bf u}_0^*=( 0.6,0.6 )^\top$. Solid lines represent solution of the true system and the (blue) circles are solutions of the approximate ones.
	}\label{fig:ex2_U}
\end{figure}

\begin{figure}[htbp]
	\centering
	\subfigure[Exact]{\includegraphics[width=0.32\textwidth]{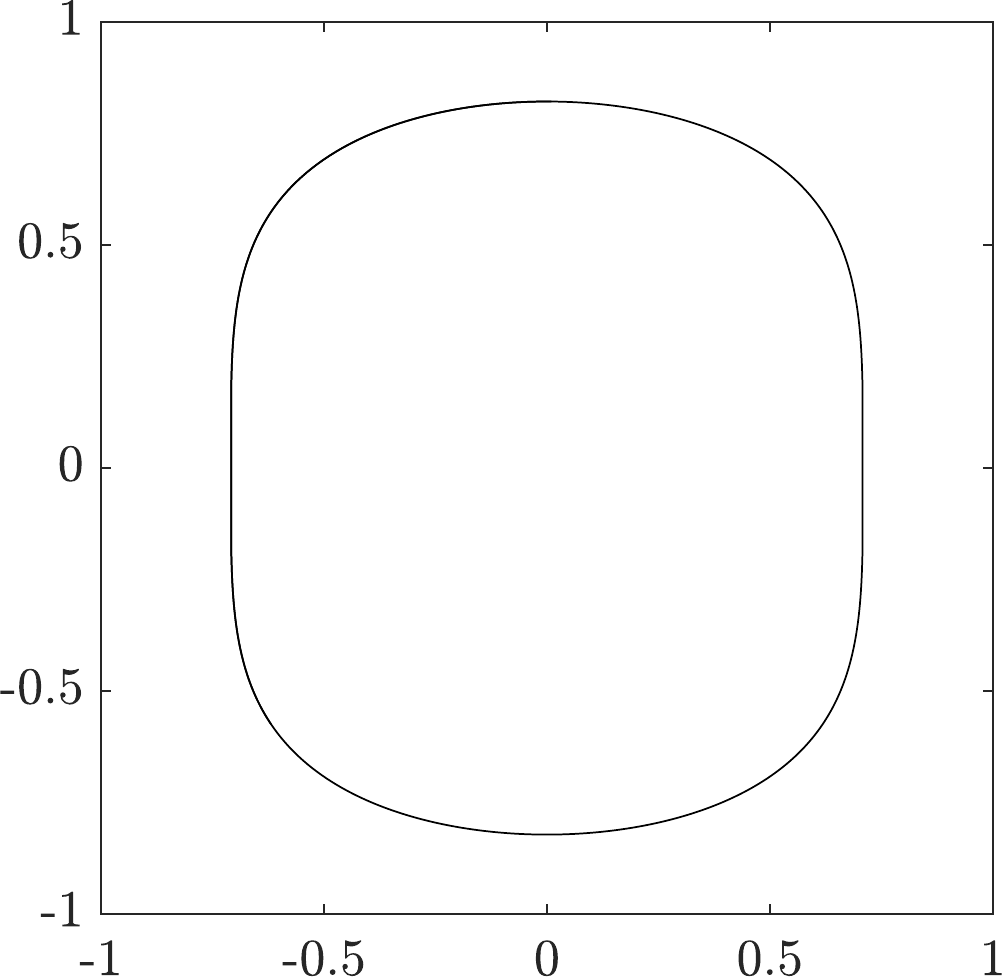}}
	\subfigure[SP]{\includegraphics[width=0.32\textwidth]{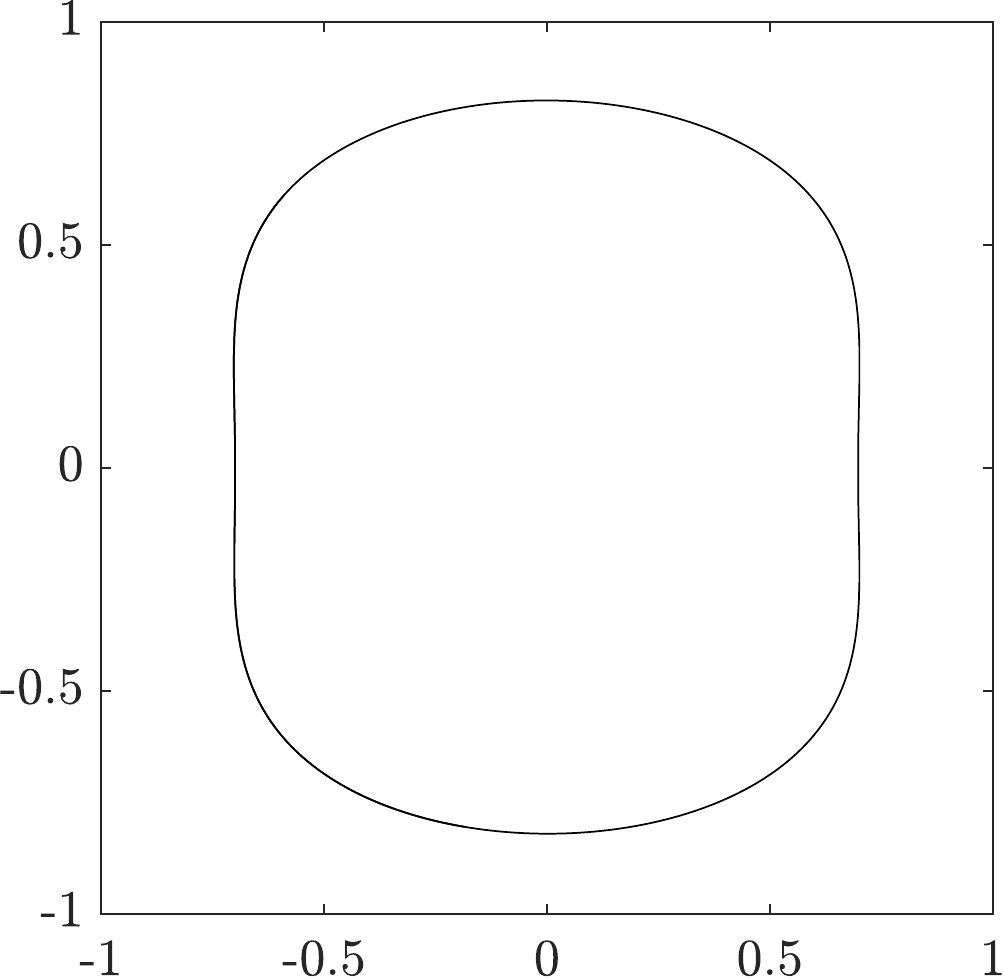}}
	\subfigure[non-SP]{\includegraphics[width=0.32\textwidth]{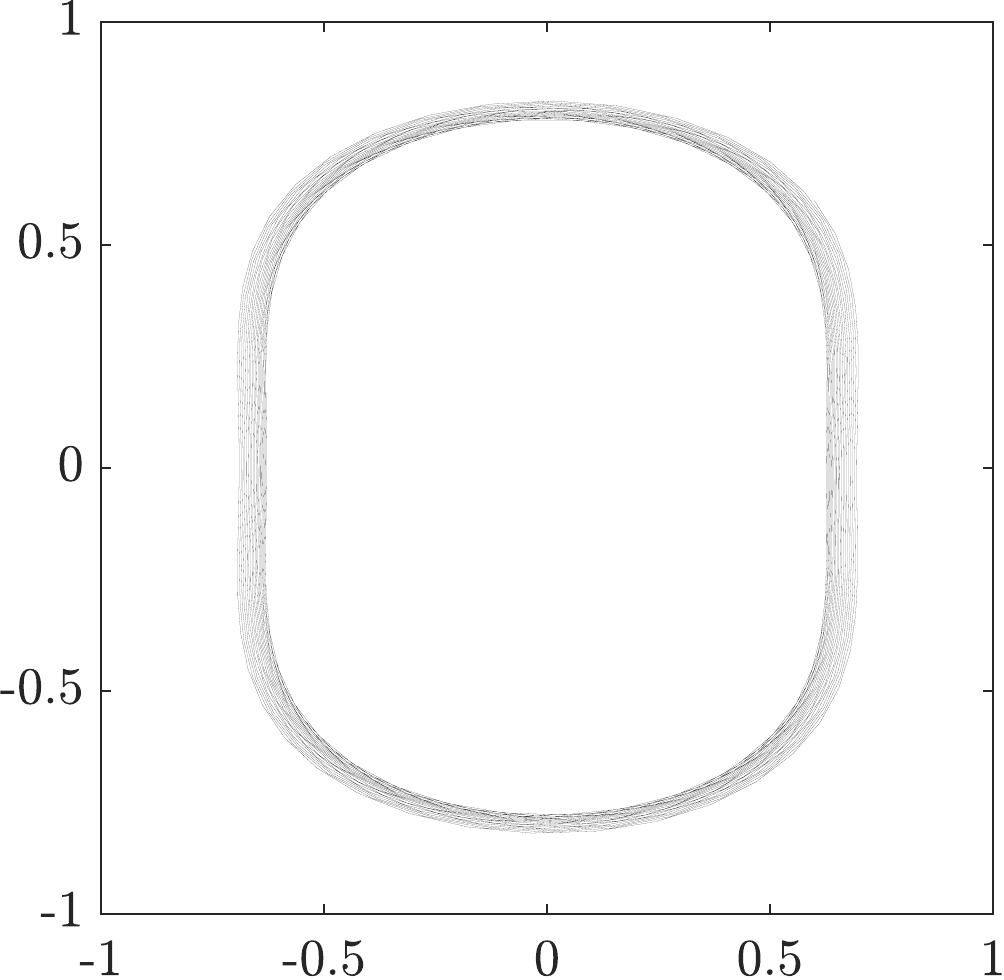}}
	\caption{\small
		Example 2: Phase plots on $p-q$ plane starting from the initial state ${\bf u}_0^*=( 0.6,0.6 )^\top$.
	}\label{fig:ex2_Phase}
\end{figure}

We set the parameters $\alpha_1=1$ and $\alpha_2=1.1$ and  the
computational domain $D$ as $[-1, 1]^2$. We use $M=300$ noiseless
short trajectory data, each of which contains $J=2$ intervals (i.e. 3
data points). The degree of the polynomials for approximation of the
Hamiltonian is $n=6$.  The reconstructed system is solved with
an initial state ${\bf u}_0^*=( 0.6,0.6 )^\top$ and compared against
the solution of the true system. The relative numerical errors in the solutions
of the 
SP algorithm (denoted as $\tu$) and non-SP algorithm from \cite{WuXiu_JCPEQ18} (denoted as $\hu$) are plotted in
\figref{fig:ex2_HErrU}, along with the time evolution of the
reconstructed Hamiltonian. The higher accuracy of the SP algorithm is
again evident from the plot,  as it induces smaller errors over long-term
integration and preserves the approximate Hamiltonian $\widetilde{H}$
along its trajectory.
In \figref{fig:ex2_U} and \figref{fig:ex2_Phase}, we present the
trajectories and the phase plots generated by the reconstructed
system.
The advantage of the new
SP algorithm is again notable, as it is able to preserves both the
phase and amplitude of the solution much better over long-term integration.

\subsubsection*{Example 3: H\'enon-Heiles problem}

\begin{figure}[htbp]
	\centering
	\subfigure[Evolution of relative errors]{\includegraphics[width=0.47\textwidth]{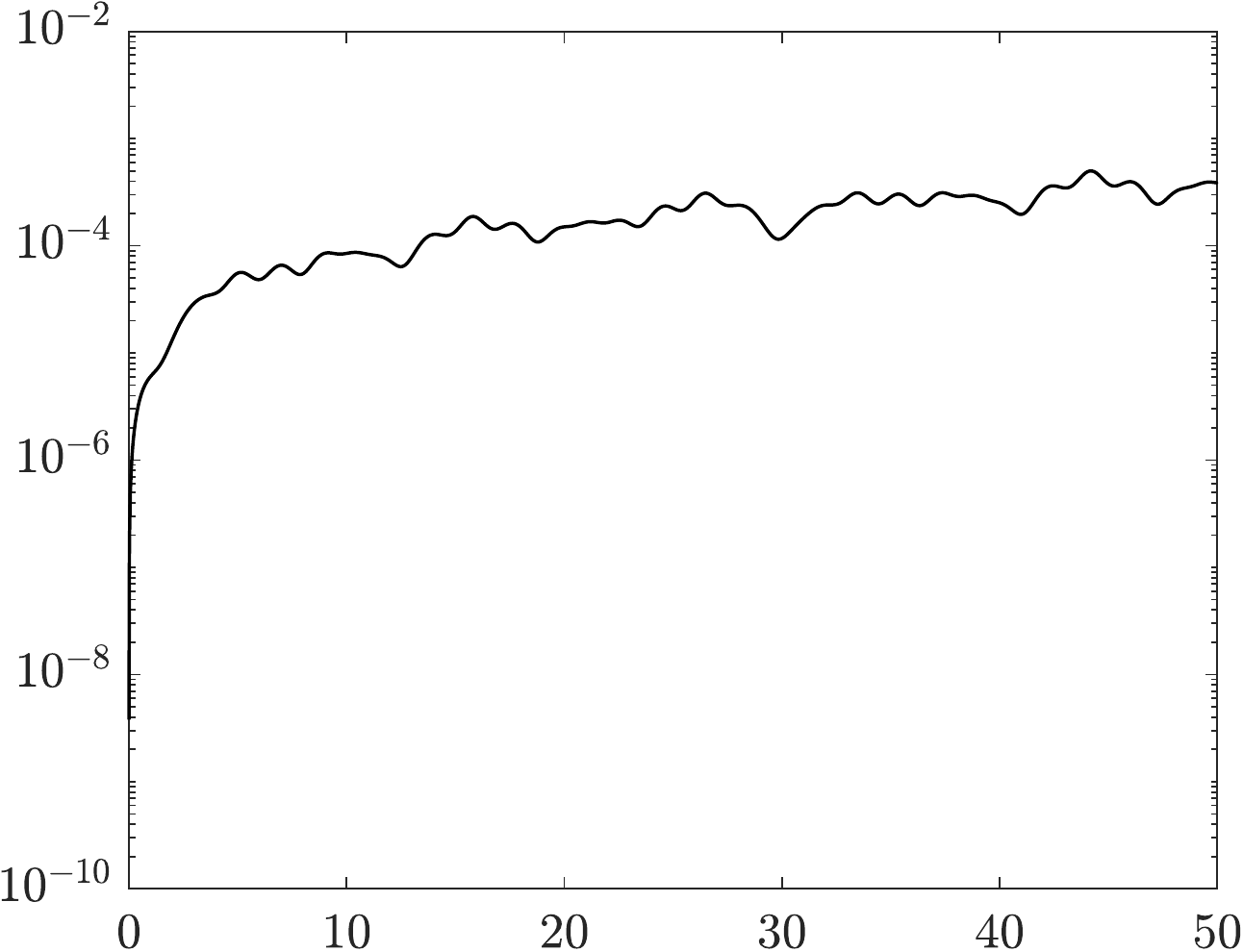}}
	\subfigure[Evolution of Hamiltonian]{\includegraphics[width=0.49\textwidth]{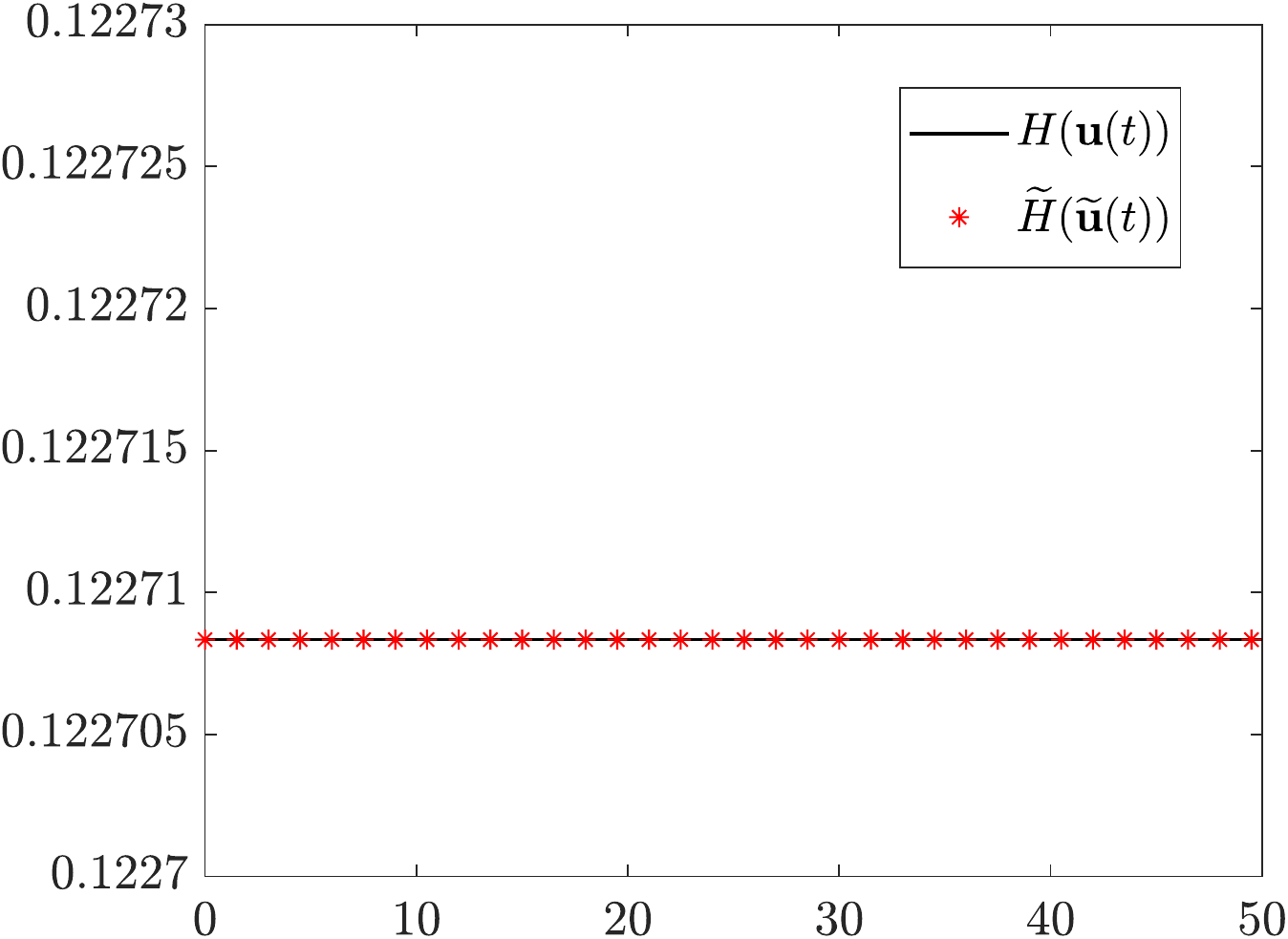}}
	\caption{\small
		Example 3: Solutions of the reconstructed systems with
                an initial state ${\bf u}_0^*=( 0.3,-0.25,0.2,-0.25 )^\top$. Left:
		time evolution of the relative errors; Right: time
		evolution of the Hamiltonian.
	}\label{fig:ex3_HErrU}
\end{figure}

We now consider the H\'enon-Heiles system \cite{henon1964applicability},
\begin{equation}
	\label{eq:example3}
	\begin{cases}
		\dot {p}_1 = -q_1 - 2 q_1 q_2,
		\\
		\dot {p}_2 = -q_2 - q_1^2 + q_2^2,
		\\
		\dot {q}_1 = p_1,\\
		\dot {q}_2 = p_2,
	\end{cases}
\end{equation}
where the Hamiltonian is
$$
H(p_1,p_2,q_1,q_2) =  \frac12 \big( p_1^2 + p_2^2 \big) + \frac12 \big( q_1^2 + q_2^2 \big) + q_1^2 q_2 - \frac13 q_2^3.
$$
This system is used to describe the motion of stars around a galactic
center. Chaotic behavior of the solution will appear when the
Hamiltonian is larger than $1/8$ (\cite{henon1964applicability}).
For  our numerical tests, we choose the computational domain $D$ to be
$[-1, 1]^4$ and employ $M=500$ trajectories, each of which contain
$J=2$ intervals. Polynomials of degree up to $n=3$ are used to
approximate the Hamiltonian. The reconstructed system is solved with
an initial state ${\bf u}_0^*=( 0.3,-0.25,0.2,-0.25 )^\top$ and
compared against the true solution. The time evolution of the
reconstructed Hamiltonian and the
numerical error in the solution are plotted in
\figref{fig:ex3_HErrU}. We observe sufficiently small and stable
numerical errors and good conservation of the Hamiltonian over
relatively long-term integration.

\begin{figure}[htbp]
	\centering
	\subfigure[$p_1(t)$]{\includegraphics[width=0.48\textwidth]{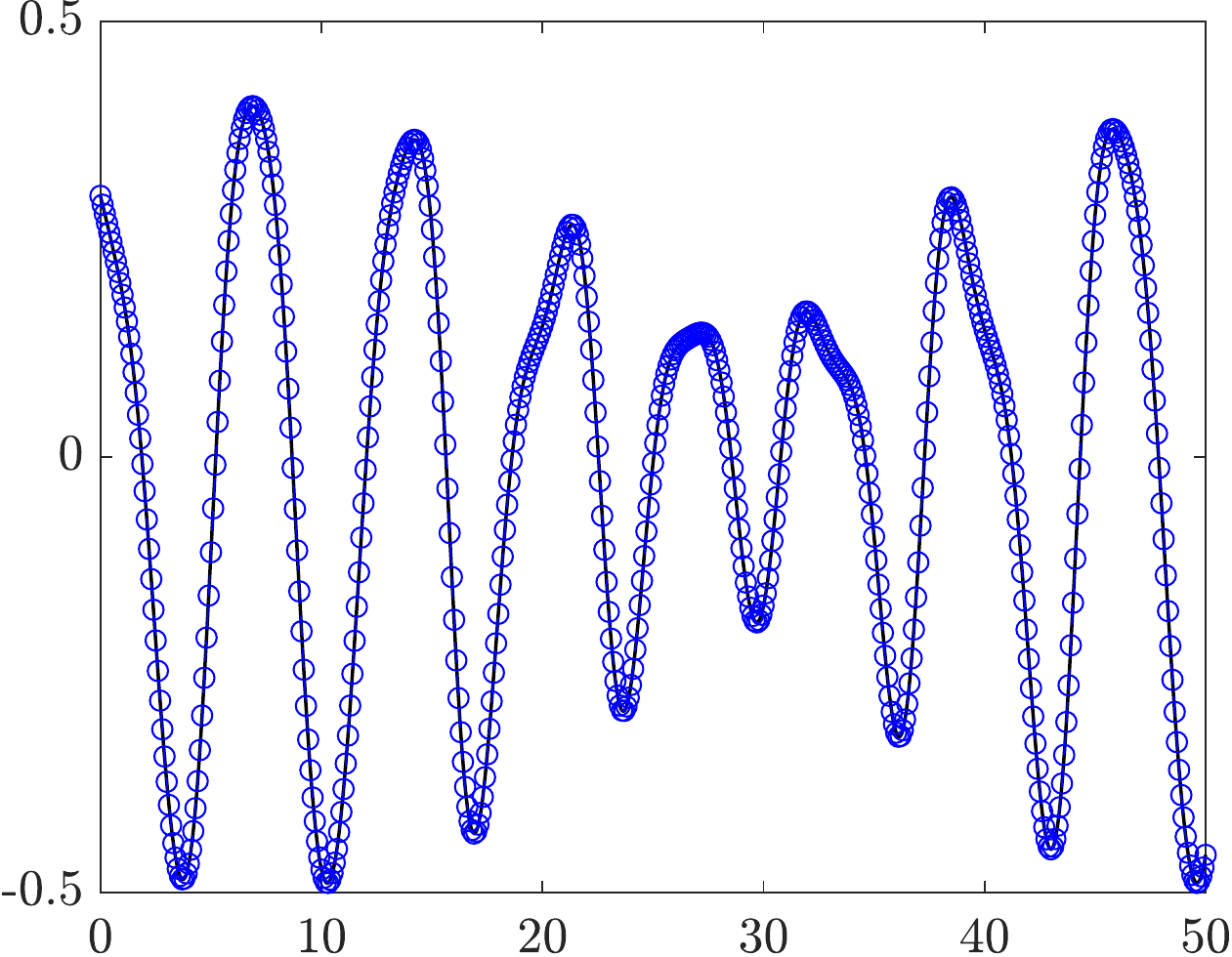}}
	\subfigure[$p_2(t)$]{\includegraphics[width=0.48\textwidth]{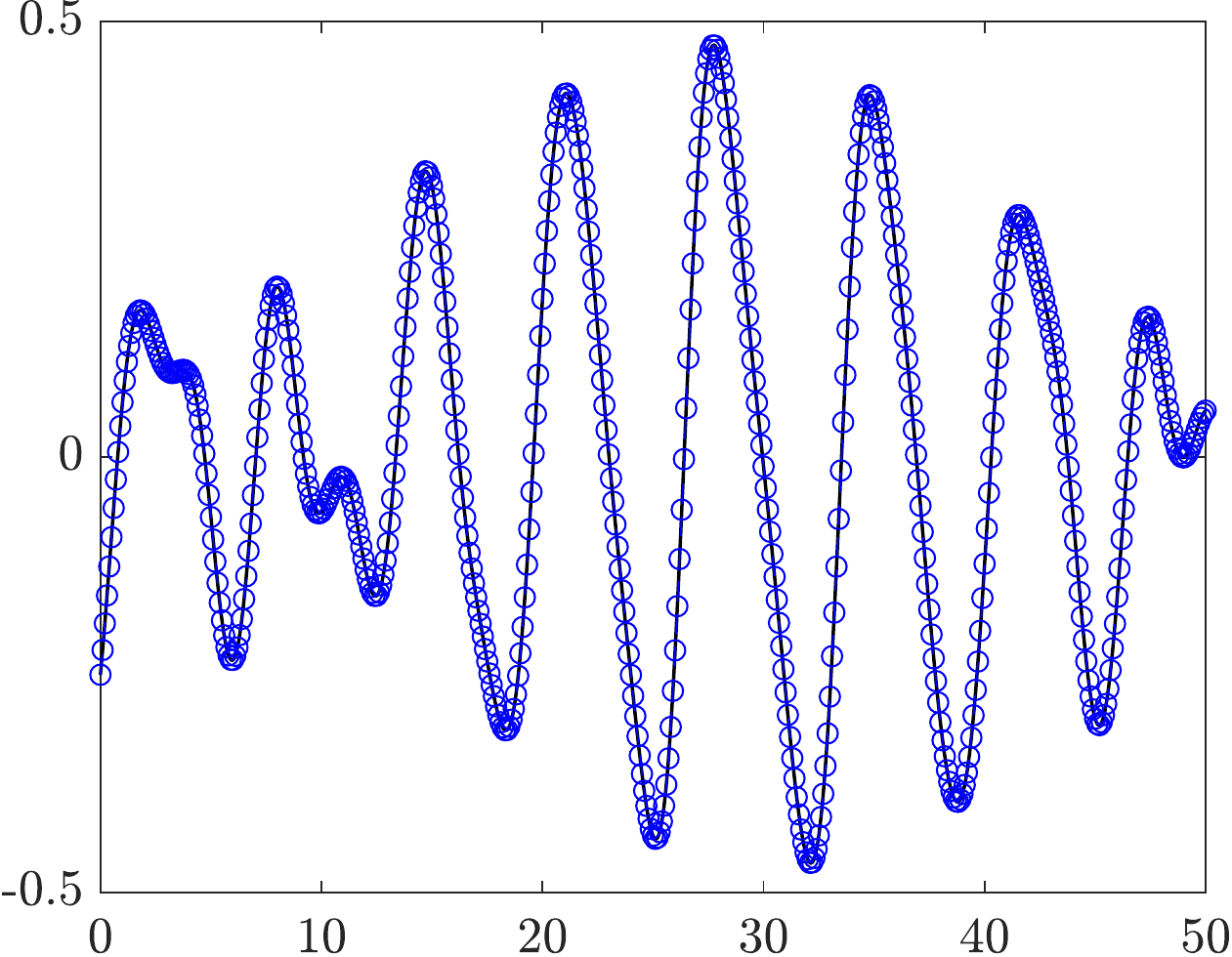}}
	\subfigure[$q_1(t)$]{\includegraphics[width=0.48\textwidth]{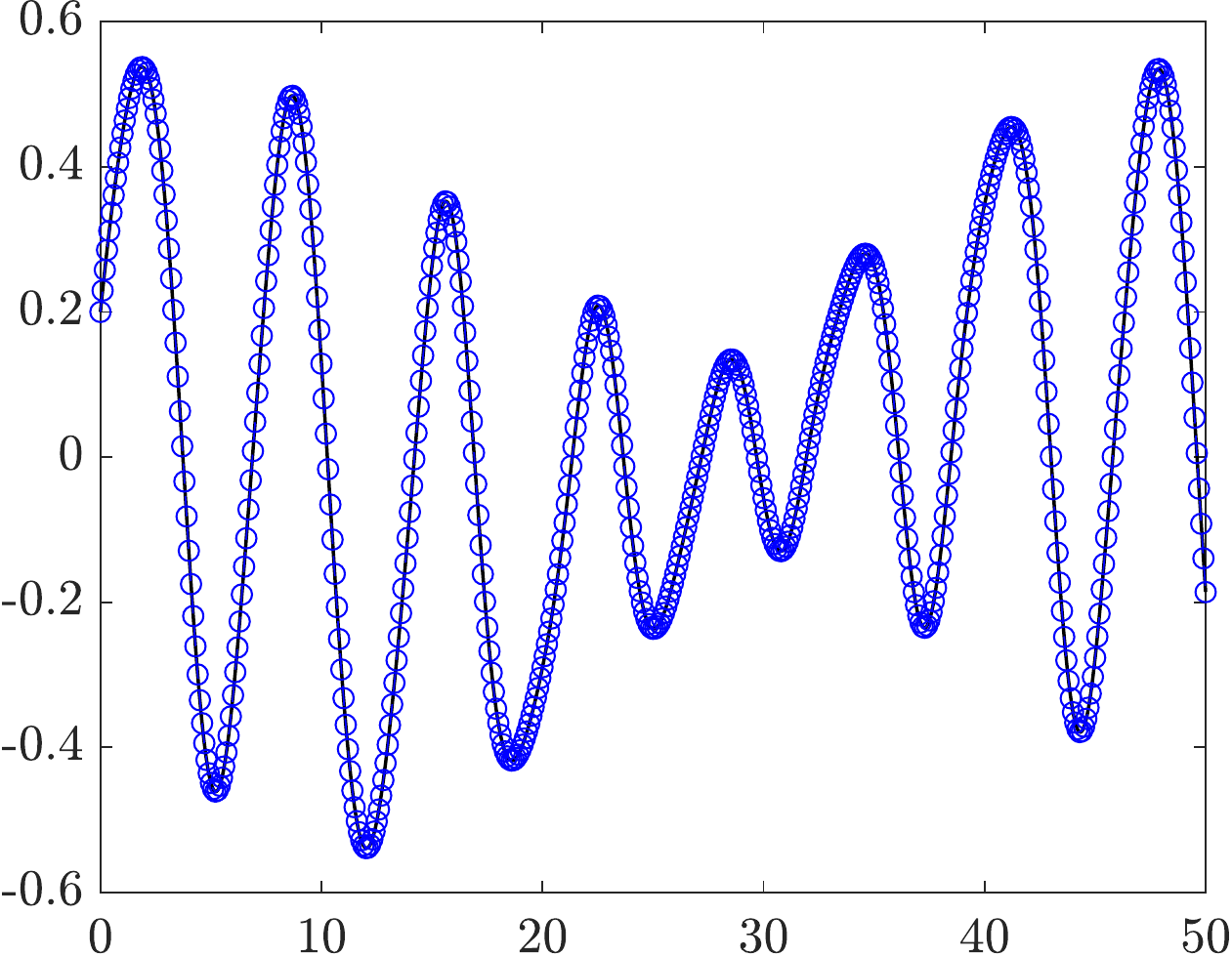}}
	\subfigure[$q_2(t)$]{\includegraphics[width=0.48\textwidth]{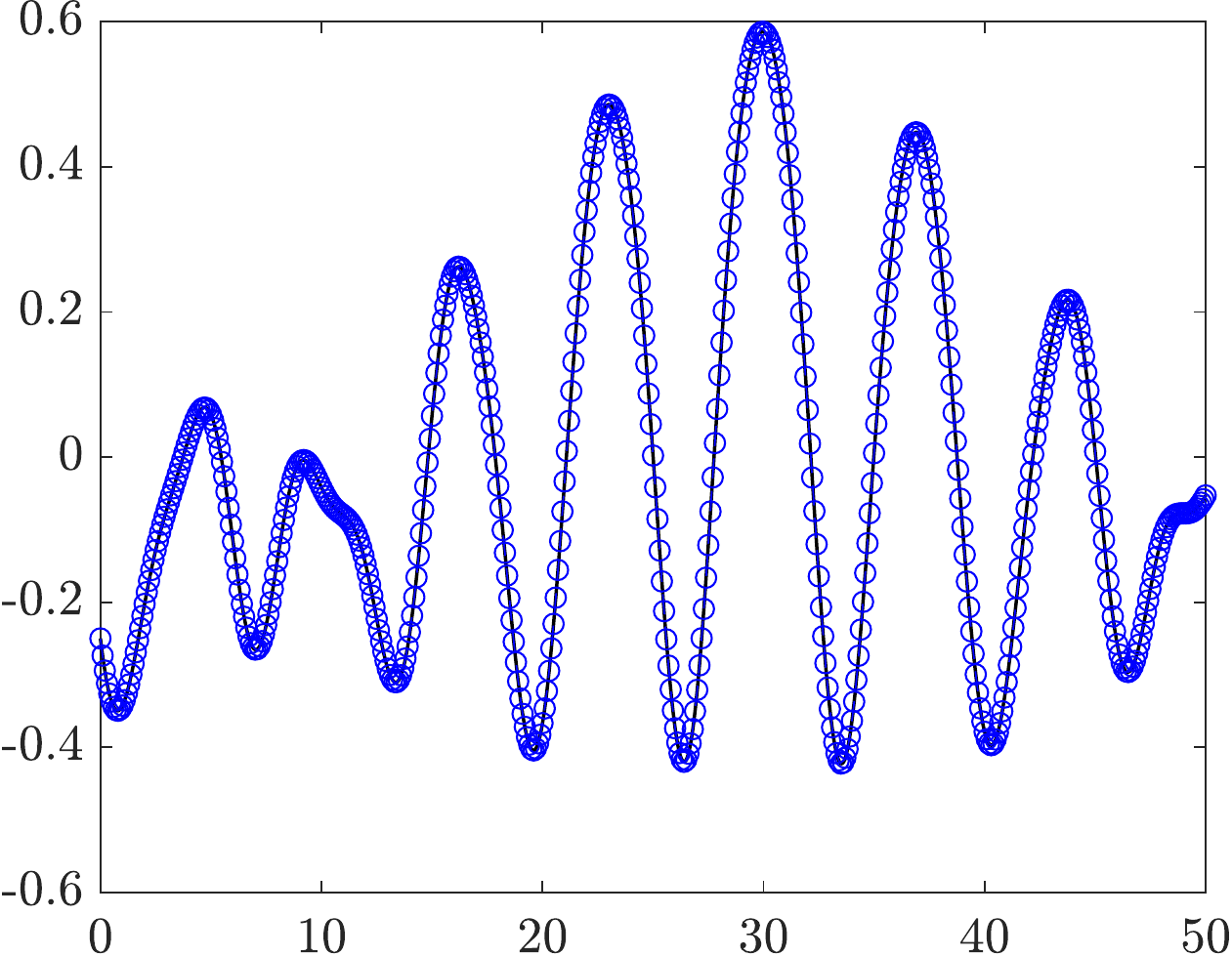}}
	\caption{\small
		Example 3: Comparison of the solutions of the
                reconstructed system (circles) with those of the 
		true system (solid lines), for the same initial state ${\bf u}_0^*=( 0.3,-0.25,0.2,-0.25 )^\top$. 
	}\label{fig:ex3_U} 
\end{figure}

\begin{figure}[htbp]
	\centering
	\subfigure[Plase plots on $q_1$-$p_1$ plane]{\includegraphics[width=0.48\textwidth]{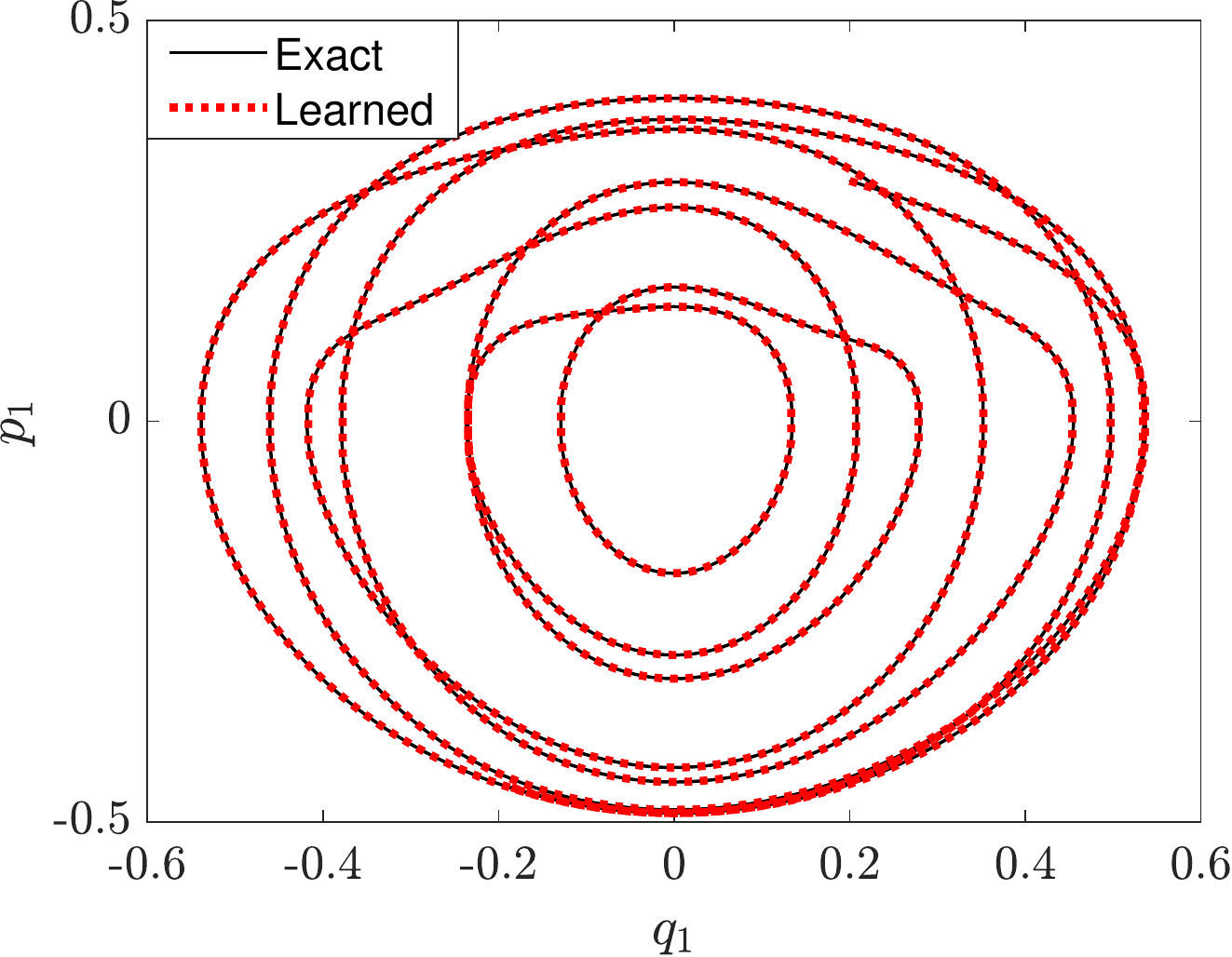}}
	\subfigure[Plase plots on $q_2$-$p_2$ plane]{\includegraphics[width=0.48\textwidth]{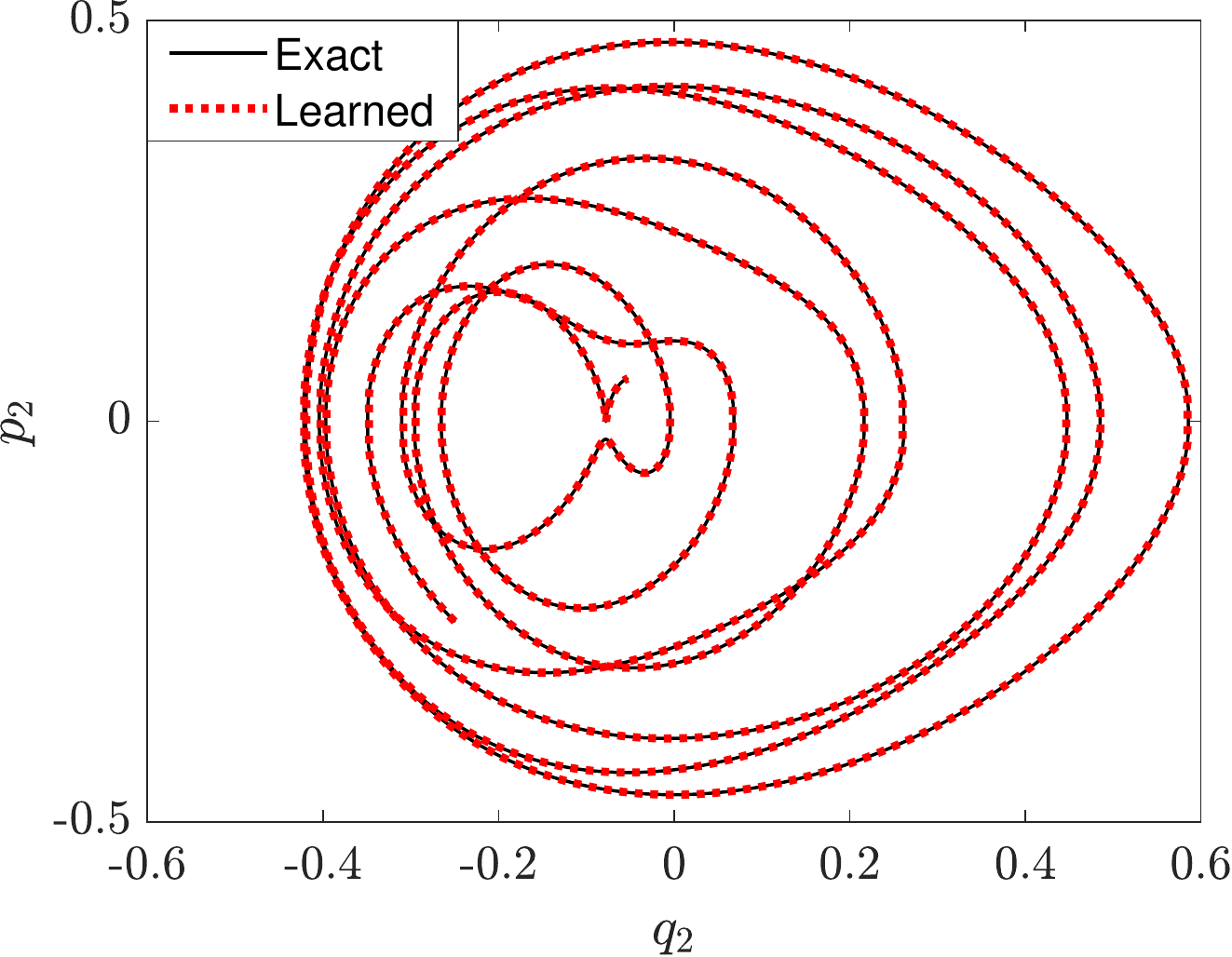}}
	\caption{\small
		Example 3: Phase plots starting from the initial state ${\bf u}_0^*=( 0.3,-0.25,0.2,-0.25 )^\top$.
	}\label{fig:ex3_Phase}
\end{figure}

In \figref{fig:ex3_U} and \figref{fig:ex3_Phase}, the trajectory plots
and phase plots for the reconstructed system using the new SP
algorithm are presented, along with those from the  true system as
reference. The solutions exhibit non-trivial behavior. And the
reconstructed system is able to accurately produce the solutions.

\subsubsection*{Example 4: Cherry problem}

We now consider the Cherry Hamiltonian system \cite{cherry1928v},
\begin{equation}
\label{eq:example4}
\begin{cases}
\dot {p}_1 = -q_1+p_2 q_1+ q_2 p_1,\\
\dot {p}_2 = 2 q_2+q_1 p_1,\\
\dot {q}_1 = p_1+p_2 p_1-q_1 q_2,  \\
\dot {q}_2 = -2 p_2+ \frac12(p_1^2-q_1^2),
\end{cases}
\end{equation}
whose  true Hamiltonian is
$$
H(p_1,p_2,q_1,q_2) =  \frac12 (q_1^2 + p_1^2) - (q_2^2 + p_2^2) +\frac12 p_2(p_1^2-q_1^2)-q_1q_2p_1.
$$

\begin{figure}[htbp]
	\centering
	\subfigure[Evolution of relative errors]{\includegraphics[width=0.48\textwidth]{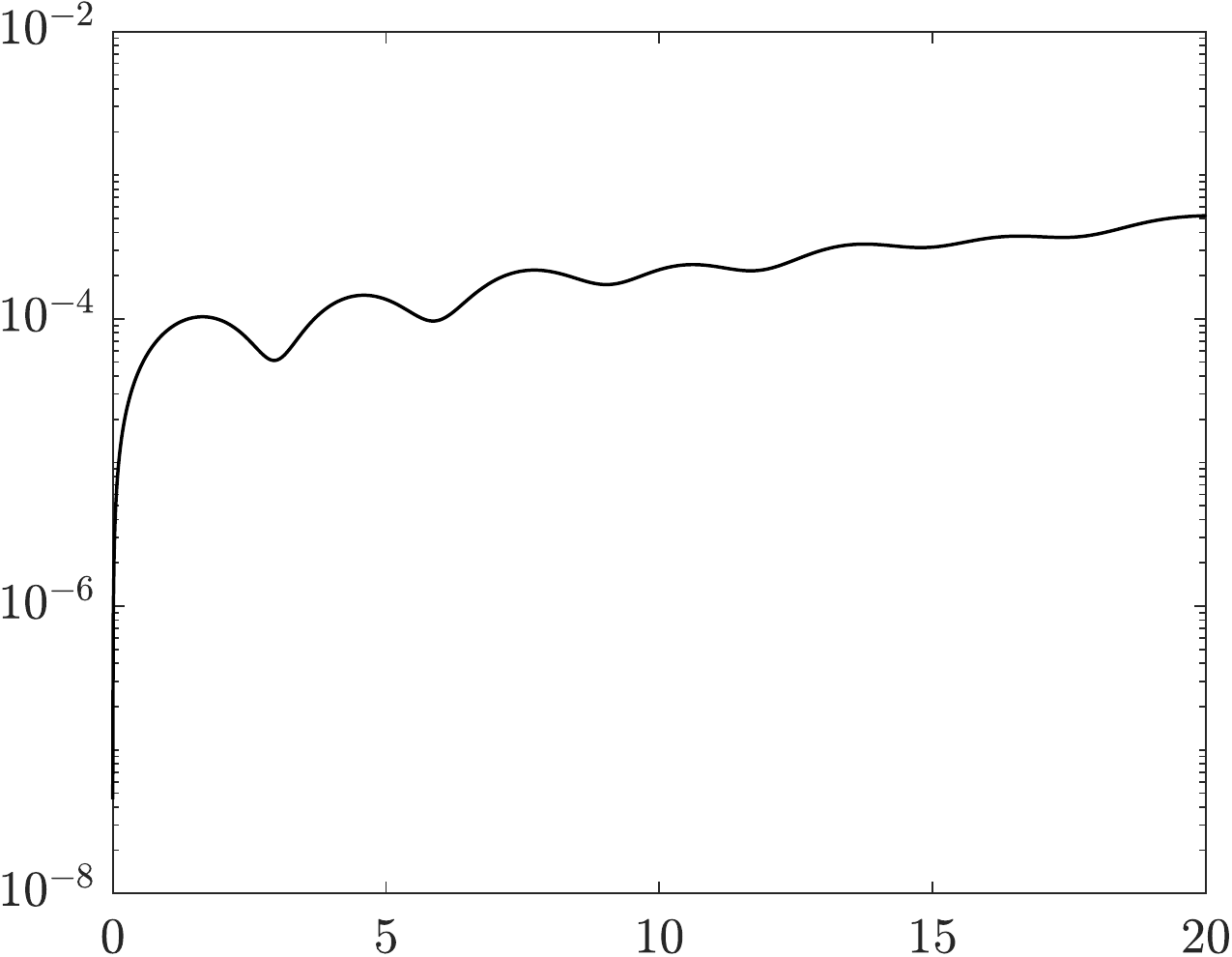}}
	\subfigure[Evolution of Hamiltonian]{\includegraphics[width=0.48\textwidth]{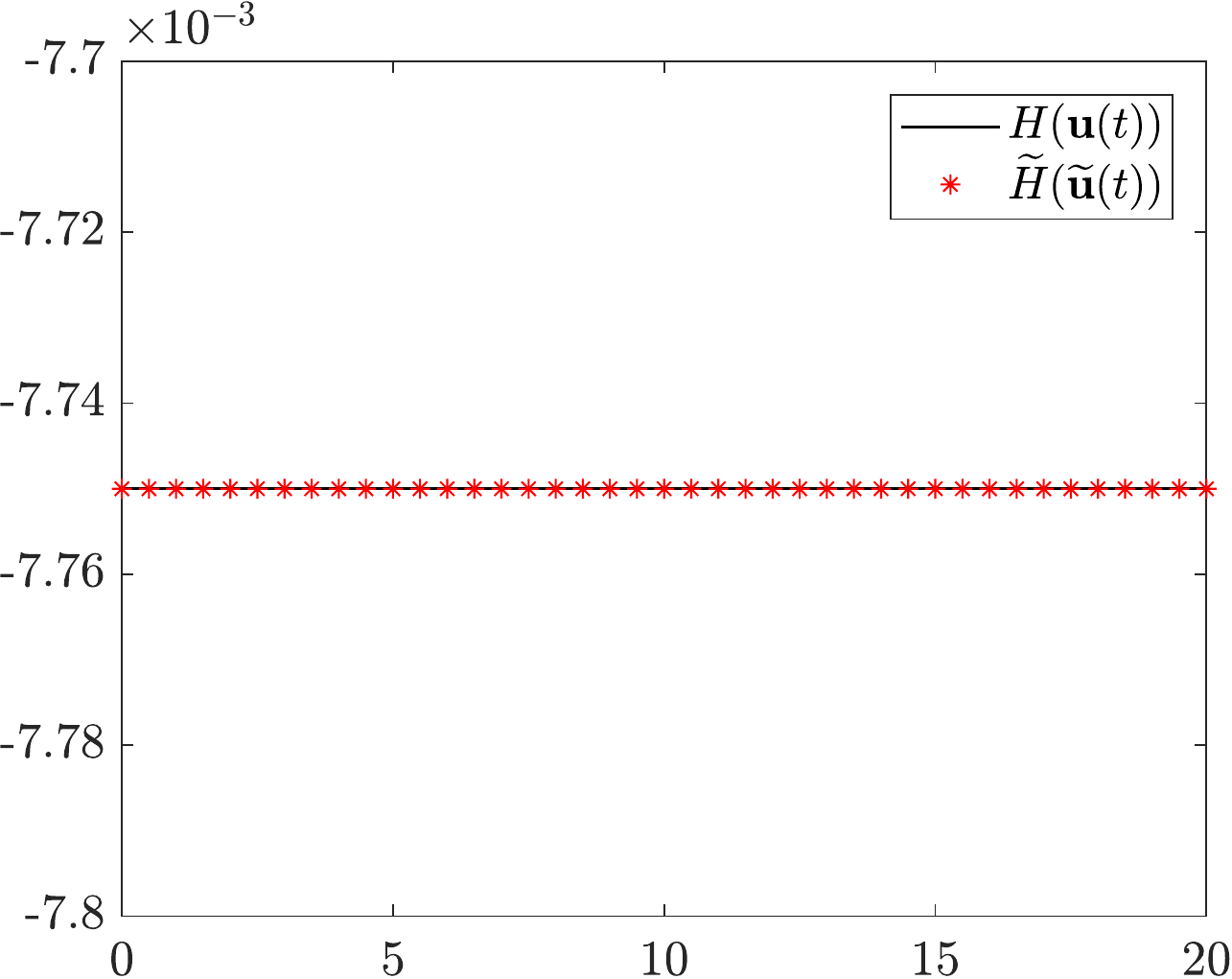}}
	\caption{\small
		Example 4: Solutions of the reconstructed system with an initial state ${\bf u}_0^*=( -0.05, 0.1, 0.15, 0.1 )^\top$.  Left:
		time evolution of the relative error; Right: time
		evolution of the Hamiltonian.
	}\label{fig:ex4_HErrU}
\end{figure}

\begin{figure}[htbp]
	\centering
	\subfigure[$p_1(t)$]{\includegraphics[width=0.48\textwidth]{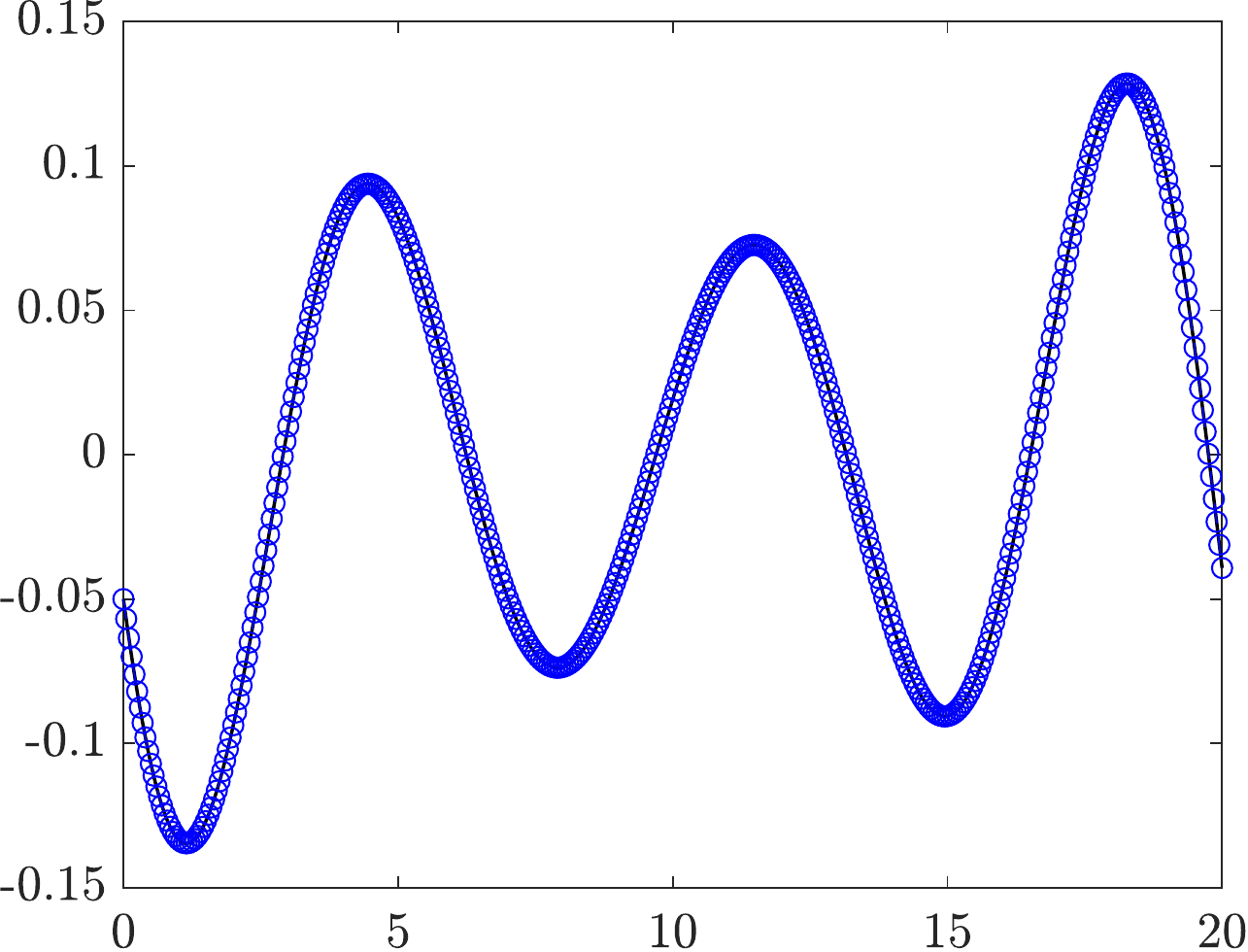}}
	\subfigure[$p_2(t)$]{\includegraphics[width=0.48\textwidth]{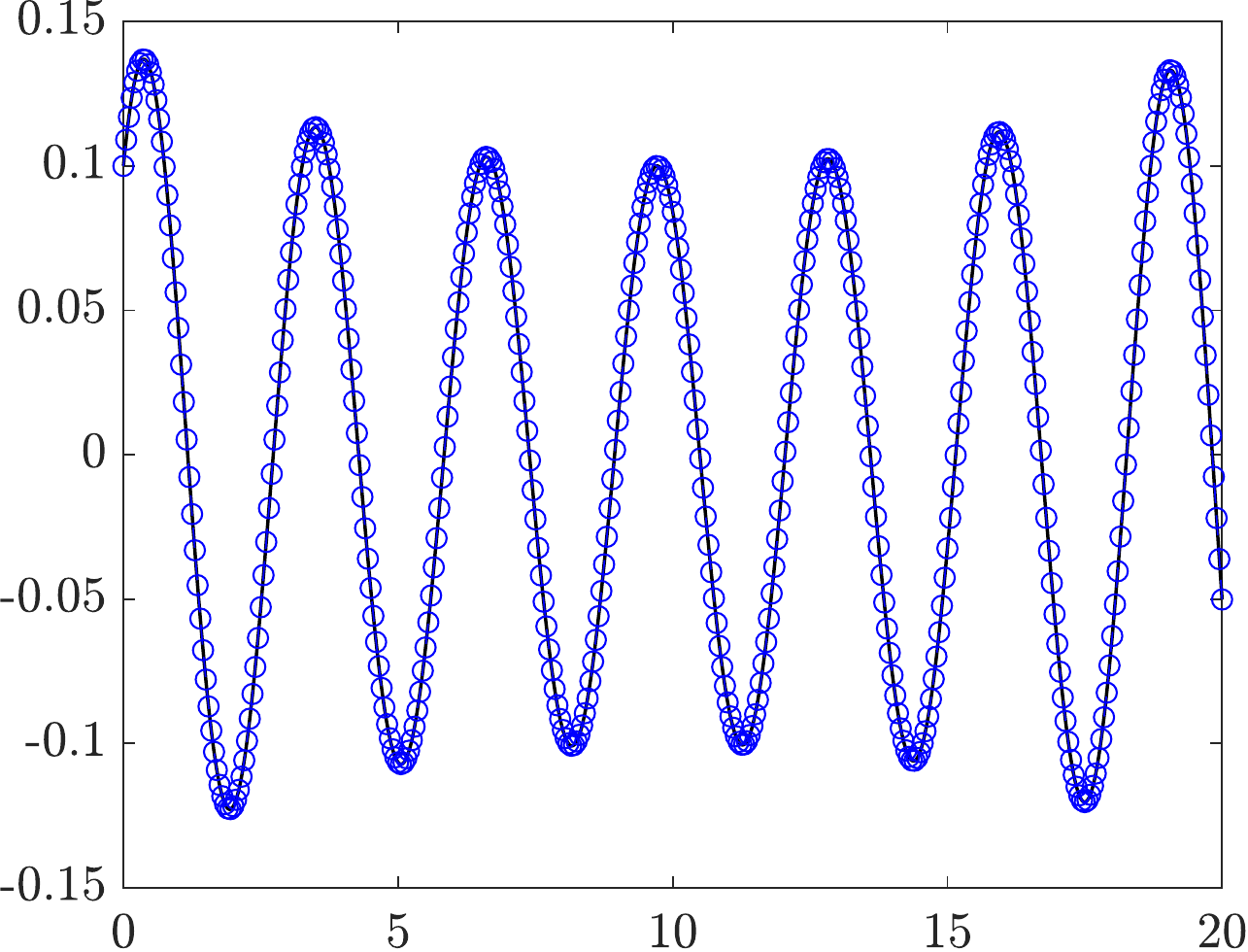}}
	\subfigure[$q_1(t)$]{\includegraphics[width=0.48\textwidth]{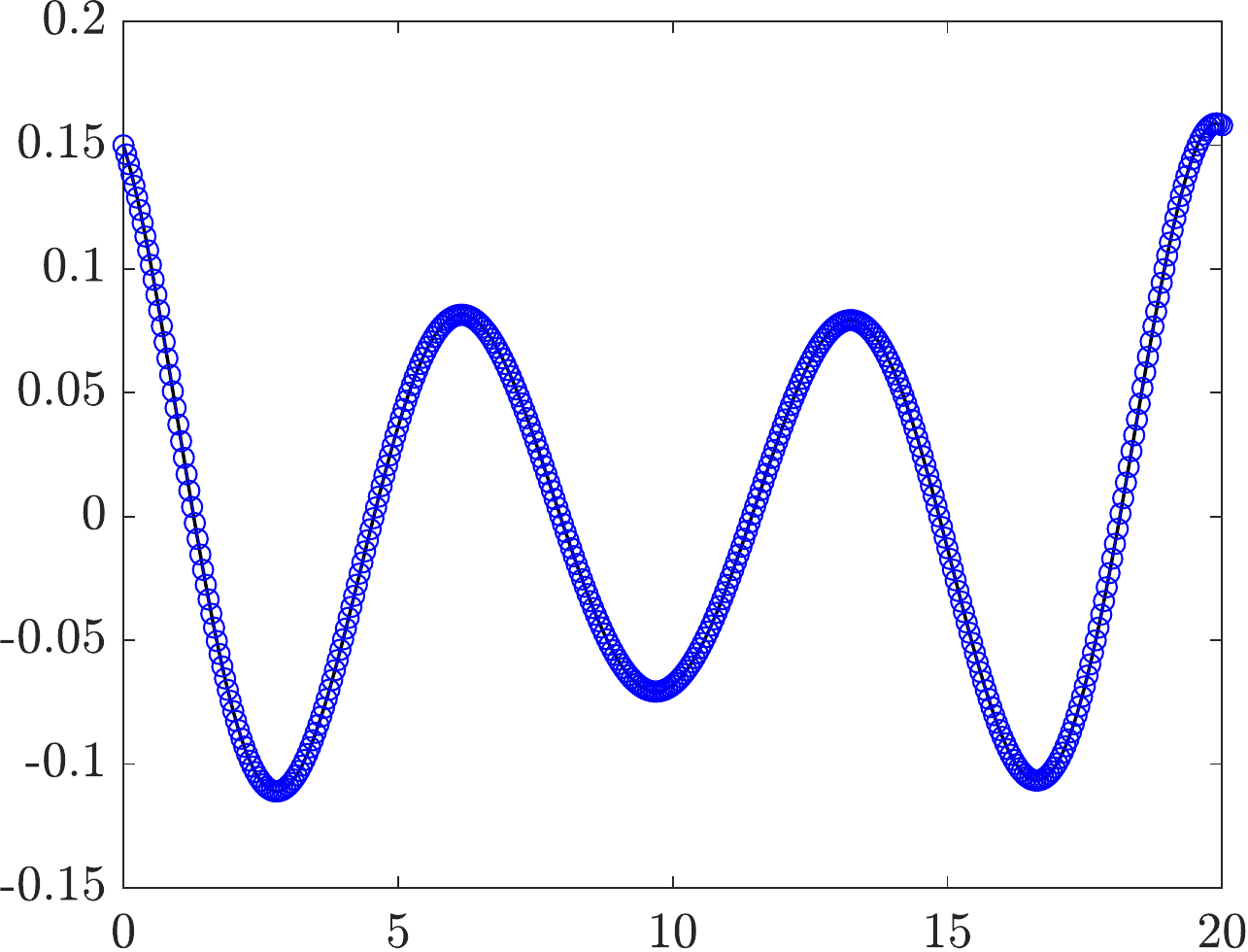}}
	\subfigure[$q_2(t)$]{\includegraphics[width=0.48\textwidth]{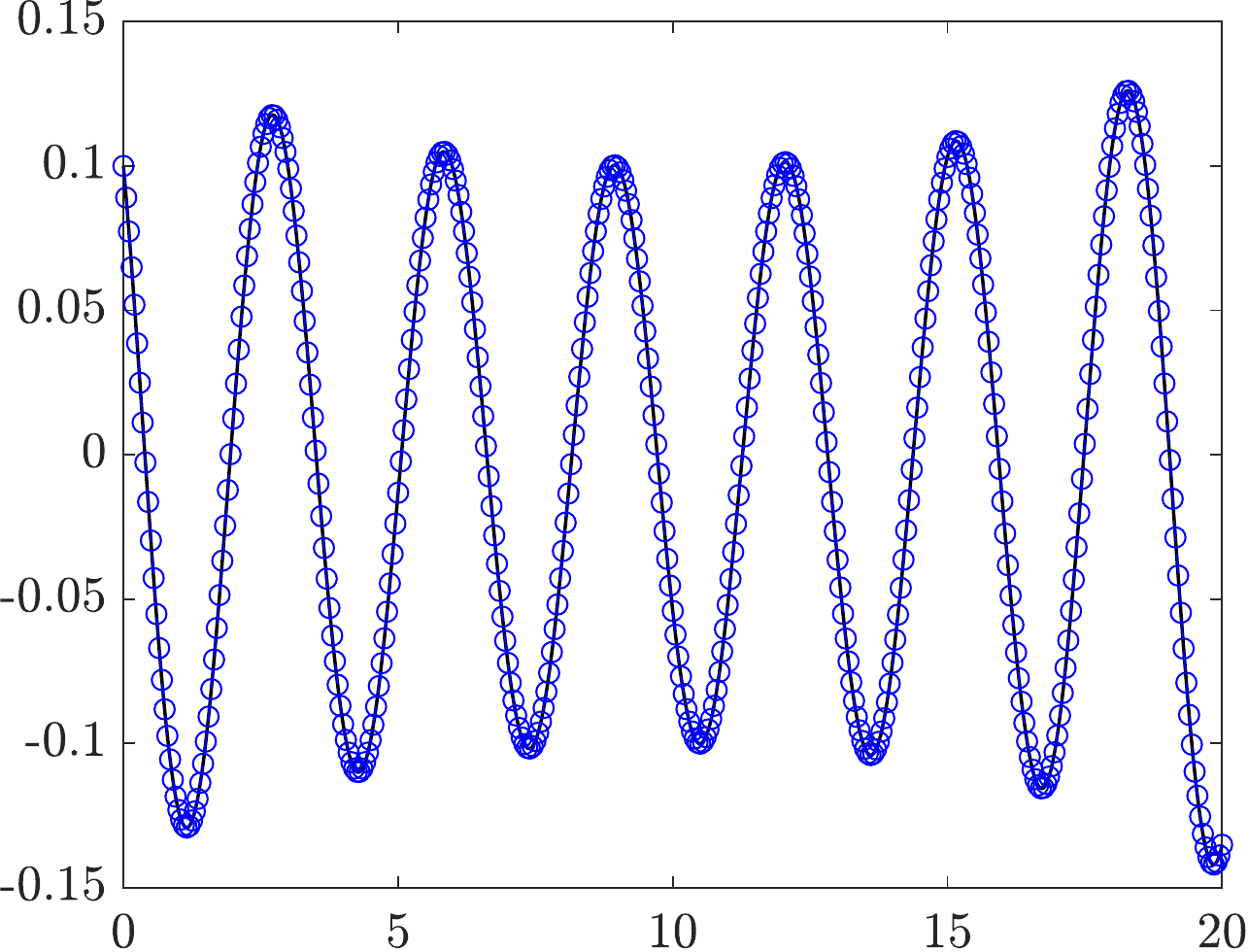}}
	\caption{\small
		Example 4: Comparison of solutions of the reconstructed
                systems (circles)  with the solution of 
		the true system (solid lines), with the initial state ${\bf u}_0^*=( -0.05, 0.1, 0.15, 0.1 )^\top$. 
	}\label{fig:ex4_U}
\end{figure}

\begin{figure}[htbp]
	\centering
	\subfigure[Plase plots on $p_1$-$p_2$ plane]{\includegraphics[width=0.48\textwidth]{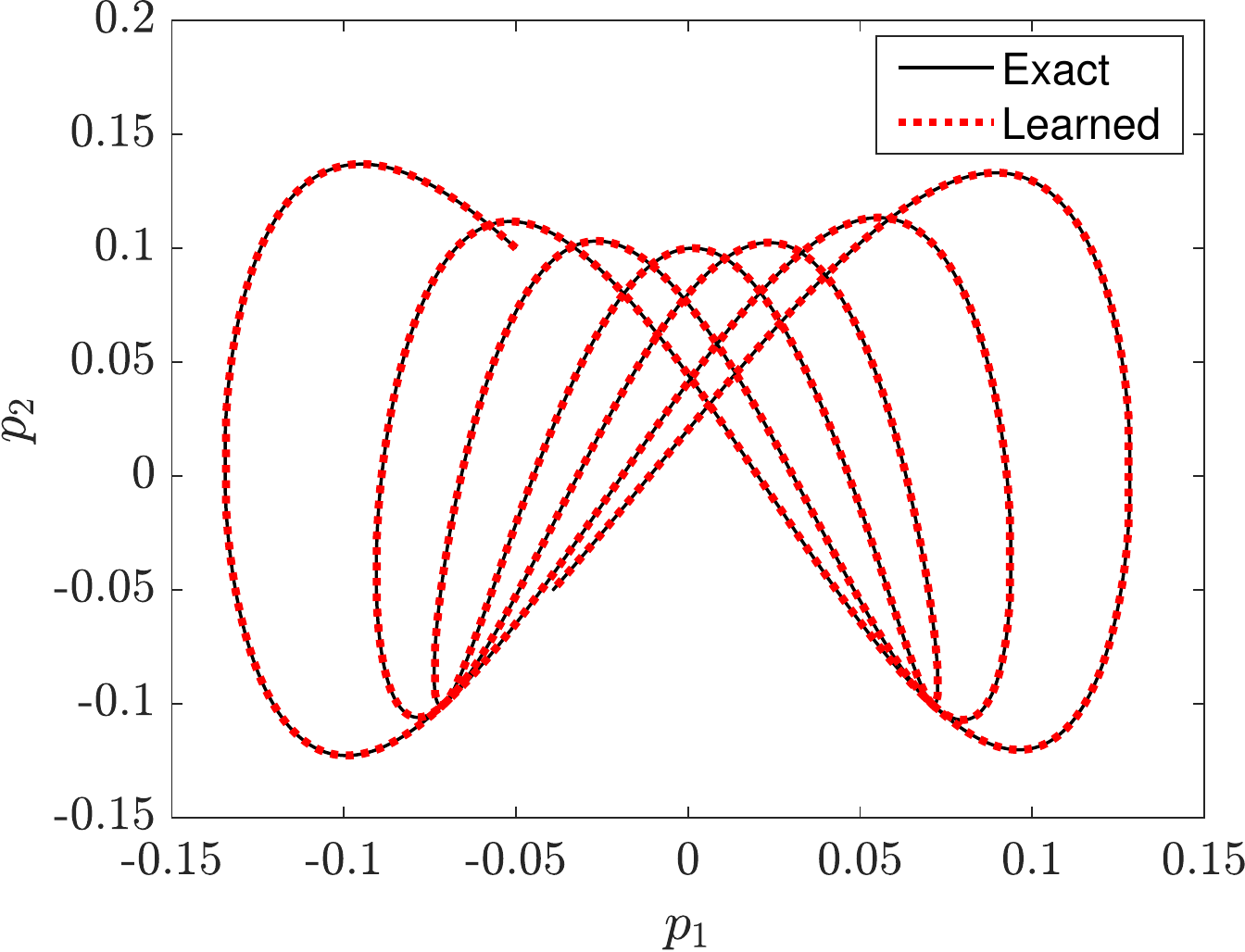}}
	\subfigure[Plase plots on $q_1$-$q_2$ plane]{\includegraphics[width=0.48\textwidth]{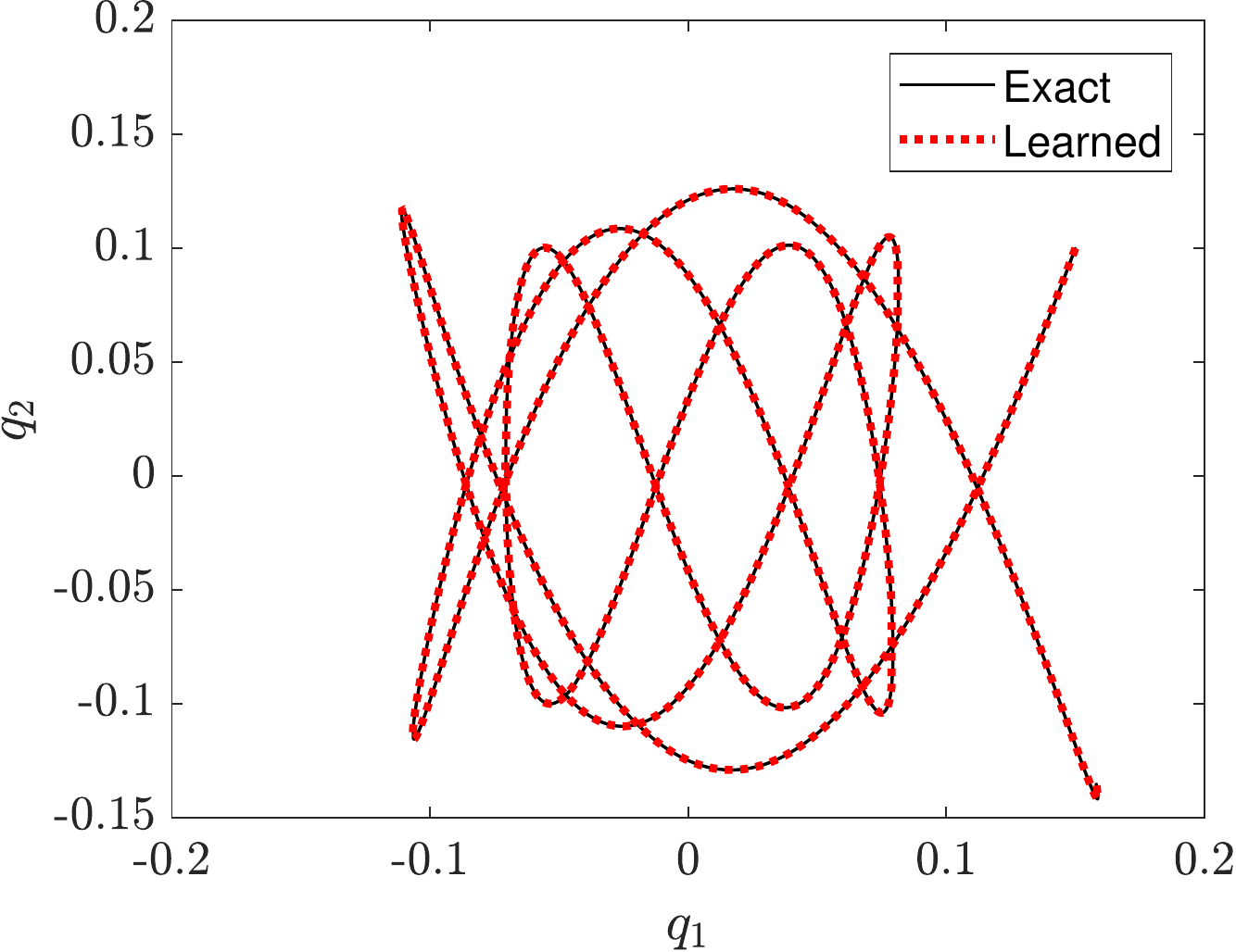}}
	\caption{\small
		Example 4: Phase plots starting from the initial state ${\bf u}_0^*=( -0.05, 0.1, 0.15, 0.1 )^\top$.
	}\label{fig:ex4_Phase}
\end{figure}

We take the computational domain as $D= (-2,2)\times (-1,2) \times
(-2,1) \times (-1,1)$ and use $M=500$ short trajectory data, each of
which contain $J=2$ intervals. Polynomials of degree up to $n=3$ are
employed to approximate the Hamiltonian. The reconstructed system is
then solved using an arbitrarily chosen initial state ${\bf u}_0^*=(
-0.05, 0.1, 0.15, 0.1 )^\top$. The solutions are then compared against
those from the true system with the same initial state. The relative errors in
the numerical prediction are plotted in \figref{fig:ex4_HErrU}, along
with the time evolution of the reconstructed Hamiltonian
$\widetilde{H}$ and the true Hamiltonian $H$. We again observe good
accuracy by the SP algorithm and conservation of the approximate
Hamiltonian.
The solution states are plotted in Figs.\ \ref{fig:ex4_U}, and their
phase plots in \ref{fig:ex4_Phase}. The numerical solutions agree with
the true solutions well.

\subsubsection*{Example 5: Double pendulum}

\begin{figure}[h]
	\centering
	\includegraphics[width=0.3\textwidth]{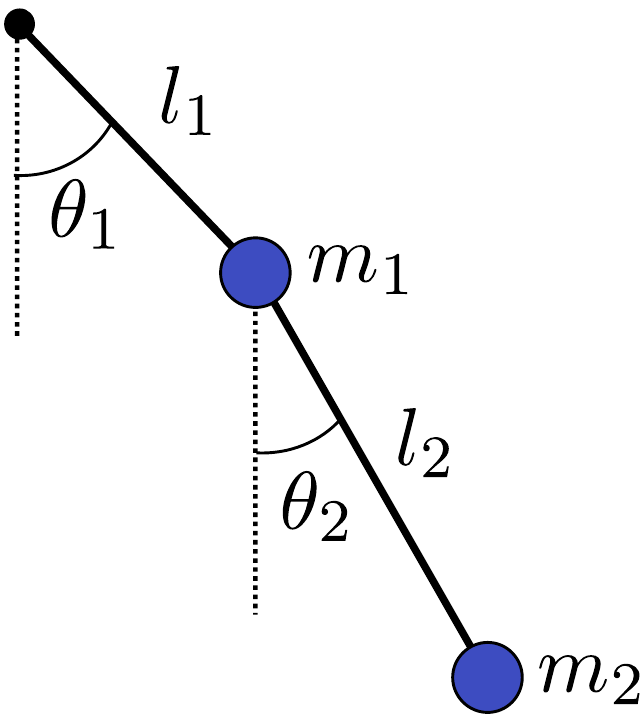}
	\caption{\small
		A diagram of the double pendulum.
	}\label{fig:ex5_DIAG}
      \end{figure}

Finally, we consider a double pendulum problem, as illustrated in \figref{fig:ex5_DIAG}.
      Two masses $m_1$ and $m_2$ are connected via massless rigid rods
      of length $l_1$ and $l_2$, and
$\theta_1$ and $\theta_2$ are
the angles of the two rods with respect to the vertical direction.
We define the canonical momenta of the system as
\begin{align*}
& p_1 := (m_1+m_2)\, l_1^2\, \dot \theta_1 + m_2\, l_1 \,l_2\, \dot \theta_2 \cos(\theta_1 - \theta_2),
\\
& p_2 := m_2\,l_2^2\, \dot \theta_2 + m_2\, l_1\, l_2\, \dot \theta_1 \cos(\theta_1 - \theta_2).
\end{align*}
By letting $q_1 = \theta_1$ and $q_2 = \theta_2$, the Hamiltonian of the system is 
\begin{align*}
H(p_1,p_2,q_1,q_2) & = \frac{ m_2\, l_2^2\, p_1^2 + (m_1 + m_2)\,l_1^2\, p_2^2 - 2m_2\, l_1\, l_2\, p_1\, p_2 
\cos(q_1 - q_2) }{ 2 m_2\, l_1^2\, l_2^2\, \big[ m_1 + m_2\, \sin^2(q_1 - q_2)  \big] }
\\
& \quad  - (m_1 + m_2)\, g\,l_1 \cos q_1 - m_2\, g\, l_2 \cos q_2,
\end{align*}
where $g$ is the gravitational constant. The governing equations of
the system are 
\begin{equation}
\label{eq:example5}
\begin{cases}
\dot {p}_1 = -(m_1 + m_2) g\,l_1 \sin q_1 - C_1 + C_2 \sin( 2 (q_1 - q_2) ),\\
\dot {p}_2 = -m_2\, g\, l_2 \sin q_2 +C_1 -C_2 \sin ( 2(q_1 -q_2) ), \\[2pt]
\dot {q}_1 = \dfrac{ l_2\, p_1 - l_1\, p_2 \cos(q_1 -q_2) }{ l_1^2\, l_2 \big[ m_1 + m_2 \sin^2 (q_1 - q_2)  \big] },  \\[10pt]
\dot {q}_2 =  \dfrac{ -m_2\, l_2\, p_1 \cos( q_1 -q_2 ) + (m_1 + m_2)\, l_1\, p_2 }{ m_2\, l_1\, l_2^2 \big[ m_1 + m_2 \sin^2 (q_1 - q_2)  \big] },
\end{cases}
\end{equation}
where
\begin{align*}
& C_1(\bm p,\bm q) := \frac{ p_1\, p_2 \sin(q_1 - q_2) }{ l_1\, l_2 \big[ m_1 + m_2 \sin^2 ( \theta_1 - \theta_2)  \big] },
\\
& C_2(\bm p,\bm q) := \frac{ m_2\, l_2^2\, p_1^2 + (m_1 + m_2)l_1^2\, p_2^2 - 2\,m_2\, l_1\, l_2\, p_1\, p_2 \cos ( q_1 - q_2) }{ 2l_1^2\, l_2^2\, \big[ m_1 + m_2 \sin^2 ( q_1 - q_2 ) \big] }.
\end{align*}

\begin{figure}[htbp]
	\centering
	{\includegraphics[width=0.48\textwidth]{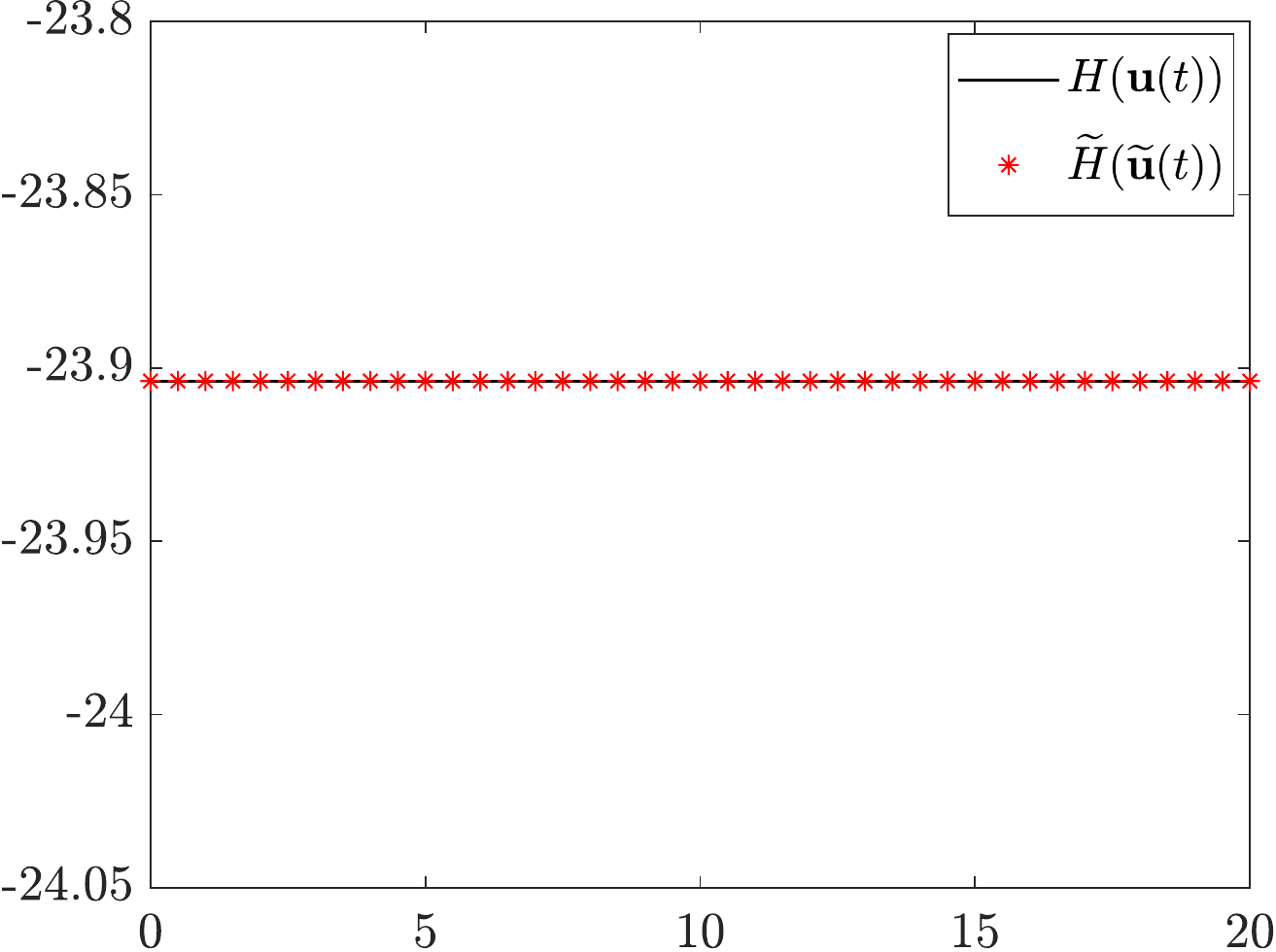}}
	\caption{\small
		Example 5: Time
		evolution of the Hamiltonian of the reconstructed system with
		an initial state ${\bf u}_0^*=( 0,0, \frac{\pi}6, \frac{\pi}4 )^\top$. 
	}\label{fig:ex5_HU}
\end{figure}

\begin{figure}[htbp]
	\centering
	\subfigure[$p_1(t)$]{\includegraphics[width=0.48\textwidth]{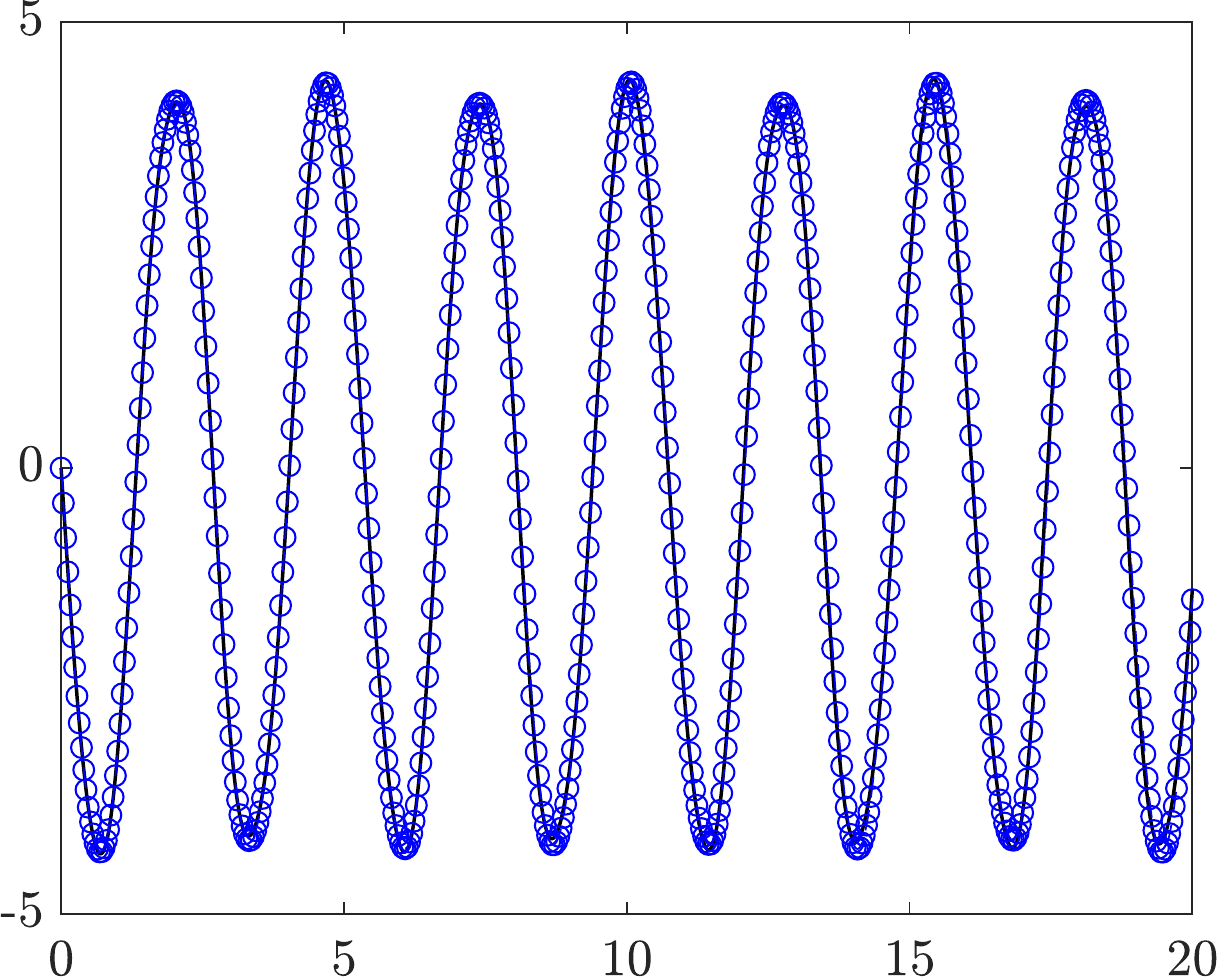}}
	\subfigure[$p_2(t)$]{\includegraphics[width=0.48\textwidth]{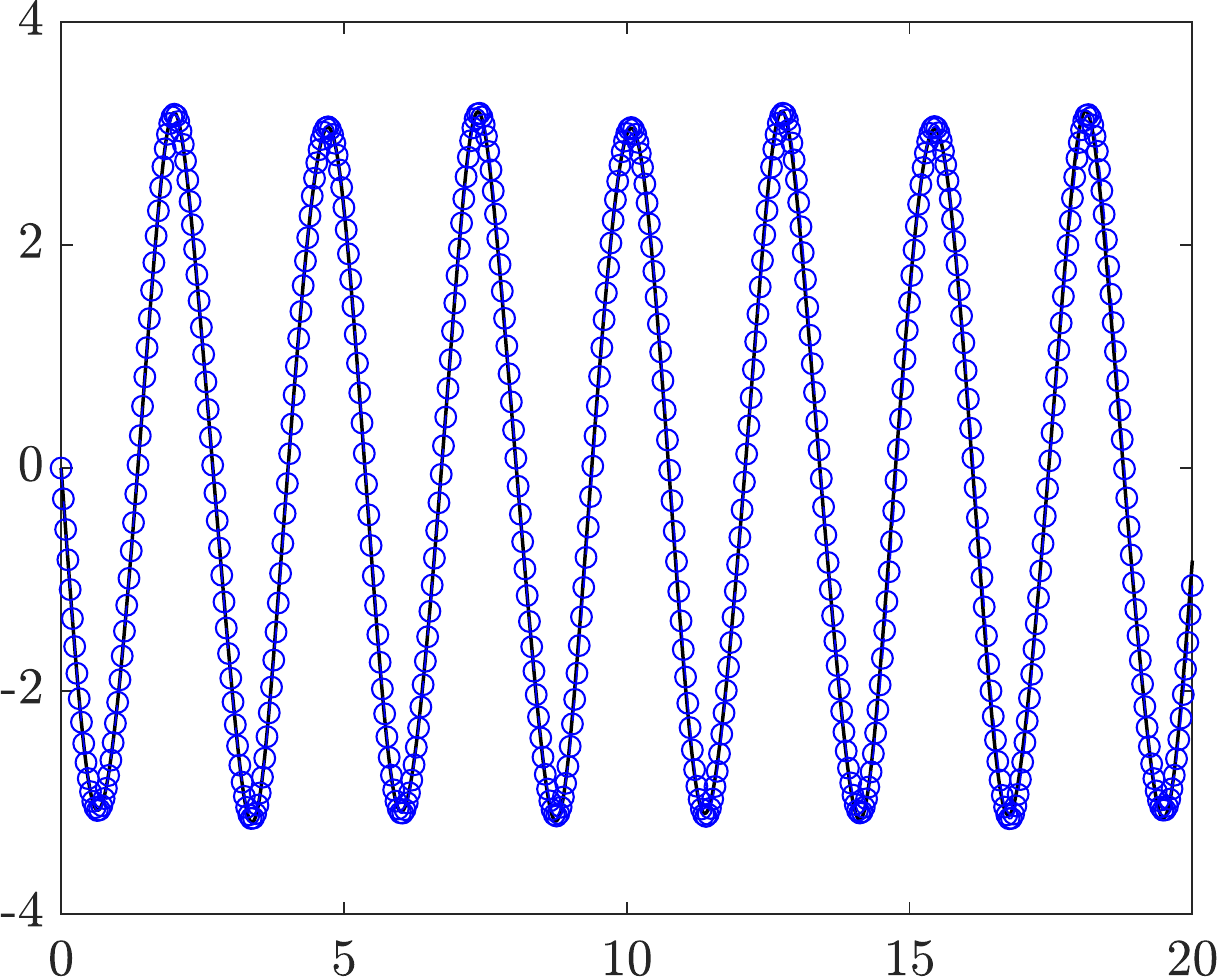}}
	\subfigure[$q_1(t)$]{\includegraphics[width=0.48\textwidth]{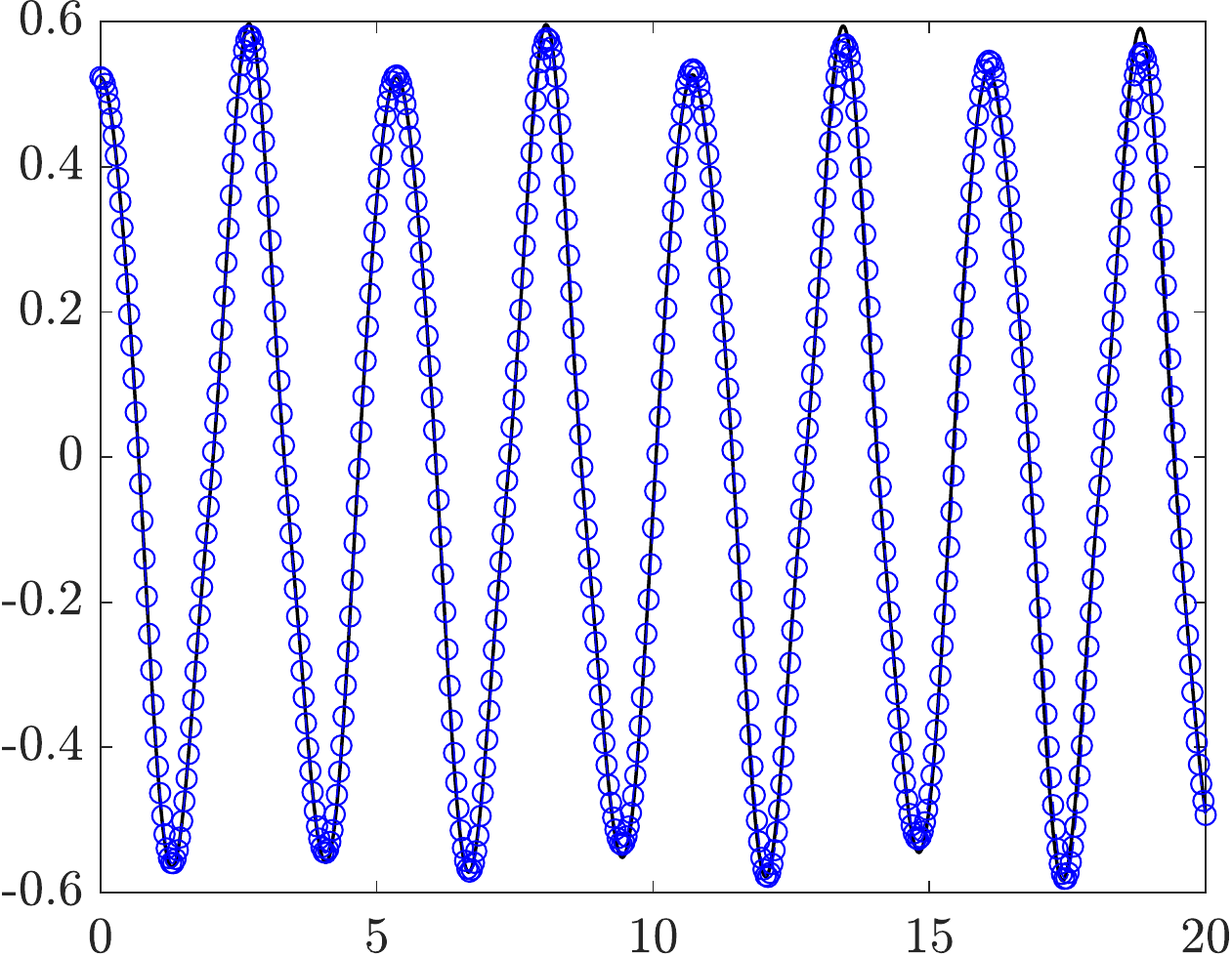}}
	\subfigure[$q_2(t)$]{\includegraphics[width=0.48\textwidth]{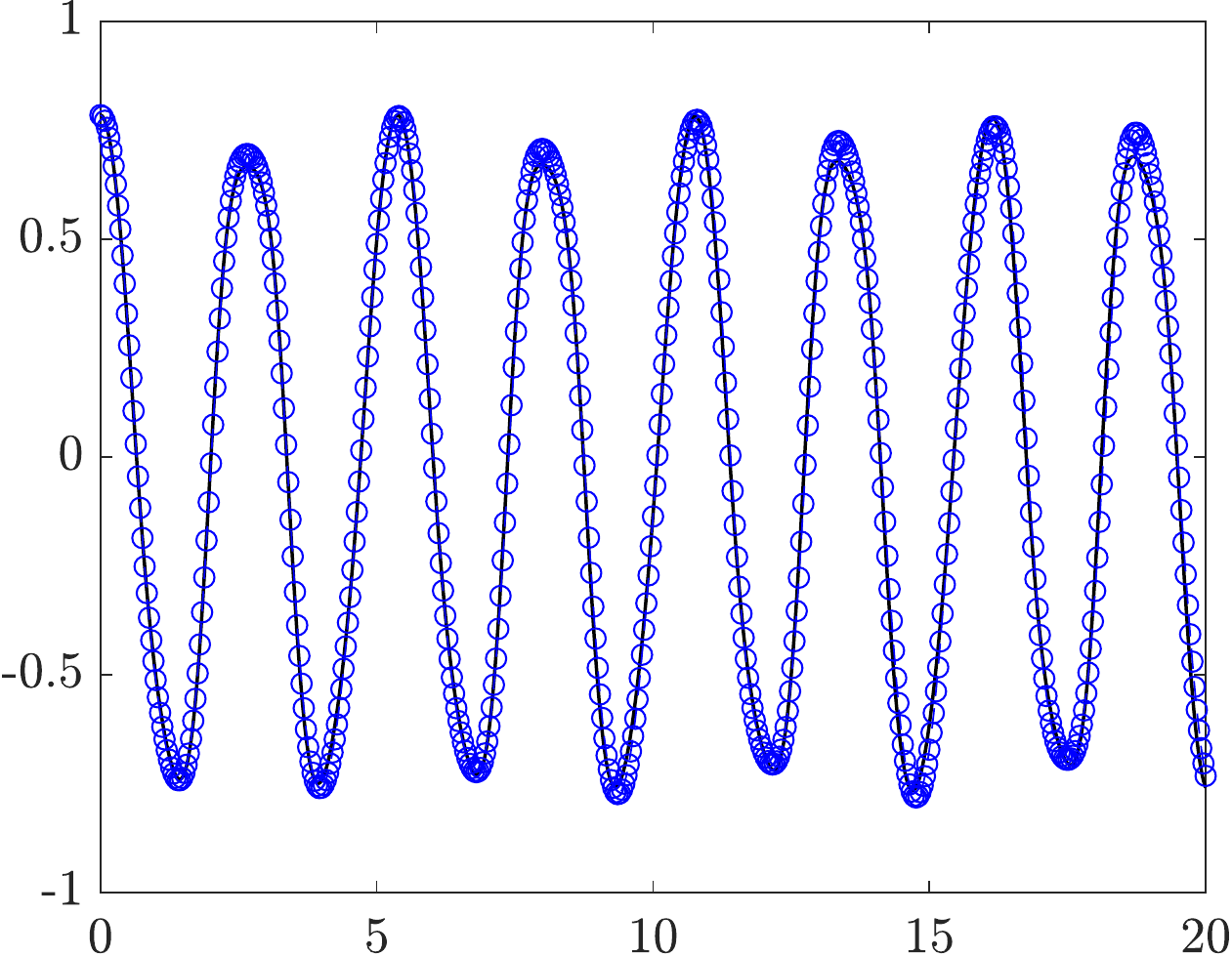}}
	\caption{\small
		Example 5: Comparison of solutions of the
                reconstructed system (circles) with the solution of 
		the true system (solid lines), for the same initial state ${\bf u}_0^*=( 0,0, \frac{\pi}6, \frac{\pi}4 )^\top$. 
	}\label{fig:ex5_U}
\end{figure}

\begin{figure}[htbp]
	\centering
	{\includegraphics[width=0.48\textwidth]{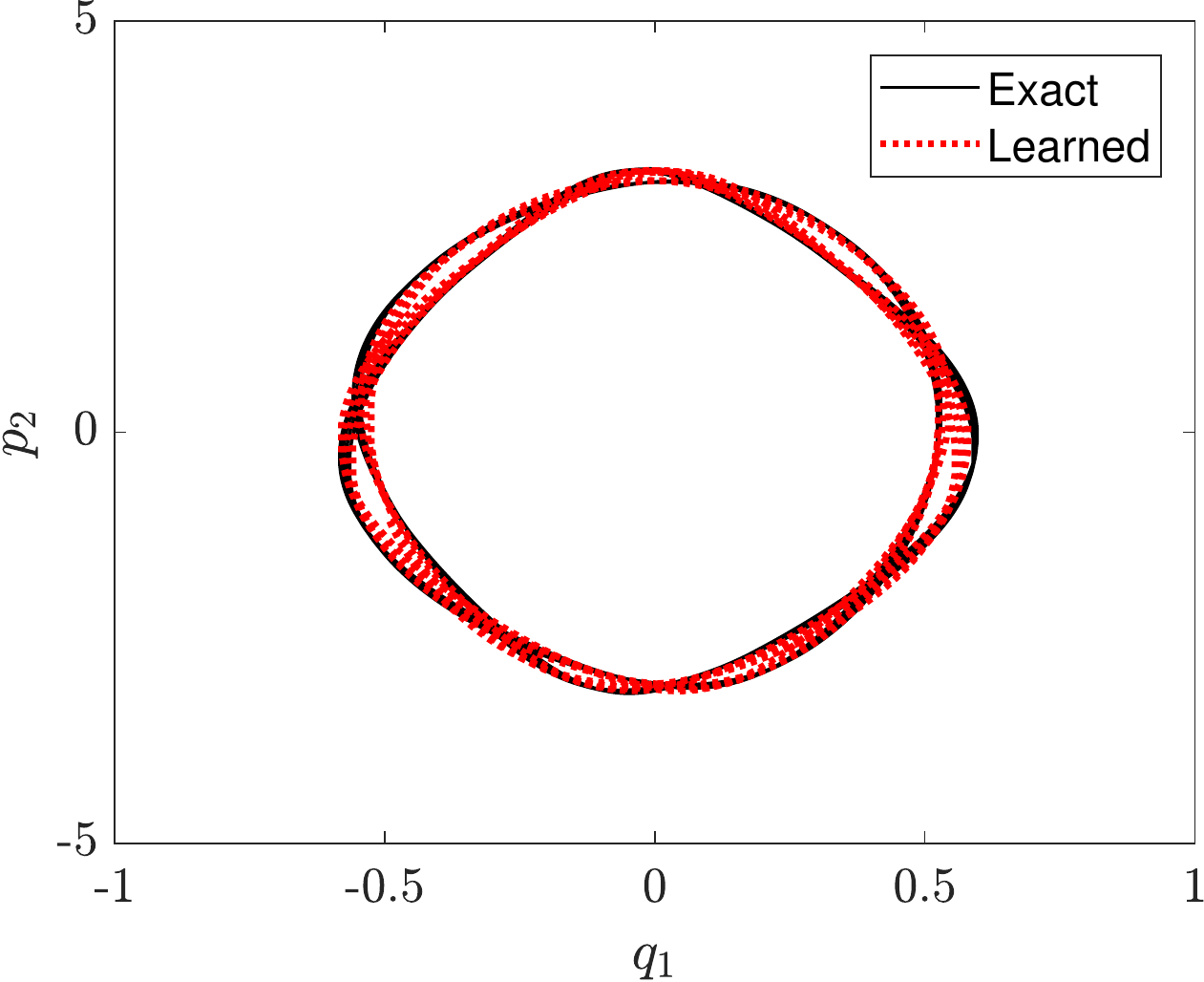}}
	{\includegraphics[width=0.48\textwidth]{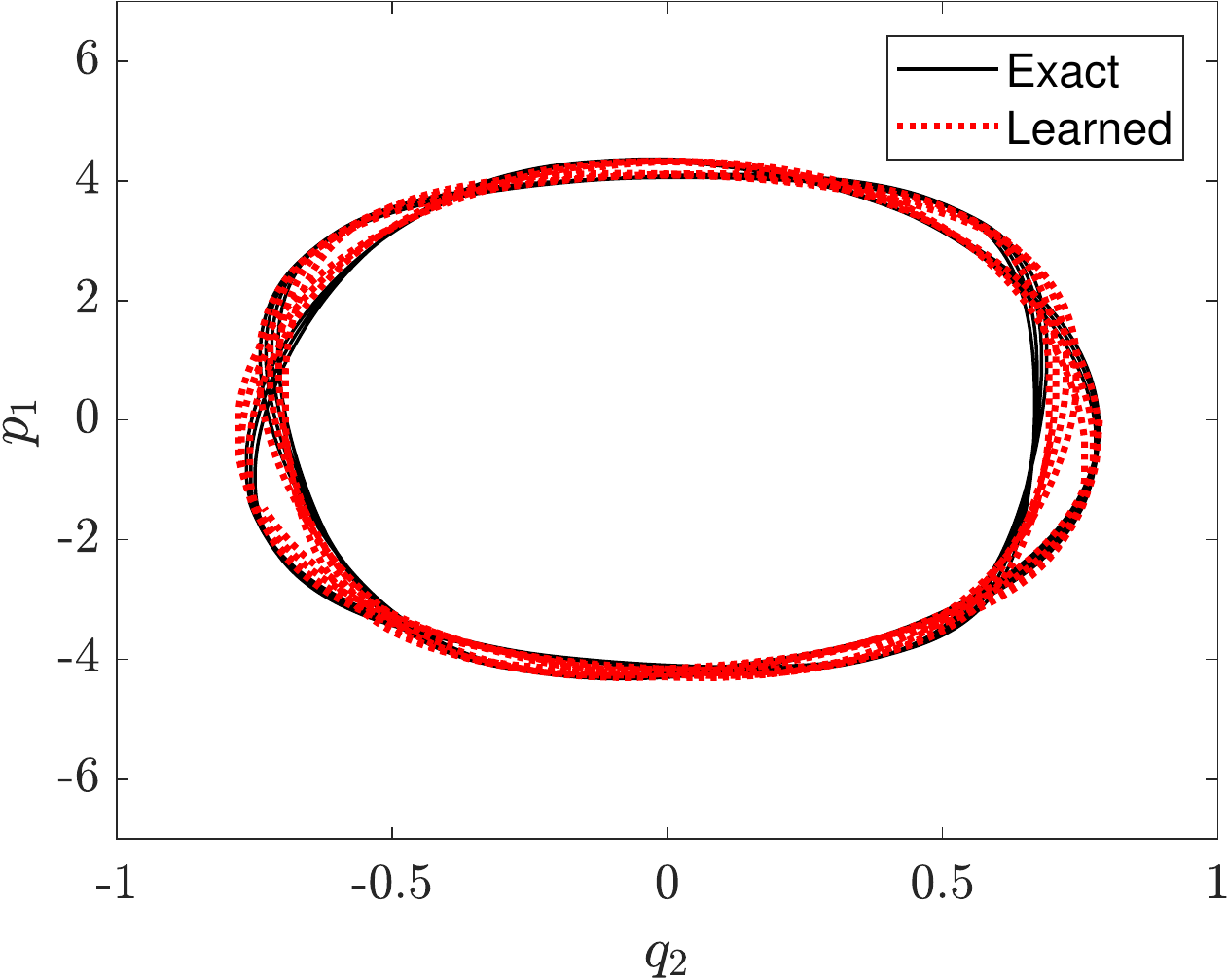}}	
	{\includegraphics[width=0.48\textwidth]{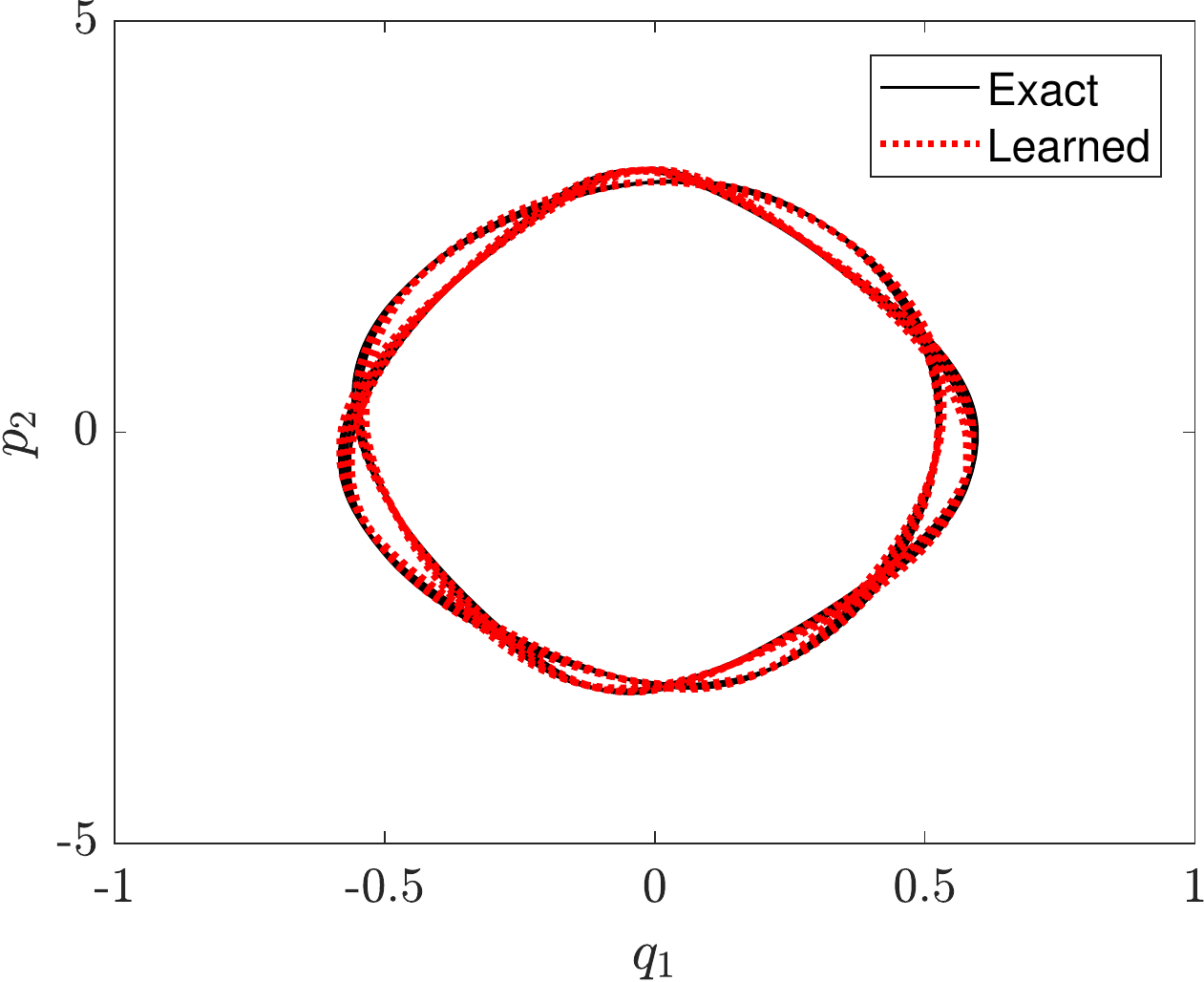}}
	{\includegraphics[width=0.48\textwidth]{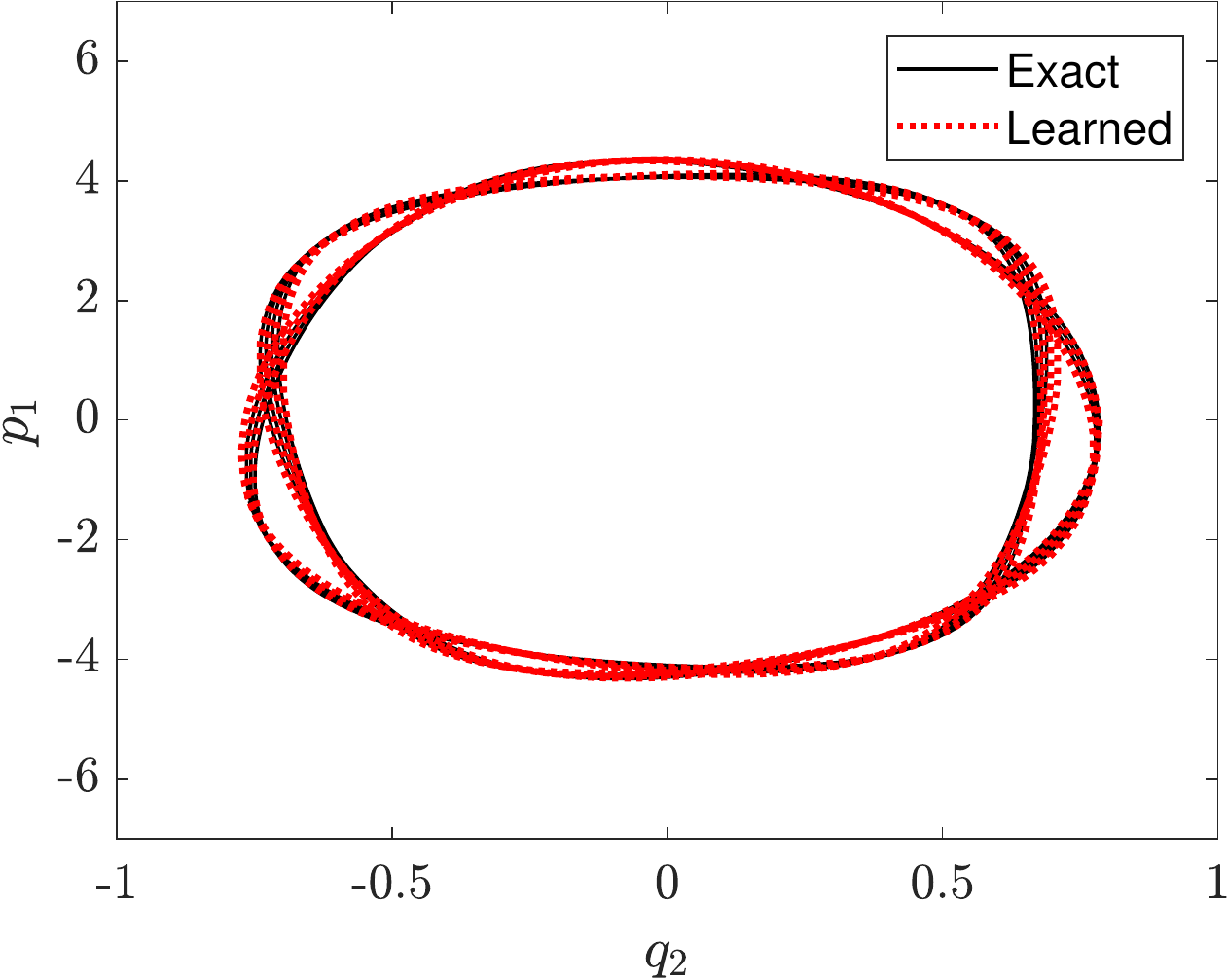}}	
	\subfigure[Plase plots on $q_1$-$p_2$ plane]{\includegraphics[width=0.48\textwidth]{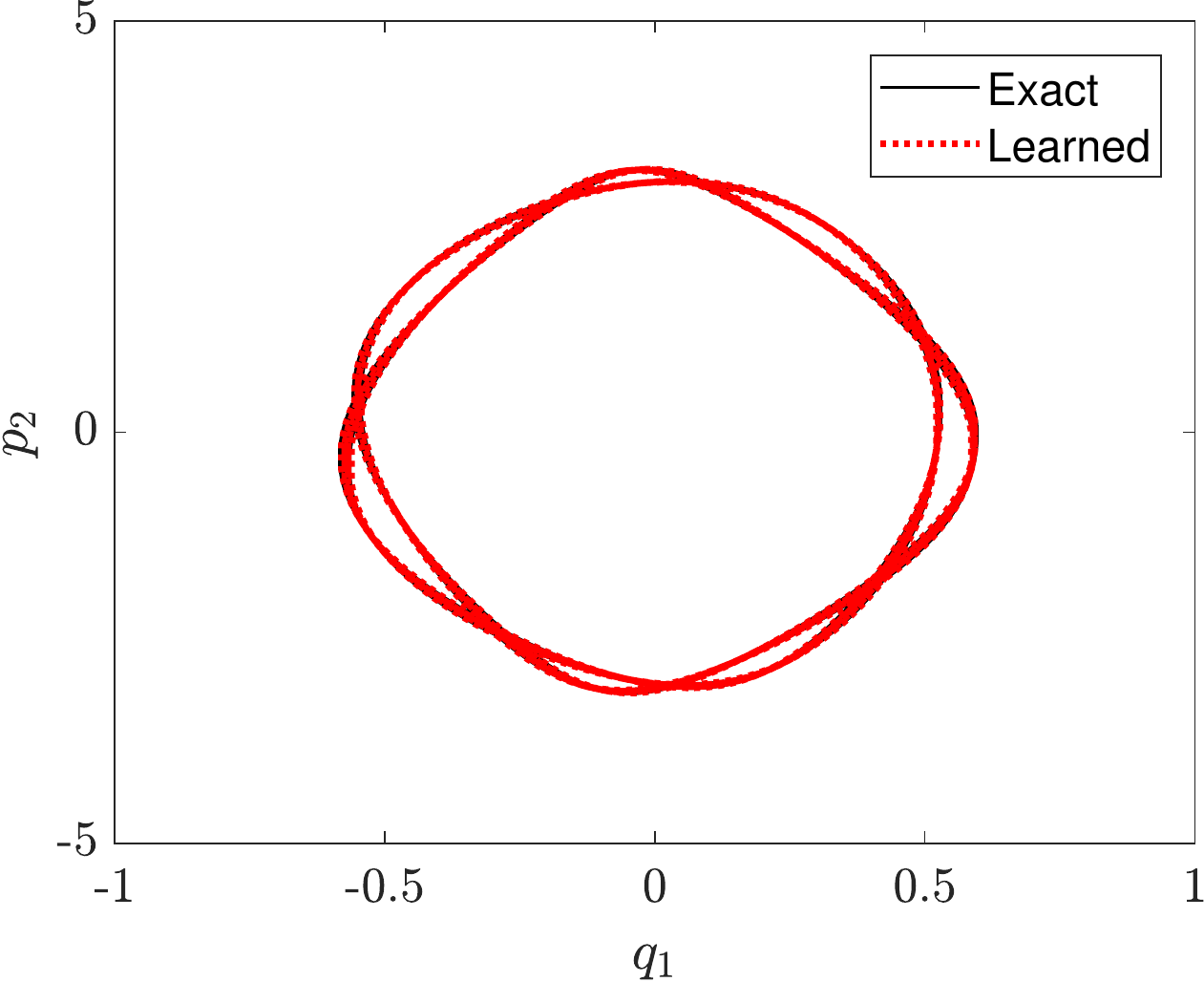}}
	\subfigure[Plase plots on $q_2$-$p_1$ plane]{\includegraphics[width=0.48\textwidth]{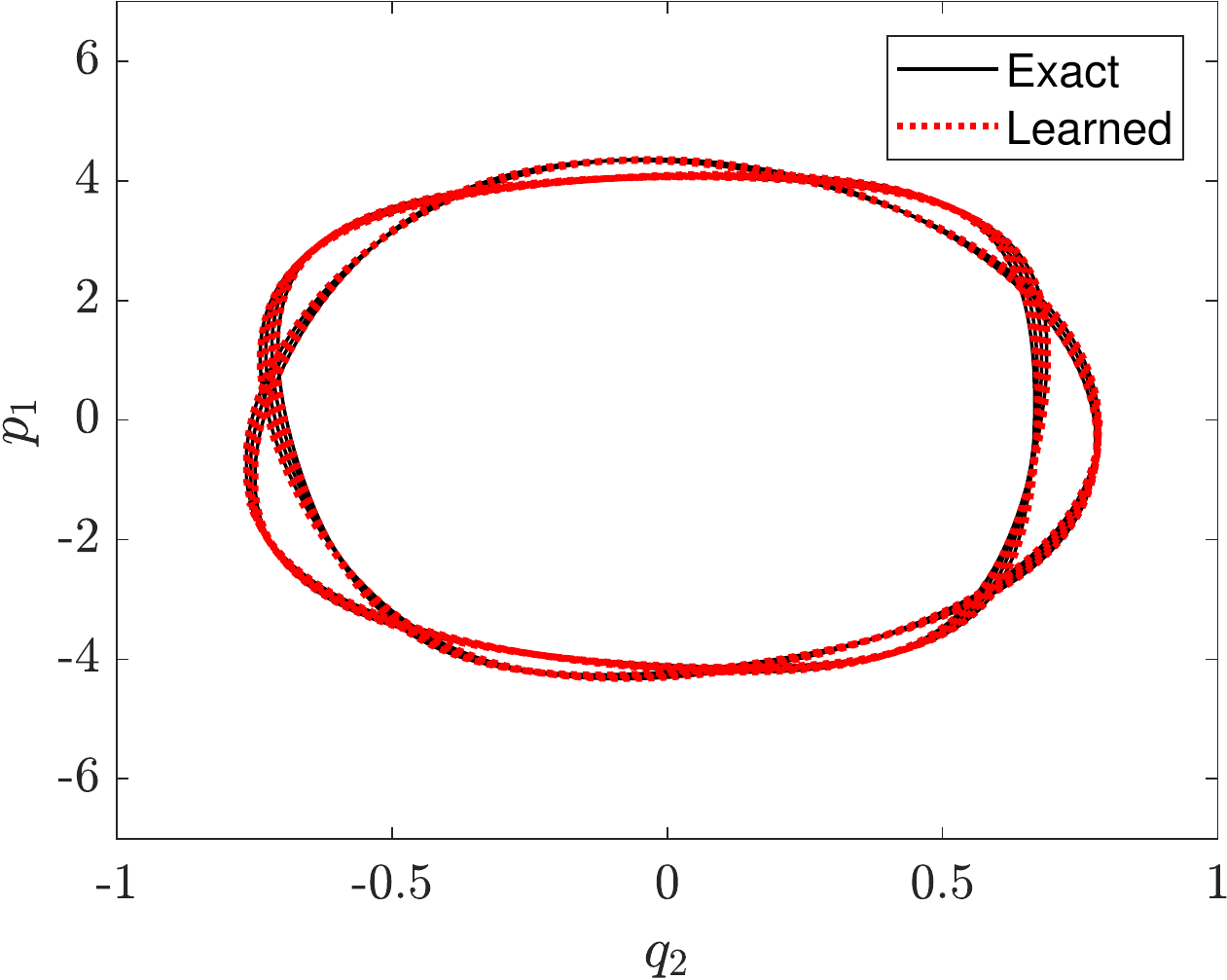}}
	\caption{\small
		Example 5: Phase plots 
		 starting from the initial state ${\bf u}_0^*=( 0,0, \frac{\pi}6, \frac{\pi}4 )^\top$. 
		Top: $n=15, M=20,000$; middle: $n=16, M=20,000$; bottom: $n=18, M=60,000$.
	}\label{fig:ex5_Phase}
\end{figure}

\begin{figure}[htbp]
	\centering
	\includegraphics[width=0.48\textwidth]{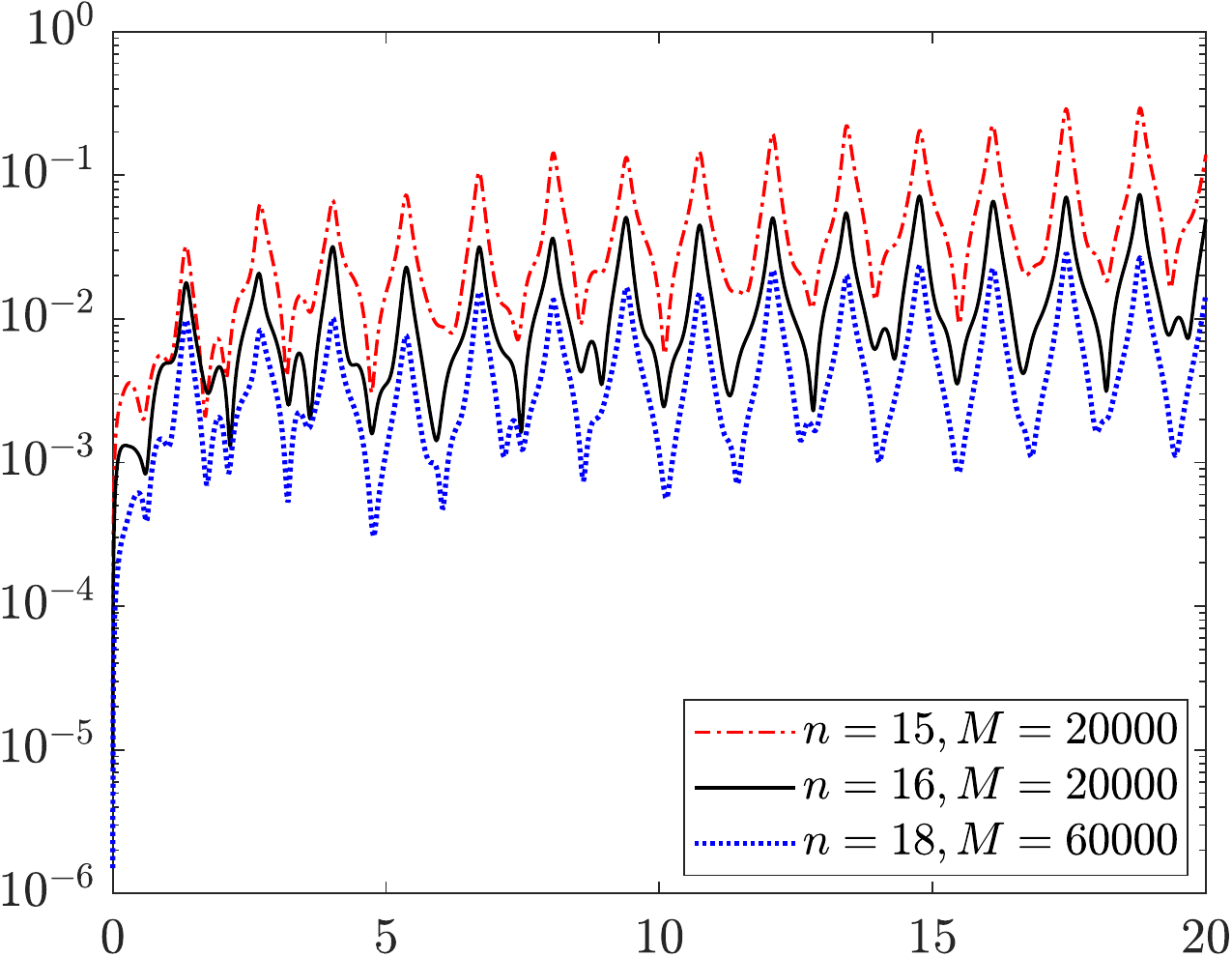}
	\caption{\small
		Example 5: Time evolution of the relative errors in solutions of the reconstructed system with
                an initial state ${\bf u}_0^*=( 0,0, \frac{\pi}6, \frac{\pi}4 )^\top$. 
	}\label{fig:ex5_ErrU}
\end{figure}

In the numerical experiment, we set $m_1 = m_2 =l_1 = l_2=1$ and
$g=9.8$, and set the computational domain as $D= (-5,5)\times (-4,4)
\times (-1,1) \times (-1,1)$.
The Hamiltonian in this example is notably more complicated than the
ones in the previous examples. Consequently, we employ a higher order
polynomial, of degree up to  $n=15$, to conduct the
approximation. The data set include $M=20,000$ short trajectories, each
of which contains $J=2$ intervals. The reconstructed Hamiltonian
system is then solved with an arbitrarily chosen initial state ${\bf
  u}_0^*=( 0,0, \frac{\pi}6, \frac{\pi}4 )^\top$ for up to $T=20$. Its
solution is compared against the reference solutions
from the true system with the same initial state. 
\figref{fig:ex5_HU} shows the evolution of the reconstructed 
Hamiltonian and the true one, which remain constant as expected. 
The time evolution of the solution states are plotted in
\figref{fig:ex5_U}. We observe good agreement with the true
solution.
The corresponding phase plots are further displayed in
\figref{fig:ex5_Phase}, where we also plot the results obtained by 
using higher order
polynomial approximations ($n=16$ with $M=20,000$ and $n=18$ with $M=60,000$) 
to show the convergence. 
We see that the evolution of trajectories is more accurately predicted by 
the reconstructed Hamiltonian system obtained by using higher degree polynomials.  
The relative numerical errors, shown in \figref{fig:ex5_ErrU}, further validate 
the convergence behavior.



%% file: Conclusion.tex
\section{Conclusion} \label{sec:conclusions}

We presented a structure-preserving numerical method for reconstructing unknown Hamiltonian systems using observation data. The key ingredient of the method is to approximate the unknown Hamiltonian first and then derive the approximate equations using the reconstructed Hamiltonian. By doing so, the reconstructed system is able to preserve the approximate Hamiltonian along its trajectories. This is an important property often desired by many practical applications. We presented the algorithm, its error estimate in a special case and used a variety of examples to demonstrate the effectiveness of the approach. In its current form,
polynomials are used to construct the approximation. Other forms of approximation, such as neural networks, will be explored in a separate work.

%% file: Appendix.tex
\appendix

\section{Proof of Theorem \ref{thm:error1}}
\label{app:proofthm2}
The technique used in the following proof is similar to
  the proof of Theorem 3 of \cite{Cohen2013}. However, our Theorem
  \ref{thm:error1} applies to vector-valued function $\nabla H$ in
  gradient function space. This prevents direct use,
  in component-by-component manner, of the result from \cite{Cohen2013}, which applies
  only to scalar-valued function.
  Also, our analysis incorporates numerical errors induced by 
  estimating time derivatives $\dot{\bf x}_k$. These numerical errors
  often do not follow random distribution. Consequently, we do not
  employ i.i.d. assumption on the errors, as opposed to the work of
  \cite{Cohen2013}.
  Due to these subtle, and yet significant, differences, we include
  the proof of Theorem \ref{thm:error1}  here for completeness of the
  paper.

\begin{proof}
	Let $d \omega^K = \otimes^K d \omega$ be the probability
        measure of the random sequence $\{ \x_k \}_{k=1}^K$. Let $\Omega$
        be the set of all possible draws, which is divided into the
        set  
	$\Omega_+$ of all draws such that 
	\begin{equation}\label{keyA12}
	\| {\bf A} - {\bf I} \| \le \frac12,
	\end{equation}
	and the complement set $\Omega_- :=\Omega \setminus \Omega_+$.
	We consider the following splitting
	\begin{equation}\label{eq:proofWKL1}
	\begin{aligned}
	\mathbb E \Big( \big\| \nabla H -   {\bf T}_L  ( \nabla \widetilde H )  \big\|_{2,\mathbb L^2_\omega}^2  \Big) &= \int_{\Omega_+} \big\| \nabla H -   {\bf T}_L  ( \nabla \widetilde H ) \big\|_{2,\mathbb L^2_\omega}^2  d \omega^K 
	\\
	& \quad 
	+ \int_{\Omega_-} \big\| \nabla H -  {\bf T}_L  ( \nabla \widetilde H ) \big\|_{2,\mathbb L^2_\omega}^2  d \omega^K 
	\\
	& =: I_1 + I_2.
	\end{aligned}
	\end{equation}
	We now estimate the upper bounds for $I_1$ and $I_2$.
	
	Let us first consider $I_2$. Based on Corollary
        \ref{coro:Prob} and under the condition \eqref{condition1}, we
        have
	\begin{equation}\label{eq:Prob111}
	{\rm Prob} \{ \Omega_- \} = \int_{\Omega_-} d \omega^K \le 2 K^{-r}.
	\end{equation}
	Note that
	\begin{align*}
	\big\| \nabla H - {\bf T}_L \big( \nabla \widetilde H \big) \big\|_{2,\mathbb L^2_\omega}^2
	&= \int_D \big\| \nabla H (\x)- {\bf T}_L \big( \nabla \widetilde H (\x) \big) \big\|_{2}^2 d \omega 
	\\
	&\le \int_D 2 \Big( \big\| \nabla H (\x) \big\|_2^2 + \big\| {\bf T}_L \big( \nabla \widetilde H (\x) \big) \big\|_{2}^2 \Big) d \omega 
	\\
	& \le 2(L^2+L^2) \int_D d \omega  = 4L^2.
	\end{align*}
	Therefore, we have 
	\begin{equation}\label{eq:I2}
	I_2 = \int_{\Omega_-} \big\| \nabla H - {\bf T}_L \big( \nabla \widetilde H \big) \big\|_{2,\mathbb L^2_\omega}^2  d \omega^K \le 4L^2 {\rm Prob} \{ \Omega_- \}  \le 8L^2 K^{-r}.
	\end{equation}

	We now consider $I_1$. 
	For every $\x \in D$, if $\| \nabla \widetilde H (\x) \|_2 \le L$, then 
	$$
	{\bf T}_L \big( \nabla \widetilde H (\x) \big)= \nabla \widetilde H (\x),
	$$
	so that 
	$$
	\big\| \nabla H (\x) - {\bf T}_L \big( \nabla \widetilde H (\x) \big) \big\|_{2}^2 =
	\big\| \nabla H (\x) - \nabla \widetilde H (\x)  \big\|_{2}^2.
	$$
	For almost every $\x \in D$ with respect to $\omega(\x)$, if $\| \nabla \widetilde H (\x) \|_2 > L$, then 
	$$
	{\bf T}_L \big( \nabla \widetilde H (\x) \big) = \frac{L}{ \| \nabla \widetilde H (\x) \|_2} \nabla \widetilde H (\x),
	$$	
	which implies
	\begin{align*}
	&\big\| \nabla H (\x) - {\bf T}_L \big( \nabla \widetilde H (\x) \big) \big\|_{2}^2 - 
	\big\| \nabla H (\x) - \nabla \widetilde H (\x)  \big\|_{2}^2
	\\
	& = \| {\bf T}_L \big( \nabla \widetilde H (\x) \big)  \|_2^2 -  \| \nabla \widetilde H (\x)  \|_2^2
	+ 2 \nabla H (\x)  \cdot \Big( \nabla \widetilde H (\x) - {\bf T}_L \big( \nabla \widetilde H (\x) \big)   \Big)
	\\
	& \le  L^2 -  \| \nabla \widetilde H (\x)  \|_2^2
	+ 2 L \Big\| \nabla \widetilde H (\x) - {\bf T}_L \big( \nabla \widetilde H (\x) \big)   \Big\|_2
	\\
	& = -\big( \| \nabla \widetilde H (\x) \|_2 - L  \big)^2 < 0.
	\end{align*}
	Therefore, we have 
	\begin{equation*}
	\big\| \nabla H (\x) - {\bf T}_L \big( \nabla \widetilde H (\x) \big) \big\|_{2}^2 \le
	\big\| \nabla H (\x) - \nabla \widetilde H (\x)  \big\|_{2}^2,
	\end{equation*}
	for almost every $\x \in D$ with respect to $\omega(\x)$.
	It follows that 
	\begin{equation}\label{eq:proofWKL3}
	\begin{aligned}
	I_1 &= \int_{\Omega_+} \big\| \nabla H - {\bf T}_L \big( \nabla \widetilde H \big) \big\|_{2,\mathbb L^2_\omega}^2  d \omega^K
	\\
	&\le \int_{\Omega_+} \big\| \nabla H -  \nabla \widetilde H  \big\|_{2,\mathbb L^2_\omega}^2  d \omega^K.
	\end{aligned}
	\end{equation}
	Let us rewrite the derivative data as 
	$$\dot{\bf x}_k =
	{\bf J}^{-1} \nabla H ( {\bf x}_k ) + {\bm \tau}_{k},
	$$	
	where ${\bm \tau}_{k}$ denotes the error in the estimated derivative. 
	Define $ \nabla G:= \nabla H-\Pi_{\mathbb V} (\nabla H)$. 
	Similar to \cite{Cohen2013}, one can write
	$$
	\nabla H- \nabla \widetilde H = \nabla G - \Pi_{\mathbb V}^K (\nabla G) - \nabla \widehat H,
	$$
	where 
	$$
	\Pi_{\mathbb V}^K (\nabla G) := \mathop{\rm argmin}\limits_{ \nabla h \in \mathbb V } \sum_{k=1}^K 
	\left\| \nabla G ( {\bf x}_k  )   - \nabla h ( {\bf x}_k  )   \right\|^2_2,
	$$
	and 
	$$
	\nabla \widehat  H := \mathop{\rm argmin}\limits_{ \nabla h \in \mathbb V } \sum_{k=1}^K 
	\left\| {\bf J} {\bm \tau}_{k}    - \nabla h ( {\bf x}_k  )   \right\|^2_2.
	$$
	Then, we have 
	\begin{align*}
	\big\| \nabla H - \nabla \widetilde H  \big\|_{2,\mathbb L^2_\omega}^2  &= 
	\big\| \nabla G  \big\|_{2,\mathbb L^2_\omega}^2  + \big\|  \Pi_{\mathbb V}^K (\nabla G) + \nabla \widehat H  \big\|_{2,\mathbb L^2_\omega}^2
	\\
	& \le \big\| \nabla G  \big\|_{2,\mathbb L^2_\omega}^2  + 2 \big\|  \Pi_{\mathbb V}^K (\nabla G)   \big\|_{2,\mathbb L^2_\omega}^2 + 2 \big\|   \nabla \widehat H  \big\|_{2,\mathbb L^2_\omega}^2
	\\
	& = \big\| \nabla G  \big\|_{2,\mathbb L^2_\omega}^2  + 2 \sum_{j =1}^{N} \xi_j^2 
	+ 2 \sum_{j =1}^{N} \eta_j^2,
	\end{align*}
	where $\bm \xi =( \xi_1, \dots, \xi_N )^\top$ and $\bm \eta =( \eta_1, \cdots, \eta_N )^\top$ are respectively the solutions of the two systems 
	$$
	{\bf A} \bm \xi = {\bf y}, \qquad \qquad {\bf A} \bm \eta = {\bf z},
	$$
	with the matrix $\bf A$ defined in \eqref{eq:DefA}, ${\bf y}=(y_1,\dots, y_{N})^\top$, ${\bf z}=(z_1,\dots, z_{N})^\top$, and
	\begin{align*}
	&y_i = \frac{1}{K} \sum_{k=1}^K  \Big( \nabla G( {\x}_k ) 
	\cdot \nabla \phi_{i} \big( {\x}_k  \big) \Big),
	\\
	&z_i = \frac{1}{K} \sum_{k=1}^K  \Big( \left( {\bf J} \bm \tau_k \right)
	\cdot \nabla \phi_{i} \big( {\x}_k  \big) \Big), \quad 1\le i \le N.
	\end{align*}
	When the draw $\{ \x_k \}_{k=1}^K$ belong to $\Omega_+$, we have \eqref{keyA12}, which yields
	$\| {\bf A}^{-1} \| \le 2$ and 
	\begin{align*}
	\sum_{j =1}^{N} \xi_j^2 \le 4 \sum_{j=1}^{N} y_i^2, \qquad \quad 
	\sum_{j =1}^{N} \eta_j^2 \le 4 \sum_{j=1}^{N} z_i^2. 
	\end{align*}	
	Hence 
	\begin{equation}\label{eq:proofWKL4}
	\begin{aligned}
	I_1 &\le \int_{\Omega_+} \left(
	\big\| \nabla G  \big\|_{2,\mathbb L^2_\omega}^2 + 8 \sum_{j =1}^N y_j^2 
	+ 8 \sum_{j =1}^N z_j^2
	\right) d \omega^K
	\\
	& \le  \big\| \nabla G  \big\|_{2,\mathbb L^2_\omega}^2
	+ 8 \sum_{j =1}^{N} \mathbb E \big( y_j^2  \big)
	+ 8 \sum_{j =1}^{N} \mathbb E \big( z_j^2  \big).
	\end{aligned}
	\end{equation}
	For each $1\le j \le N$, we estimate $\mathbb E \big( y_j^2  \big)$ as follows:
	\begin{align*}
	\mathbb E \big( y_j^2  \big) & = \frac{1}{K^2} \sum_{k=1}^K \sum_{l=1}^K 
	\mathbb E \left[ \left( \nabla G ( \x_k) \cdot \nabla \phi_{j} ( \x_k ) 
	\right) \left( \nabla G ( \x_l) \cdot \nabla \phi_{j} ( \x_l ) 
	\right) \right]	
	\\
	& = \frac{1}{K^2} \left[
	K(K-1)  \big| \mathbb E  \big( \nabla G ( \x) \cdot \nabla \phi_{j} ( \x ) \big) \big|^2  
	+ K   \mathbb E  \Big( \big| \nabla G ( \x) \cdot \nabla \phi_{j} ( \x )  \big|^2   \Big)
	\right] 
	\\
	& = \left(1- \frac{1}{K} \right) \left| \int_D \nabla G ( \x) \cdot \nabla \phi_{j} ( \x ) d \omega  \right|^2 + 
	\frac{1}{K} \int_{D}  \big| \nabla G ( \x) \cdot \nabla \phi_{j} ( \x )  \big|^2 d \omega
	\\
	& = \frac{1}{K} \int_{D}  \big| \nabla G ( \x) \cdot \nabla \phi_{j} ( \x )  \big|^2 d \omega 
	\le \frac{1}{K} \int_{D}  \| \nabla G ( \x) \|_2^2 \| \nabla \phi_{j} ( \x )  \|_2^2 d \omega,
	\end{align*}
	where the Cauchy–Schwarz inequality has been used in the inequality. 
	It follows that 
	\begin{align*}
	\sum_{j =1}^{N} \mathbb E \big( y_j^2  \big)  & \le \frac{1}{K} \int_{D}  \| \nabla G ( \x) \|_2^2 
	\left( \sum_{j=1}^{N} \| \nabla \phi_{j} ( \x )  \|_2^2 \right) d \omega
	\\
	& \le \frac{{\mathscr K}_N}{K} \int_{D}  \| \nabla G ( \x) \|_2^2  d \omega  \le \frac{\lambda}{\log K} \big\| \nabla G  \big\|_{2,\mathbb L^2_\omega}^2.
	\end{align*}
	We now estimate $\mathbb E \big( z_j^2  \big)$ for each $1\le j \le N$.
	\begin{align*}
	\mathbb E \big( z_j^2  \big) &= \frac{1}{K^2} 
	\mathbb E \left( \sum_{k=1}^K  \Big( \left( {\bf J}\bm \tau_k  \right)
	\cdot \nabla \phi_{j} \big( {\bf x}_k  \big) \Big) \right)^2
	\\
	& \le \frac{1}{K^2}   \mathbb E \left( \sum_{k=1}^K  \Big\| {\bf J} \bm \tau_k \Big\|_2^2  \sum_{k=1}^K 
	\Big\| \nabla \phi_{j} \big( {\bf x}_k  \big) \Big\|_2^2 \right) 
	\\
	& \le \frac{1}{K^2}   \mathbb E \left( \sum_{k=1}^K  \tau_\infty^2  \sum_{k=1}^K
	\Big\| \nabla \phi_{j} \big( {\bf x}_k  \big) \Big\|_2^2 \right)  
	\\
	&= \tau_\infty^2  \big \| \nabla \phi_{j}  \big \|_{2,\mathbb L^2_\omega}^2 = \tau_\infty^2.
	\end{align*}	
	It follows from \eqref{eq:proofWKL4} that 
	$$
	I_1 \le \left( 1 + \frac{8\lambda}{\log K} \right) \big\| \nabla G  \big\|_{2,\mathbb L^2_\omega}^2 + 8N \tau_\infty^2, 
	$$
	which together with \eqref{eq:I2} complete the proof.
\end{proof}